\documentclass{amsart}
\usepackage{amssymb}

\usepackage{amscd}
\usepackage{thmdefs}

\input{tcilatex}
\theoremstyle{definition}
\theoremstyle{remark}
\numberwithin{equation}{section}

\input tcilatex

\begin{document}
\title[$C^{*}$-Algebras Generated by Partial Isometries]{$C^{*}$-Algebras Generated by Partial Isometries}
\author{Ilwoo Cho and Palle Jorgensen}
\address{Saint Ambrose Univ., Dep. of Math, 421 Ambrose Hall, 518 W. Locust Street,
Davenport, Iowa, 52308, U. S. A.\\
\strut \strut \strut \\
University of Iowa, Dep, of Math, McLean Hall, Iowa City, Iowa, U. S. A.}
\email{chowoo@sau.edu\\
\strut \\
jorgensen@math.uiowa.edu}
\thanks{}
\date{Dec., 2007}
\subjclass{05C62, 05C90, 17A50, 18B40, 46K10, 47A67, 47A99, 47B99}
\keywords{Finite Directed Graphs, Partial Isometries, Unitary parts, Shift Parts,
Graph Groupoids, $C^{*}$-Algebras, $*$-Isomorphic Index, Algebraic Free
Product, Algebraic Reduced Free Product, Topological Free Product,
Topological Reduced Free Product.}
\dedicatory{}
\thanks{}
\maketitle

\begin{abstract}
We prove a structure theorem for a finite set $\mathcal{G}$ of partial
isometries in a fixed countably infinite dimensional complex Hilbert space $%
H.$ Our result is stated in terms of the $C^{*}$-algebra generated by $%
\mathcal{G}.$ The result is new even in the case of a single partial
isometry which is not an isometry or a co-isometry; and in this case, it
extends the Wold decomposition for isometries. We give applications to
groupoid $C^{*}$-algebras generated by graph groupoids, and to partial
isometries which have finite defect indices and which parametrize the
extensions of a fixed Hermitian symetric operator with dense domain on the
Hilbert space $H.$ Our classification parameters for the ``Wold decomposed''
set $\mathcal{G}_{W}$ of our finite set $\mathcal{G}$ of partial isometries
involve infinite and explicit Cartesian product sets, and they are
computationally attractive. Moreover, our classification labels generalize
the notion of defect indices in the special case of the family $\mathcal{G}%
_{W}$ of partial isometries from the Cayley transform theory and Hermitian
extensions of unbounded Hermitian operators with dense domain.
\end{abstract}

\strut

\textbf{Table of Contents}

\strut

$
\begin{array}{lll}
\text{{\small 1. Introduction}} &  & \;2 \\ 
\text{{\small 2. Preliminaries}} &  &  \\ 
\qquad \text{{\small 2.1. Partial Isometries on a Hilbert Space}} &  & \;5
\\ 
\qquad \text{{\small 2.2. Directed Graphs and Graph Groupoids}} &  & \;5 \\ 
\qquad \text{{\small 2.3.}}\ C^{*}\text{{\small -Algebras}}\ \text{{\small %
Induced}}\ \text{{\small by}}\ \text{{\small Graphs}} &  & \;9 \\ 
\qquad \text{{\small 2.4.}}\ \text{{\small Groupoids}}\ \text{{\small and}}\ 
\text{{\small Groupoid Actions}} &  & 12 \\ 
\text{{\small 3.}}\ C^{*}\text{{\small -Subalgebras of}}\ B(H)\ \text{%
{\small Generated by Partial Isometries}} &  &  \\ 
\qquad \text{{\small 3.1.}}\ *\text{{\small -Isomorphic Indices}}\ \text{%
{\small of Partial Isometries}} &  & 13 \\ 
\qquad \text{{\small 3.2. Graph Groupoids Induced}}\ \text{{\small by
Partial Isometries}} &  & 22 \\ 
\qquad \qquad \text{{\small 3.2.1. Construction of Corresponding Graphs of
Partial Isometrie}}s &  & 25 \\ 
\qquad \qquad \text{{\small 3.2.2. The Glued Graph}}\
G_{1}^{v_{1}}\#^{v_{2}}G_{2}\ \text{{\small of}}\ G_{1}\ \text{{\small and}}%
\ G_{2} &  & 36 \\ 
\qquad \qquad \text{{\small 3.2.3. Conditional Itrerated Gluing on}}\ 
\mathcal{G}_{G} &  & 27 \\ 
\qquad \qquad \text{{\small 3.2.4. The Vertex of the}}\ \mathcal{G}\text{%
{\small -Grap}}h &  & 33 \\ 
\qquad \text{{\small 3.3. A Representation of the}}\ \mathcal{G}\text{%
{\small -Groupoi}}d &  & 34 \\ 
\qquad \text{{\small 3.5.}}\ C^{*}\text{{\small -Subalgebras Generated by
Partial Isometrie}}s &  & 36
\end{array}
$

{\small \qquad \qquad \qquad \qquad }

$
\begin{array}{lll}
\text{{\small 4. Block Structures of}}\ C^{*}(\mathcal{G}) & \qquad \qquad
\qquad \qquad &  \\ 
\qquad \text{{\small 4.1. Topological Reduced Free Product}}\ \text{{\small %
Algebras}} & \qquad \qquad & 38 \\ 
\qquad \text{{\small 4.2. Topological Free Block Structures on}}\ C^{*}(%
\mathcal{G}) &  & 43 \\ 
\qquad \text{{\small 4.3. Example}}s &  & 52 \\ 
\text{{\small 5. Extensions of Unbounded Operator}}s &  & 59
\end{array}
$

\strut

\strut

\strut

\section{Introduction}

\strut

\strut

In this paper, all Hilbert spaces are separable infinite dimensional. We
will consider $C^{*}$-subalgebras of $B(H)$ generated by finitely many
partial isometries on $H$, where $B(H)$ is the operator algebra consisting
of all (bounded linear) operators. First, we characterize the $C^{*}$%
-subalgebras generated by a single partial isometry characterized by a $*$%
-isomorphism. We show that they are characterized by a $*$-isomorphism index
in $(\Bbb{N}_{0}^{\infty })^{4},$ where $\Bbb{N}_{0}^{\infty }$ $=$ $\Bbb{N}$
$\cup $ $\{0,$ $\infty \}$, up to finite numbers. This shows that the $C^{*}$%
-subalgebras of $B(H)$ generated by a single partial isometry are
characterized by quadruples of numbers in $\Bbb{N}_{0}^{\infty }.$ More
generally, we consider $C^{*}$-subalgebras of $B(H)$ generated by finitely
many partial isometries. We prove that they are determined by certain
combinatorial objects, called the conditional iterated glued graphs and
their corresponding graph groupoids. By observing the block structures of
the graph groupoids, we can have the block structures of a $C^{*}$-algebra
generated by finitely many partial isometries in $B(H).$

\strut

Let $H$ be a Hilbert space and let $a$ $\in $ $B(H)$ be an operator, where $%
B(H)$ is the operator algebra consisting of all bounded linear operators on $%
H.$ We say that this operator $a$ is a \emph{partial isometry} if the
operator $a^{*}$ $a$ is a projection on $H,$ where $a^{*}$ is the adjoint of 
$a.$ Recall that an operator $p$ $\in $ $B(H)$ is a \emph{projection} if $p$
satisfies that $p^{2}$ $=$ $p$ $=$ $p^{*}.$ We study $C^{*}$-subalgebras
generated by finitely many partial isometries in $B(H).$ We will extend the
classification results of [13].

\strut

Directed graphs, both finite and infinite, play a role in a myriad of areas
of applications; electrical network of resistors, the internet, and random
walk in probability theory, to mention only a few. Some of the results in
the subject stress combinatorial aspects of the graphs, while others have a
more analytic slant. Applications to quantum theory fall in the latter
group, and that is where Hilbert space and non-commuting operators play a
role. To see this, consider a given graph $G$, and some matrix operation
which transfers data encoded in one site (or vertex) of $G$ onto another. If
there is a notion of energy, or some other conserved quantity for the entire
graph, then there is also a naturally associated Hilbert space, and with the
energy now representing the norm-squared. Under these conditions, the matrix
transfer between vertices will take the form of partial isometries: As a
result, we will then be assigning partial isometries to the directed edges
in $G$. In this paper we will derive some of the $C^{*}$-algebraic
invariants associated with such an approach. Hence, the analytic approach
asks for representations of the directed graphs, but each application
motivates different representations.

\strut

Before turning to our main theorem, we recall some basic definitions and
results in the literature which are natural precursors to our present model.
A partial isometry $a$ in a Hilbert space $H$ is a linear mapping $a$ $:$ $H$
$\rightarrow $ $H,$ which restricts to an isometry between two closed
subspaces in $H$, the initial space $H_{init}^{a}$ and the final space $%
H_{fin}^{a}$. In addition, it is assumed that $a$ maps the orthogonal
complement of $H_{init}^{a}$ into $\{0\}$. The dimensions of the
ortho-complements of the two spaces are called the defect indices of $a$. If
the first defect index is $0$, we say that $a$ is an isometry, and if the
second is 0, a co-isometry. While there are good classification theorems for
isometries, e.g., [27]; the case of partial isometries when both the defect
indices are non-zero have so far resisted classification. See however [28].
From the pioneering work of Akhiezer and Glazman [26], we know that the
Cayley transform of unbounded densely defined skew-symmetric ordinary
differential operators are partial isometries with finite defect indices. In
Section 3.1, we generalize such indices.

\strut

The simplest and most primitive case, of course would be asking for a
universal $C^{*}$-algebra which is generated by a single partial isometry.
Using a certain grading on this algebra, the author of [33] constructs an
associated analog of the Cuntz algebras [31]. Recall that the Cuntz algebra $%
O_{d}$ is the $C^{*}$-algebra on $d$-isometries on $H$ which have mutually
orthogonal final spaces adding up to $H$.

\strut

A popular invariant is the $KK$-group, and with his universal algebra,
Kandelaki gives a homotopical interpretation of $KK$-groups (See [33] and
[37]). He proves that his universal $C^{*}$-algebra is homotopically
equivalent to $M_{2}(\Bbb{C})$ up to stabilization by $(2$ $\times $ $2)$%
-matrices. Therefore, these Kandelaki algebras are $KK$-isomorphic.

\strut

This paper provides not only a characterization of $C^{*}$-subalgebras of $%
B(H)$ generated by a single partial isometry (See Section 3.1) but also a
topological reduced free block structures for the analysis of $C^{*}$%
-subalgebras of $B(H)$ generated by finitely many partial isometries on $H.$
In particular, if $a$ is a partial isometry, then we have the so-called Wold
decomposition $u$ $+$ $s$ of $a,$ where $u$ is the unitary part of $a$ and $%
s $ is the shift part of $a,$ and hence the $C^{*}$-subalgebra $C^{*}(\{a\})$
of $B(H)$ is the $C^{*}$-algebra $C^{*}(\{u,$ $s\})$, which is $*$%
-isomorphic to $C^{*}(\{u\})$ $\oplus $ $C^{*}(\{s\})$ inside $B(H),$
generated by $u$ and $s.$ We realize that $C^{*}(\{u\})$ and $C^{*}(\{s\})$
are characterized by certain graph groupoids $\Bbb{G}_{u}$ and $\Bbb{G}_{s},$
induced by the corresponding graphs $G_{u}$ and $G_{s}$ of $u$ and $s,$
respectively.

\strut \strut \strut

In conclusion, we show that if $\mathcal{G}$ $=$ $\{a_{1},$ ..., $a_{N}\}$
is a family of partial isometries in $B(H),$ and $\mathcal{G}_{W}$ $=$ $%
\{x_{1},$ ..., $x_{n}\}$, the corresponding Wold decomposed family of $%
\mathcal{G},$ and if $\mathcal{G}_{G}$ $=$ $\{G_{x_{1}},$ ..., $G_{x_{n}}\}$
is the family of corresponding graphs of $\mathcal{G}_{W},$ then the $C^{*}$%
-subalgebra $C^{*}(\mathcal{G})$ generated by $\mathcal{G}$ is $*$%
-isomorphic to the groupoid algebra $\mathcal{A}$ $=$ $C_{\alpha }^{*}(\Bbb{G%
}),$ as embedded $C^{*}$-subalgebras of $B(H).$ Futhermore, the $C^{*}$%
-algebra $\mathcal{A}$ has its building blocks determined by the topological
reduced free product.

\strut

Our motivations come from both pure and applied mathematics, including
mathematical physics.

\strut

On the applied side, we notice that the internet offers graphs of very large
size, hence to an approximation, infinite. This is a context where algebraic
models have been useful (e.g., [55] and the references cited in it).

\strut

Features associated with infinite models are detected especially nicely with
the geometric tools from operators on a Hilbert space. A case in point is
the kind of Transfer Operator Theory or Spectral Theory which goes into the
mathematics of internet search engines (e.g., [55]). A second instance is
the use of graph models in the study of spin models in Quantum Statistical
Mechanics (e.g., [30], [32], [49] and [50]). In both of these classes of
``infinite'' models mentioned above, what happens is there is an essential
distinction between conclusions from the use of infinite dimensional models,
as compared to the finite dimensional counterparts.

\strut

And there is the same distinction between the mathematics of a finite number
of degrees of freedom vs. the statistics of infinite models: For a finite
number of degrees of freedom, we have the Stone-von Neumann uniqueness
theorem, but not so in the infinite case. Nonetheless, isomorphism classes
of $C^{*}$-algebras are essential for the study of fields and statistics in
Physics. In fact, these math physics models lie at the foundations of
Operator Algebra Theory.

\strut

In Statistical Mechanics, we know that phase-transition may happen in
certain infinite models, but of course is excluded for the corresponding
finite models.

\strut

On the pure side, there are several as well: First, in the literature, there
is a variety of classification theorems for $C^{*}$-algebras generated by
systems of isometries, starting with the case of a single isometry; Coburn's
Theorem (See [54]). Second, in Free Probability (e.g., [16] and [18]), $%
C^{*} $-algebras built on generators and relations are the context for
``con-comutative'' random variables.

\strut \strut

\strut

\strut

\section{Preliminaries}

\strut

\strut

In this paper, we will consider a $C^{*}$-algebra $C^{*}(\mathcal{G})$
generated by a family $\mathcal{G}$ of finitely many partial isometries in $%
B(H).$ We first introduce the main objects we will need throughout the paper.

\strut

\strut

\strut

\subsection{Partial Isometries on a Hilbert Space}

\strut

\strut

We say that an operator $a$ $\in $ $B(H)$ is a \emph{partial isometry}, if
the operators $a^{*}a$ $\in $ $B(H)$ is a projection. The various equivalent
characterizations of partial isometries are well-known: the operator $a$ is
a partial isometry if and only if $a$ $=$ $aa^{*}a,$ if and only if its
adjoint $a^{*}$ is a partial isometry, too, in $B(H)$. Recall that an
operator $p$ in $B(H)$ is a projection if $p$ is a self-adjoint idempotent.
i.e., $p^{2}$ $=$ $p$ $=$ $p^{*}$ in $B(H).$ Every partial isometry $a$ has
its initial space $H_{init}^{a}$ and its final space $H_{fin}^{a}$ which are
closed subspaces of $H.$ Notice that a partial isometry $a$ is a unitary
from $H_{init}^{a}$ onto $H_{fin}^{a}.$ i.e., $a$ $\in $ $B(H_{init}^{a},$ $%
H_{fin}^{a})$ is a unitary satisfying that $a^{*}a$ $=$ $1_{H_{init}^{a}}$
and $aa^{*}$ $=$ $1_{H_{fin}^{a}},$ where $1_{K}$ means the identity
operator on a Hilbert space $K.$ Therefore, as Hilbert spaces, $H_{init}^{a}$
and $H_{fin}^{a}$ are Hilbert-space isomorphic.

\strut

\strut

\subsection{Directed Graphs and Graph Groupoids}

\strut

\strut

Countable directed graphs and their applications have been studied
extensively in Pure and Applied Mathematics (e.g., [4], [5], [6], [9], [10],
[17], [19], [24], [35] and [36]). Not only are they connected with certain
noncommutative structures but they also let us visualize such structures.
Futhermore, the visualization has a nice matricial expressions.

\strut \strut \strut \strut \strut

A \emph{graph} is a set of objects called vertices (or points or nodes)
connected by links called edges (or lines). In a \emph{directed graph}, the
two directions are counted as being distinct directed edges (or arcs). A
graph is depicted in a diagrammatic form as a set of dots (for vertices),
jointed by curves (for edges). Similarly, a directed graph is depicted in a
diagrammatic form as a set of dots jointed by arrowed curves, where the
arrows point the direction of the directed edges.

\strut

Our global results and theorems have applications to Operator Theory as well
as to a class of finite and infinite electrical networks used in Analysis on
fractals (e.g., [53]), and Physics (e.g., [49] and [50]). We have given our
results an axiometic formulation, making use here of a general Operator
Theoretic framework. The operator-theoretic approach is new, especially in
its general scope. When stated in the context of Hilbert-space operators, we
believe that our results are also of independent interest; quite apart from
their other applications. Our results have a variety of applications, e.g.,
to fractals with affine selfsimilarity (See [53]), to models in Quantum
Statistical Mechanics (See [51]), to the analysis of energy forms,
interplaying between selfsimilar measures and associated energy forms; and
to random walks on graphs (See [53]).

\strut

By a graph $G,$ we will mean a system of vertices $V(G)$ and edges $E(G)$
(See the previous axioms below). Intuitively, the points in $E(G)$ will be
``lines'' connecting prescribed pairs of points in $V(G).$ Each vertex $x$
will have a finite set $g(x)$ of edges connecting it to other points in $%
V(G).$ (For models in Physics, each set $g(x)$ will correspond to a set of
nearest neighbors for the vertex $x,$ with the vertices typically arranged
in a lattice (See [49], [50] and [51]).)

\strut

But our present graphs will be more general. Starting with $G,$ we specify a
function $g$ from $V(G)$ into the set of all finite subsets of $V(G),$
subject to certain axioms; see below for details.

\strut

An important difference between there graphs and those that are more
traditional in Geometric Analysis is that we do not ``a priori'' fix an
orientation. Nonetheless, for electrical network models, we will be using
orientations, but they will be introduced only indirectly, via a prescribed
function $R,$ positive, real valued and defined on $E(G),$ i.e., on the set
of edges on $G.$ (The letter $R$ is for resistance.) By analogy to models in
Electrical Engineering, we think of such a function $R$ as a system of
resistors, and we will be interested in the pair $(G,$ $R).$ Given a system $%
R$ of resistors, then potential theoretic problems for $(G,$ $R)$ will
induce currents in the global system, the associated current functions will
be denoted by $I,$ and each admissible (subject to Kirchhoff's laws) current
function induces an orientation.

\strut

This is a second class of functions $I$ ($I,$ for electrical current) again
functions from $E(G)$ into the real numbers. It is only when such a function 
$I$ is prescribed, that our graph $G$ will acquire an orientation: Specially
we say that an edge $e$ is positively oriented (relative to $I$), if $I(e)$ $%
>$ $0.$ And a different choice of current function $I$ will induce a
different orientation on the graph $G.$

\strut

Our results will be stated in terms of precise mathematical axioms, and we
present our theorems in the framework of Operator Theory / Operator Algebra,
which happens to have an element of motivation from the theory of infinite
systems of resistors. The main tools in our proofs will be operators on a
Hilbert space. In particular, we will solve problems in Discrete Potential
Theory with the use of a family of operators and their adjoints. Good
references to the fundamentals for operators are [44] and [48]. The second
author of this paper recalls conversations with Raul Bott with a remainder
about an analogous Hilbert-space Operator Theoretic approach to Electrical
Networks; apparently attempted in the 1950's in the engineering literature,
but we did not find details in journals. The closest we could come is the
fascinating paper [45], by Bott et al.

\strut

A motivation for our present study is a series of papers in the 1970's by
Powers which introduced infinite systems of resistors into the resolution of
an important question from Quantum Statistical Mechanics (e.g., [49], [50],
[51] and [52]). This is coupled with the emergence of these same techniques,
or their close vintage. By ``Analysis on Fractals'' we refer to a set of
Potential Theoretic tools developed recently by Kigami, Strichartz and
others (e.g., [44] and [46]).

\strut \strut

In this paper, we observe how an algebraic structure, a graph groupoid
induced by a finite directed graph, acts on a Hilbert space $H.$ In [9] and
[10], we showed that every graph groupoid maybe realized with a system of
Hilbert space operators on a so-called graph Hilbert space.

\strut

Let $G$ be a countable directed graph with its vertex set $V(G)$ and its
edge set $E(G).$ Denote the set of all finite paths of $G$ by $FP(G).$
Clearly, the edge set $E(G)$ is contained in $FP(G).$ Let $w$ be a finite
path in $FP(G).$ Then it is represented as a word in $E(G).$ If $e_{1},$
..., $e_{n}$ are connected directed edges in the order $e_{1}$ $\rightarrow $
$e_{2}$ $\rightarrow $ ... $\rightarrow $ $e_{n},$ for $n$ $\in $ $\Bbb{N},$
then we can express $w$ by $e_{1}$ ... $e_{n}$ in $FP(G).$ If there exists a
finite path $w$ $=$ $e_{1}$ ... $e_{n}$ in $FP(G),$ where $n$ $\in $ $\Bbb{N}
$ $\setminus $ $\{1\},$ we say that the directed edges $e_{1},$ ..., $e_{n}$
are admissible. The length $\left| w\right| $ of $w$ is defined to be $n,$
which is the cardinality of the admissible edges generating $w.$ Also, we
say that finite paths $w_{1}$ $=$ $e_{11}$ ... $e_{1k_{1}}$ and $w_{2}$ $=$ $%
e_{21}$ ... $e_{2k_{2}}$ are admissible, if $w_{1}$ $w_{2}$ $=$ $e_{11}$ ... 
$e_{1k_{1}}$ $e_{21}$ ... $e_{2k_{2}}$ is again in $FP(G),$ where $e_{11},$
..., $e_{1k_{1}},$ $e_{21},$ ..., $e_{2k_{2}}$ $\in $ $E(G).$ Otherwise, we
say that $w_{1}$ and $w_{2}$ are not admissible. Suppose that $w$ is a
finite path in $FP(G)$ with its initial vertex $v_{1}$ and its terminal
vertex $v_{2}.$ Then we write $w$ $=$ $v_{1}w$ or $w$ $=$ $wv_{2}$ or $w$ $=$
$v_{1}wv_{2},$ for emphasizing the initial vertex of $w,$ respectively the
terminal vertex of $w,$ respectively both the initial vertex and the
terminal vertex of $w$. Suppose $w$ $=$ $v_{1}wv_{2}$ in $FP(G)$ with $%
v_{1}, $ $v_{2}$ $\in $ $V(G).$ Then we also say that [$v_{1}$ and $w$ are
admissible] and [$w$ and $v_{2}$ are admissible]. Notice that even though
the elements $w_{1}$ and $w_{2}$ in $V(G)$ $\cup $ $FP(G)$ are admissible, $%
w_{2}$ and $w_{1}$ are not admissible, in general. For instance, if $e_{1}$ $%
=$ $v_{1}e_{1}v_{2}$ is an edge with $v_{1}$, $v_{2}$ $\in $ $V(G)$ and $%
e_{2}$ $=$ $v_{2}e_{2}v_{3}$ is an edge with $v_{3}$ $\in $ $V(G)$ such that 
$v_{3}$ $\neq $ $v_{1},$ then there is a finite path $e_{1}$ $e_{2}$ in $%
FP(G),$ but there is no finite path $e_{2}$ $e_{1}$.

\strut

The \emph{free semigroupoid} $\Bbb{F}^{+}(G)$ of $G$ is defined by a set

\strut

\begin{center}
$\Bbb{F}^{+}(G)=\{\emptyset \}\cup V(G)\cup FP(G),$
\end{center}

\strut

with its binary operation $(\cdot )$ on $\Bbb{F}^{+}(G)$, defined by

\strut

\begin{center}
$(w_{1},w_{2})\mapsto w_{1}\cdot w_{2}=\left\{ 
\begin{array}{ll}
w_{1} & \text{if }w_{1}=w_{2}\text{ in }V(G) \\ 
w_{1} & \text{if }w_{1}\in FP(G),\text{ }w_{2}\in V(G)\text{ and }%
w_{1}=w_{1}w_{2} \\ 
w_{2} & \text{if }w_{1}\in V(G),\text{ }w_{2}\in FP(G)\text{ and }%
w_{2}=w_{1}w_{2} \\ 
w_{1}w_{2} & \text{if }w_{1},\text{ }w_{2}\text{ in }FP(G)\text{ and }%
w_{1}w_{2}\in FP(G) \\ 
\emptyset & \text{otherwise,}
\end{array}
\right. $
\end{center}

\strut

where $\emptyset $ is the empty word in $V(G)$ $\cup $ $E(G).$ (Sometimes,
the free semigroupoid $\Bbb{F}^{+}(G)$ of a certain graph $G$ does not
contain the empty word $\emptyset .$ For instance, the free semigroupoid of
the one-vertex--multi-loop-edge graph does not have the empty word. But, in
general, the empty word $\emptyset $ is contained in the free semigroupoid,
whenever $\left| V(G)\right| $ $>$ $1.$ So, if there is no confusion, then
we usually assume that the empty word is contained in free semigroupoids.)
This binary operation $(\cdot )$ on $\Bbb{F}^{+}(G)$ is called the \emph{%
admissibility}. i.e., the algebraic structure $(\Bbb{F}^{+}(G),$ $\cdot )$
is the free semigroupoid of $G.$ For convenience, we denote $(\Bbb{F}%
^{+}(G), $ $\cdot )$ simply by $\Bbb{F}^{+}(G).$

\strut

For the given countable directed graph $G,$ we can define a new countable
directed graph $G^{-1}$ which is the opposite directed graph of $G,$ with

\strut

\begin{center}
$V(G^{-1})=V(G)$ \ \ and $\ \ E(G^{-1})=\{e^{-1}:e\in E(G)\}$,
\end{center}

\strut

where $e^{-1}$ $\in $ $E(G^{-1})$ is the opposite directed edge of $e$ $\in $
$E(G)$, called the \emph{shadow} of $e$ $\in $ $E(G).$ i.e., if $e$ $=$ $%
v_{1}$ $e$ $v_{2}$ in $E(G)$ with $v_{1},$ $v_{2}$ $\in $ $V(G),$ then $%
e^{-1}$ $=$ $v_{2}$ $e^{-1}$ $v_{1}$ in $E(G^{-1})$ with $v_{2},$ $v_{1}$ $%
\in $ $V(G^{-1})$ $=$ $V(G).$ This new directed graph $G^{-1}$ is said to be 
\emph{the shadow of} $G$. It is trivial that $(G^{-1})^{-1}$ $=$ $G.$ This
relation shows that the admissibility on the shadow $G^{-1}$ is oppositely
preserved by that on $G.$

\strut

A new countable directed graph $\widehat{G}$ is called the \emph{shadowed
graph} \emph{of} $G$ if it is a directed graph with

\strut

\begin{center}
$V(\widehat{G})=V(G)=V(G^{-1})$
\end{center}

and

\begin{center}
$E(\widehat{G})=E(G)\cup E(G^{-1}).$
\end{center}

\strut \strut \strut

\begin{definition}
Let $G$ be a countable directed graph and $\widehat{G},$ the shadowed graph
of $G,$ and let $\Bbb{F}^{+}(\widehat{G})$ be the free semigroupoid of $%
\widehat{G}.$ Define the reduction (RR) on $\Bbb{F}^{+}(\widehat{G})$ by

\strut 

(RR)$\ \ \ \ \ \ \ \ \ \ \ \ \ \ \ \ \ \ ww^{-1}=v$ \ \ and \ \ $%
w^{-1}w=v^{\prime },$

\strut 

whenever $w$ $=$ $vwv^{\prime }$ in $FP(\widehat{G}),$ with $v,$ $v^{\prime }
$ $\in $ $V(\widehat{G})$. The set $\Bbb{F}^{+}(\widehat{G})$ with this
reduction (RR) is denoted by $\Bbb{F}_{r}^{+}(\widehat{G}).$ And this set $%
\Bbb{F}_{r}^{+}(\widehat{G})$ with the inherited admissibility $(\cdot )$
from $\Bbb{F}^{+}(\widehat{G})$ is called the graph groupoid of $G.$ Denote $%
(\Bbb{F}_{r}^{+}(\widehat{G}),$ $\cdot )$ of $G$ by $\Bbb{G}.$ Define the
reduced finite path set $FP_{r}(\widehat{G})$ of $\Bbb{G}$ by

\strut 

\begin{center}
$FP_{r}(\widehat{G})$ $\overset{def}{=}$ $\Bbb{G}$ $\setminus $ $(V(\widehat{%
G})$ $\cup $ $\{\emptyset \}).$
\end{center}

\strut 

All elements of $FP_{r}(\widehat{G})$ are said to be reduced finite paths of 
$\widehat{G}.$
\end{definition}

\strut \strut \strut \strut \strut

Let $w_{1}$ and $w_{2}$ be reduced finite paths in $FP_{r}(\widehat{G})$ $%
\subset $ $\Bbb{G}.$ We will use the same notation $w_{1}$ $w_{2}$ for the
(reduced) product of $w_{1}$ and $w_{2}$ in $\Bbb{G}.$ But we have to keep
in mind that the product $w_{1}$ $w_{2}$ in the graph groupoid $\Bbb{G}$ is
different from the (non-reduced) product $w_{1}$ $w_{2}$ in the free
semigroupoid $\Bbb{F}^{+}(\widehat{G}).$ Suppose $e_{1}$ and $e_{2}$ are
edges in $E(\widehat{G})$ and assume that they are admissible, and hence $%
e_{1}$ $e_{2}$ is a finite path in $FP(\widehat{G}).$ Then the product $%
e_{1} $ $e_{2}$ $e_{2}^{-1}$ of $e_{1}$ $e_{2}$ and $e_{2}^{-1}$ is a
length-3 finite path in $FP(\widehat{G})$ $\subset $ $\Bbb{F}^{+}(\widehat{G}%
),$ but the product $e_{1}$ $e_{2}$ $e_{2}^{-1}$ of them is $e_{1}$ $(e_{2}$ 
$e_{2}^{-1})$ $=$ $e_{1},$ which is a length-1 reduced finite path, in $%
FP_{r}(\widehat{G})$ $\subset $ $\Bbb{G}.$

\strut

\strut

\strut

\subsection{$C^{*}$-Algebras Induced by Graphs}

\strut

\strut \strut \strut \strut \strut \strut \strut \strut

In [13], we observed the $C^{*}$- algebras generated by certain family of
partial isometries on a fixed Hilbert space $H$.\ \strut Suppose a family $%
\mathcal{G}$ $=$ $\{a_{1},$ ..., $a_{N}\}$ in $B(H)$ is a collection of
finite number of partial isometries, where $N$ $\in $ $\Bbb{N}.$ We say that
such family $\mathcal{G}$ \emph{construct a finite directed graph} $G$ if
there exists a finite directed graph $G$ such that

\strut

(i) $\left| E(G)\right| $ $=$ $\left| \mathcal{G}\right| $, and $\left|
V(G)\right| $ $=$ $\left| \mathcal{G}_{pro}\right| ,$ where

$\strut $

\begin{center}
$\mathcal{G}_{pro}$ $=$ $\{a^{*}$ $a$ $:$ $a$ $\in $ $\mathcal{G}\}$ $\cup $ 
$\{a$ $a^{*}$ $:$ $a$ $\in $ $\mathcal{G}\},$
\end{center}

\strut

(ii) the edges $e_{1}$ and $e_{2}$ create a nonempty finite path $e_{1}$ $%
e_{2}$ on $G$ if and only if $H_{init}^{a_{1}}$ $=$ $H_{fin}^{a_{2}},$ where
``$=$'' means ``being identically same in $H$''.

\strut

In this case, the family $\mathcal{G}$ is called a $G$-\emph{family of
partial isometries} in $B(H)$. The family $\widehat{\mathcal{G}}$ $=$ $%
\mathcal{G}$ $\cup $ $\mathcal{G}^{*}$ generates a $C^{*}$-subalgebra $C^{*}(%
\mathcal{G})$ of $B(H),$ where $\mathcal{G}^{*}$ $=$ $\{a_{j}^{*}$ $:$ $j$ $%
= $ $1,$ ..., $N\}.$

\strut

Let $\Bbb{G}$ be the graph groupoid of a finite directed graph $G,$ with $%
\left| V(\widehat{G})\right| $ $=$ $n$ and $\left| E(\widehat{G})\right| $ $%
= $ $2N$ (i.e., $\left| E(G)\right| $ $=$ $N$). We give a indexing on $V(G)$
by $\{1,$ ..., $n\}$. i.e., we will let $V(G)$ $=$ $\{v_{1},$ ..., $v_{n}\}.$
By doing that, we can also index the elements $E(G)$ and $E(\widehat{G})$ as
follows:

\strut

\begin{center}
$E(G)$ $=$ $\left\{ e_{m:ij}\left| 
\begin{array}{c}
m=1,\text{ ..., }k_{ij},\text{ }k_{ij}\neq 0 \\ 
e_{m:ij}=v_{i}e_{m:ij}v_{j}
\end{array}
\right. \right\} $
\end{center}

and

\begin{center}
$E(\widehat{G})$ $=$ $\left\{ x_{m:ij}\left| 
\begin{array}{c}
x_{m:ij}=e_{m:ij}\text{ if }x_{m:ij}\in E(G) \\ 
x_{m:ij}=e_{m:ji}^{-1}\text{ if }x_{m:ij}\in E(G^{-1}) \\ 
m=1,...,k_{ij},\text{ }k_{ij}\neq 0
\end{array}
\right. \right\} ,$
\end{center}

\strut

where $k_{ij}$ means the cardinality of edges connecting the vertex $v_{i}$
to $v_{j}.$ By the finiteness of $G,$ $k_{ij}$ $<$ $\infty ,$ whenever $%
k_{ij}$ $\neq $ $0.$ Clearly, ``$k_{ij}$ $=$ $0$'' means that ``there is no
edge connecting $v_{i}$ to $v_{j}$''. In [13], we showed that such indexing
process on $V(G)$ $\cup $ $E(G)$ is uniquely determined up to
graph-isomorphisms. This means that if we fix an indexing, then this
indexing contains the full (admissibility) combinatorial data of $G.$

\strut

We will say that a finite directed graph $G$ is \emph{connected}, if its
graph groupoid $\Bbb{G}$ satisfies that: for any $(v_{1},$ $v_{2})$ $\in $ $%
V(G)$ $\times $ $V(G),$ where $v_{1}$ $\neq $ $v_{2},$ there exists an
element $w$ $\in $ $\Bbb{G}$ such that $w$ $=$ $v_{1}$ $w$ $v_{2}$ and $%
w^{-1}$ $=$ $v_{2}$ $w^{-1}$ $v_{1}.$ Otherwise, we say that the graph $G$
is not connected or is disconnected.

\strut \strut

Under this setting, if $G$ is a connected finite directed graph, then the
graph groupoid $\Bbb{G}$ has its matricial representation $(\mathcal{H}_{G},$
$\pi ),$ where $\mathcal{H}_{G}$ is a Hilbert space defined by $\oplus
_{j=1}^{n}$ $(\Bbb{C}$ $\xi _{v_{j}}),$ which is Hilbert-space isomorphic to 
$\Bbb{C}^{\oplus \,n},$ and where $\pi $ $:$ $\Bbb{G}$ $\rightarrow $ $B(%
\mathcal{H}_{G})$ is a groupoid action satisfying that

\strut

\begin{center}
$\pi (v_{j})$ $=$ $P_{j}$
\end{center}

and

\begin{center}
$\pi (x_{m:ij})$ $=$ $\left\{ 
\begin{array}{ll}
E_{m:ij} & \text{if }x_{m:ij}=e_{m:ij} \\ 
&  \\ 
E_{m:ji}^{*} & \text{if }x_{m:ij}=e_{m:ji}^{-1}\text{ ,}
\end{array}
\right. $
\end{center}

\strut

where $P_{j}$ is the diagonal matrix in $M_{n}(\Bbb{C}),$ having its only
nonzero $(j,$ $j)$-entry $1,$

\strut

\begin{center}
$
\begin{array}{ll}
\qquad \qquad \;j\text{-th} &  \\ 
\left( 
\begin{array}{lllll}
0 &  &  &  & 0 \\ 
& \ddots &  &  &  \\ 
&  & 1 &  &  \\ 
&  &  & \ddots &  \\ 
0 &  &  &  & 0
\end{array}
\right) & j\text{-th}
\end{array}
$
\end{center}

\strut \strut \strut

and $E_{m:ij}$ is the matrix in $M_{n}(\Bbb{C})$, having its only nonzero $%
(i,$ $j)$-entry $\omega ^{m}$, where $\omega $ is the root of unity of the
polynomial $z^{k_{ij}},$ whenever $i$ $\neq $ $j,$

\strut

\begin{center}
$
\begin{array}{ll}
\qquad \;i\text{-th} &  \\ 
\left( 
\begin{array}{lllll}
0 &  &  &  & 0 \\ 
&  &  &  &  \\ 
& \omega ^{m} &  &  &  \\ 
&  &  &  &  \\ 
0 &  &  &  & 0
\end{array}
\right) & j\text{-th}
\end{array}
$
\end{center}

\strut \strut \strut

or it is the diagonal matrix in $M_{n}(\Bbb{C})$, having its only nonzero $%
(j,$ $j)$-entry $e^{i\,\theta _{m:jj}},$ where $\theta _{m:jj}$ $\in $ $\Bbb{%
R}$ $\setminus $ $\{0\},$ satisfying that $\theta _{m_{1}:jj}$ $\neq $ $%
\theta _{m_{2}:jj},$ whenever $m_{1}$ $\neq $ $m_{2}.$

\strut

\begin{center}
$
\begin{array}{ll}
\qquad \qquad \;j\text{-th} &  \\ 
\left( 
\begin{array}{lllll}
0 &  &  &  & 0 \\ 
& \ddots &  &  &  \\ 
&  & e^{i\theta _{jj}} &  &  \\ 
&  &  & \ddots &  \\ 
0 &  &  &  & 0
\end{array}
\right) & j\text{-th}
\end{array}
$
\end{center}

\strut \strut

Then we can see that the graph groupoid $\Bbb{G}$ generates the $C^{*}$-
algebra $\mathcal{M}_{G}$ which is a $C^{*}$-subalgebra of $B(\mathcal{H}%
_{G})$ $\overset{*\text{-isomorphic}}{=}$ $M_{n}(\Bbb{C}).$ This $C^{*}$%
-subalgebra $\mathcal{M}_{G}$ induced by $G$ in $M_{n}(\Bbb{C})$ is called
the \emph{matricial graph }$C^{*}$\emph{-algebra} of a connected finite
directed graph $G.$

\strut

If the family $\mathcal{G}$ of partial isometries in $B(H)$ constructs a
connected finite directed graph $G,$ then, since $\mathcal{G}$ generates a
groupoid $\Bbb{G}_{\mathcal{G}},$ under the operator multiplication on $%
B(H), $ and it is groupoid-isomorphic to a certain graph groupoid $\Bbb{G},$
we can check that the $C^{*}$-subalgebra $C^{*}(\mathcal{G})$ of $B(H)$ is $%
* $-isomorphic to the affiliated matricial graph $C^{*}$- algebra

$\strut $

\begin{center}
$\mathcal{M}_{G}(H_{0})$ $\overset{def}{=}$ $(\Bbb{G}$ $\cdot $ $1_{H_{0}})$ 
$\otimes _{\Bbb{C}}$ $\mathcal{M}_{G}$ $\overset{C^{*}\text{-subalgebra}}{%
\subseteq }$ $M_{n}(B(H_{0})),$
\end{center}

\strut

for certain Hilbert space $H_{0}$ which is embedded in $H,$ where

\strut

\begin{center}
$M_{n}\left( B(H_{0})\right) $ $=$ $\left\{ \left. \left( 
\begin{array}{lll}
T_{11} & \cdots & T_{1n} \\ 
\,\,\,\vdots & \ddots & \,\,\,\vdots \\ 
T_{n1} & \cdots & T_{nn}
\end{array}
\right) \right| T_{ij}\in B(H_{0})\right\} .$
\end{center}

\strut

We conclude this section with two results which will be extended later.

\strut

\begin{theorem}
(See [13]) Let $\mathcal{G}$ be a finite family of partial isometries in $%
B(H),$ and assume that $\mathcal{G}$ constructs a connected finite directed
graph $G.$ Then the $C^{*}$-subalgebra $C^{*}(\mathcal{G})$ of $B(H)$ is $*$%
-isomorphic to the affiliated matricial algebra $\mathcal{M}_{G}(H_{0}),$
where $H_{0}$ is a Hilbert space which is Hilbert-space isomorphic to $p$ $H,
$ for all $p$ $\in $ $\mathcal{G}_{pro}.$ $\square $
\end{theorem}

\strut

\begin{corollary}
(Also See [13]) Let $\mathcal{G}$ be a finite family of partial isometries
in $B(H),$ and assume that $\mathcal{G}$ constructs a finite directed graph $%
G,$ having its connected components $G_{1},$ ..., $G_{t},$ for $t$ $\in $ $%
\Bbb{N}.$ Then the $C^{*}$-subalgebra $C^{*}(\mathcal{G})$ of $B(H)$ is $*$%
-isomorphic to the direct product algebra $\oplus _{j=1}^{t}$ $\mathcal{M}%
_{G}(H_{j})$ of the affiliated matricial graph $C^{*}$-algebras $\mathcal{M}%
_{G_{j}}(H_{j})$ $\overset{def}{=}$ $(\Bbb{C}$ $\cdot $ $1_{H_{j}})$ $%
\otimes _{\Bbb{C}}$ $\mathcal{M}_{G_{j}},$ where $\mathcal{M}_{G_{j}}$'s are
the matricial $C^{*}$-algebras induced by $G_{j},$ for $j$ $=$ $1,$ ..., $N.$
$\square $
\end{corollary}

\strut

The above theorem and corollary are the main motivation of this paper. They
provide a certain connection between partial isometries and graph groupoids.
We will extend the above results. We realize that, in general, we have much
more complicated structures than the above special cases.

\strut

\strut

\subsection{Groupoids and Groupoid Actions}

\strut

\strut

As we observed in [9], [10] and [12], every graph groupoid is indeed a
(categorial) groupoid. This means that graph groupoids have rough but rich
algebraic structures. The following definition is inspired by [19] and [25].

\strut

\begin{definition}
We say an algebraic structure $(\mathcal{X},$ $\mathcal{Y},$ $s,$ $r)$ is a
(categorial) groupoid if it satisfies that (i) $\mathcal{Y}$ $\subset $ $%
\mathcal{X},$ (ii) for all $x_{1},$ $x_{2}$ $\in $ $\mathcal{X},$ there
exists a partially-defined binary operation $(x_{1},$ $x_{2})$ $\mapsto $ $%
x_{1}$ $x_{2},$ for all $x_{1},$ $x_{2}$ $\in $ $\mathcal{X},$ depending on
the source map $s$ and the range map $r$ satisfying that:

\strut \strut 

(ii-1) the operation $x_{1}$ $x_{2}$ is well-determined, whenever $r(x_{1})$ 
$s(x_{2})$ $\in $ $\mathcal{Y}$,

\strut 

(ii-2) the operation $(x_{1}$ $x_{2})$ $x_{3}$ $=$ $x_{1}$ $(x_{2}$ $x_{3})$
is defined, if the constituents are well-determined in the sense of (ii-1),
for $x_{1},$ $x_{2},$ $x_{3}$ $\in $ $\mathcal{X},$

\strut 

(ii-3) if $x$ $\in $ $\mathcal{X},$ then there exist $y,$ $y^{\prime }$ $\in 
$ $\mathcal{Y}$ such that $s(x)$ $=$ $y$ and $r(x)$ $=$ $y^{\prime },$
satisfying $x$ $=$ $y$ $x$ $y^{\prime },$

\strut 

(ii-4) if $x$ $\in $ $\mathcal{X},$ then there exists a unique
groupoid-inverse $x^{-1}$ satisfying $x$ $x^{-1}$ $=$ $s(x)$ and $x^{-1}$ $x$
$=$ $r(x),$ in $\mathcal{Y}.$
\end{definition}

\strut \strut

By definition, we can conclude that every group is a groupoid $(\mathcal{X},$
$\mathcal{Y},$ $s,$ $r)$ with $\left| \mathcal{Y}\right| $ $=$ $1$ (and
hence $s$ $=$ $r$ on $\mathcal{X}$). This subset $\mathcal{Y}$ of $\mathcal{X%
}$ is said to be the base of $\mathcal{X}$. Remark that we can naturally
assume that there exists the empty element $\emptyset $ in a groupoid $%
\mathcal{X}.$ Roughly speaking, the empty element $\emptyset $ means that
there are products $x_{1}$ $x_{2}$ which are not well-defined, for $x_{1},$ $%
x_{2}$ $\in $ $\mathcal{X}.$ Notice that if $\left| \mathcal{Y}\right| $ $=$ 
$1$ (equivalently, if $\mathcal{X}$ is a group), then the empty word $%
\emptyset $ is not contained in the groupoid $\mathcal{X}.$ However, in
general, whenever $\left| \mathcal{Y}\right| $ $\geq $ $2,$ a groupoid $%
\mathcal{X}$ always contain the empty word. So, if there is no confusion, we
will automatically assume that the empty element $\emptyset $ is contained
in $\mathcal{X}.$

\strut

It is easily checked that our graph groupoid $\Bbb{G}$ of a finite directed
graph $G$ is indeed a groupoid with its base $V(\widehat{G}).$ i.e., every
graph groupoid $\Bbb{G}$ of a countable directed graph $G$ is a groupoid $(%
\Bbb{G},$ $V(\widehat{G}),$ $s$, $r)$, where $s(w)$ $=$ $s(v$ $w)$ $=$ $v$
and $r(w)$ $=$ $r(w$ $v^{\prime })$ $=$ $v^{\prime },$ for all $w$ $=$ $v$ $%
w $ $v^{\prime }$ $\in $ $\Bbb{G}$ with $v,$ $v^{\prime }$ $\in $ $V(%
\widehat{G}).$

\strut

Let $\mathcal{X}_{k}$ $=$ $(\mathcal{X}_{k},$ $\mathcal{Y}_{k},$ $s_{k},$ $%
r_{k})$ be groupoids, for $k$ $=$ $1,$ $2.$ We say that a map $f$ $:$ $%
\mathcal{X}_{1}$ $\rightarrow $ $\mathcal{X}_{2}$ is a \emph{groupoid
morphism} if (i) $f(\mathcal{Y}_{1})$ $\subseteq $ $\mathcal{Y}_{2},$ (ii) $%
s_{2}\left( f(x)\right) $ $=$ $f\left( s_{1}(x)\right) $ in $\mathcal{X}%
_{2}, $ for all $x$ $\in $ $\mathcal{X}_{1}$, and (iii) $r_{2}\left(
f(x)\right) $ $=$ $f\left( r_{1}(x)\right) $ in $\mathcal{X}_{2},$ for all $%
x $ $\in $ $\mathcal{X}_{1}.$ If a groupoid morphism $f$ is bijective, then
we say that $f$ is a \emph{groupoid-isomorphism}, and that the groupoids $%
\mathcal{X}_{1}$ and $\mathcal{X}_{2}$ are \emph{groupoid-isomorphic}.

\strut

Notice that, if two countable directed graphs $G_{1}$ and $G_{2}$ are
graph-isomorphic, via a graph-isomorphism $g$ $:$ $G_{1}$ $\rightarrow $ $%
G_{2},$ in the sense that (i) $g$ is bijective from $V(G_{1})$ onto $%
V(G_{2}),$ (ii) $g$ is bijective from $E(G_{1})$ onto $E(G_{2}),$ (iii) $g(e)
$ $=$ $g(v_{1}$ $e$ $v_{2})$ $=$ $g(v_{1})$ $g(e)$ $g(v_{2})$ in $E(G_{2}),$
for all $e$ $=$ $v_{1}$ $e$ $v_{2}$ $\in $ $E(G_{1}),$ with $v_{1},$ $v_{2}$ 
$\in $ $V(G_{1}),$ then the graph groupoids $\Bbb{G}_{1}$ and $\Bbb{G}_{2}$
are groupoid-isomorphic. More generally, if $G_{1}$ and $G_{2}$ have
graph-isomorphic shadowed graphs $\widehat{G}$ and $\widehat{G_{2}},$ then $%
\Bbb{G}_{1}$ and $\Bbb{G}_{2}$ are groupoid-isomorphic (See Section 3, in
detail).

\strut

Let $\mathcal{X}$ $=$ $(\mathcal{X},$ $\mathcal{Y},$ $s,$ $r)$ be a
groupoid. We say that this groupoid $\mathcal{X}$ acts on a set $X$ if there
exists a groupoid action $\pi $ of $\mathcal{X}$ such that $\pi (x)$ $:$ $X$ 
$\rightarrow $ $X$ is a well-determined function. Sometimes, we call the set 
$X,$ a $\mathcal{X}$-set. We are interested in the case where a $\mathcal{X}$%
-set $X$ is a Hilbert space. A nice example of groupoid actions acting on a
Hilbert space are graph-representations defined and observed in [9]. In this
paper, we will consider a new kind of graph-groupoid action on certain
Hilbert spaces.

\strut \strut

\strut \strut \strut \strut

\strut

\section{$C^{*}$-Subalgebras of $B(H)$ Generated by Partial Isometries}

\strut

\strut

\strut

The section has two parts, one centered around our Theorem 3.3. This result
generalizes the more traditional notion of deficiency indices (See Section
5). The second part shows that our $C^{*}$-algebras can be realized in
representations of groupoids. 

\strut 

In [13], we considered the case where the collection of finite partial
isometries are graph-families. How about the general cases where the family
of finite partial isometries on $H$ is not a graph-family in the sense of
[13]? We provide an answer of this question. As before, we will fix a
Hilbert space $H$ and the corresponding operator algebra $B(H).$ Let $%
\mathcal{G}$ $=$ $\{a_{1},$ ..., $a_{N}\}$ be a family of partial isometries 
$a_{1},$ ..., $a_{N}$ of $B(H),$ where $N$ $\in $ $\Bbb{N}.$ Compared with
[13], we will have much more complicated results.

\strut

\strut

\strut

\subsection{$*$-Isomorphic Indices of Partial Isometries}

\strut

\strut

In this section, we will classify the set of all partial isometries of $B(H)$%
. For convenience, we denote the subset of all partial isometries of $B(H)$
by $PI(H).$ Let $a$ $\in $ $PI(H)$ with its initial space $H_{init}^{a}$ and
its final space $H_{fin}^{a}.$ Recall that $H_{init}^{a}$ $=$ $(a^{*}$ $a)$ $%
H$ and $H_{fin}^{a}$ $=$ $(a$ $a^{*})$ $H.$ Then we can find the subspaces

$\strut $

\begin{center}
$(H_{init}^{a})^{\perp }$ $=$ $H$ $\ominus $ $H_{init}^{a}$ $=$ $(1_{H}$ $-$ 
$a^{*}a)$ $H$ $=$ $\ker a$
\end{center}

and

\begin{center}
$(H_{fin}^{a})^{\perp }$ $=$ $H$ $\ominus $ $H_{fin}^{a}$ $=$ $(1_{H}$ $-$ $%
a $ $a^{*})$ $H$ $=$ $\ker a^{*}.$
\end{center}

\strut

Assume that

\strut

\begin{center}
$\varepsilon _{+}$ $\overset{def}{=}$ $\dim \left( (H_{init}^{a})^{\perp
}\right) $ and $\varepsilon _{-}$ $\overset{def}{=}$ $\dim $ $\left(
(H_{fin}^{a})^{\perp }\right) ,$
\end{center}

\strut

where $\varepsilon _{+},$ $\varepsilon _{-}$ $\in $ $\Bbb{N}_{0}^{\infty }$ $%
=$ $\Bbb{N}$ $\cup $ $\{0,$ $\infty \}.$ Then we can determine a subset $%
PI_{H}(\varepsilon _{+},$ $\varepsilon _{-})$ of $PI(H)$ by

\strut

\begin{center}
$PI_{H}(\varepsilon _{+},$ $\varepsilon _{-})$ $\overset{def}{=}$ $\left\{
a\in PI(H)\left| 
\begin{array}{c}
\dim \left( (H_{init}^{a})^{\perp }\right) =\varepsilon _{+}, \\ 
\dim \left( (H_{fin}^{a})^{\perp }\right) =\varepsilon _{-}
\end{array}
\right. \right\} .$
\end{center}

\strut

Note that

\strut

\begin{center}
$PI(H)$ $=$ $\underset{(\varepsilon _{+},\text{ }\varepsilon _{-})\in (\Bbb{N%
}_{0}^{\infty })^{2}}{\cup }$ $\left( PI_{H}(\varepsilon _{+},\text{ }%
\varepsilon _{-})\right) ,$
\end{center}

\strut

set-theoretically, in $B(H),$ where $\Bbb{N}_{0}^{\infty }$ $=$ $\Bbb{N}$ $%
\cup $ $\{0,$ $\infty \}.$ For instance, if $U$ is the well-known unilateral
shift (which is an isometry) on the Hilbert space $l^{2}(\Bbb{N}_{0}),$
sending $(\xi _{0},$ $\xi _{1},$ $\xi _{2},$ ...$)$ $\in $ $l^{2}(\Bbb{N}%
_{0})$ to

\strut

\begin{center}
$U\left( (\xi _{0},\text{ }\xi _{1},\text{ }\xi _{2},\text{ ...})\right) $ $%
\overset{def}{=}$ $(0,$ $\xi _{0},$ $\xi _{1},$ $\xi _{2},$ ...$)$ $\in $ $%
l^{2}(\Bbb{N}_{0}),$
\end{center}

\strut

where $\Bbb{N}_{0}$ $\overset{def}{=}$ $\Bbb{N}$ $\cup $ $\{0\}$, then it is
contained in $PI_{l^{2}(\Bbb{N}_{0})}(0,$ $1),$ since $H_{init}^{U}$ $=$ $%
l^{2}(\Bbb{N}_{0})$ and $H_{fin}^{U}$ $=$ $l^{2}(\Bbb{N}_{0})$ $\ominus $ $%
\Bbb{C}.$ Notice that every operator $U^{\varepsilon }$ on $l^{2}(\Bbb{N}%
_{0})$ can be regarded as elements in $PI_{l^{2}(\Bbb{N}_{0})}(0,$ $%
\varepsilon ),$ for all $\varepsilon $ $\in $ $\Bbb{N},$ and every unitaries
on an arbitrary Hilbert space $K$ can be regarded as elements in $PI_{K}(0,$ 
$0).$

\strut

Let $u_{1}$ $\neq $ $u_{2}$ $\in $ $PI_{H}(0,$ $0)$ be unitaries, satisfying
that $spec(u_{1})$ $=$ $spec(u_{2})$ in $\Bbb{T}$ $\subset $ $\Bbb{C},$
where $\Bbb{T}$ is the unit circle in $\Bbb{C}$ and $spec(x)$ means the
spectrum of $x$, for all $x$ $\in $ $B(H).$ Then we can get that:

\strut

\begin{lemma}
If the subspaces $H_{u_{1}}$ and $H_{u_{2}}$ of $H$ are Hilbert-space
isomorphic, where $H_{u_{k}}$ $=$ $(u_{k}^{*}$ $u_{k})$ $H,$ for $k$ $=$ $1,$
$2,$ and if the spectrums $spec(u_{1})$ and $spec(u_{2})$ are identical in $%
\Bbb{T}$ $\subset $ $\Bbb{C},$ then the $C^{*}$-subalgebras $C^{*}(\{u_{1}\})
$ and $C^{*}(\{u_{2}\})$ of $B(H)$ are $*$-isomorphic.
\end{lemma}

\strut

\begin{proof}
By [13], we have that the $C^{*}$-algebras $C^{*}(\{u_{k}\})$ are $*$%
-isomorphic to $(\Bbb{C}$ $\cdot $ $1_{H_{u_{k}}})$ $\otimes _{\Bbb{C}}$ $%
C\left( spec(u_{k})\right) ,$ as embedded $C^{*}$-subalgebra of $B(H),$ for $%
k$ $=$ $1,$ $2,$ where $C(X)$ means the $C^{*}$-algebra consisting of all
continuous functions on a compact subset $X$ of $\Bbb{C}.$ Since the
subspaces $H_{u_{1}}$ and $H_{u_{2}}$ are Hilbert-space isomorphic in $H$
and since the spectrums $spec(u_{1})$ and $spec(u_{2})$ are same in $\Bbb{T}%
, $ the $C^{*}$-algebras $(\Bbb{C}$ $\cdot $ $1_{H_{k}})$ $\otimes _{\Bbb{C}%
} $ $C\left( spec(u_{k})\right) $'s are $*$-isomorphic $C^{*}$\strut
-subalgebras of $B(H).$
\end{proof}

\strut \strut \strut \strut

If $a$ $\in $ $PI(H)$ is a partial isometry on $H$, then we can regard $a$
as an isometry on $H$ $\ominus $ $\ker a.$ Recall that if $x$ $\in $ $B(K)$
is an isometry on a Hilbert space $K,$ then $x$ has its Wold decomposition $%
x $ $=$ $u_{x}$ $+$ $s_{x},$ where $u_{x}$ is the unitary part of $x$ and $%
s_{x}$ is the shift part of $x.$ Futhermore, the Hilbert space $K$ is
decomposed by $K$ $=$ $K_{u_{x}}$ $\oplus $ $K_{s_{x}},$ where $K_{u_{x}}$
(resp., $K_{s_{x}}$) is the subspace of $K$ where $u_{x}$ (resp., $s_{x}$)
is acting on, as a unitary (resp., as a shift).

\strut

\begin{definition}
Suppose $a$ $\in $ $PI(H)$ is a partial isometry. We say that $u_{a}$ $+$ $%
s_{s}$ $\in $ $B(H)$ is the Wold decomposition of $a,$ if $u_{a}$ is the
unitary part of $a$ on $H_{a}$ and $s_{a}$ is the shift part of $a$ on $%
H_{a},$ by regarding $a$ as an isometry on $H_{a},$ where $H_{a}$ $=$ $(a^{*}
$ $a)$ $H$ $=$ $H$ $\ominus $ $\ker a.$ Equivalently, the Wold decomposition
of $a$ $\mid _{H_{a}}$ of the isometry $a$ $\mid _{H_{a}}$ is said to be the
Wold decomposition of the partial isometry $a.$
\end{definition}

\strut

Let $a$ $\in $ $PI_{H}(\varepsilon _{+},$ $\varepsilon _{-}).$ Then, by the
Wold decomposition, we can decompose $a$

\strut

\begin{center}
$a$ $=$ $u_{a}$ $+$ $s_{a},$
\end{center}

\strut

where $u_{a}$ is the unitary part on $H_{u_{a}}$ $=$ $(u_{a}^{*}$ $u_{a})$ $%
H $ of $a$; and $s_{a}$ is the shift part on $H_{s_{a}}$ $=$ $(s_{a}^{*}$ $%
s_{a})$ $H$ of $a.$ Note that the subspaces $H_{u_{a}}$ and $H_{s_{a}}$ of $%
H $ are also the subspaces of the initial space $H_{init}^{a}$ of $a,$ and
they decompose $H_{init}^{a}$. i.e., $H_{init}^{a}$ $=$ $H_{u_{a}}$ $\oplus $
$H_{s_{a}}.$ We will call the subspaces $H_{u_{a}}$ and $H_{s_{a}}$ of $H$,
the \emph{unitary part of} $H$ and the\emph{\ shift part of} $H,$ \emph{in
terms of} $a.$

\strut

By the Wold decomposition of $a,$ we then get the following diagram:

\strut

\begin{center}
$
\begin{array}{lllllll}
H & = & H_{u_{a}} & \oplus & \,\,\,\,\,\,\,\,\,\,\,\,\,\,\,H_{s_{a}} & \oplus
& \ker a \\ 
\,\,\downarrow _{a} &  & \,\,\downarrow _{u_{a}} &  & \,\,\,\,\,\,\,\,\,\,\,%
\,\,\,\,\downarrow _{s_{a}} &  &  \\ 
H & = & H_{u_{a}} & \oplus & (\ker s_{a}^{*}\oplus H_{s_{a}}) & \oplus & 
\ker a^{*}\ominus \ker s_{a}^{*}
\end{array}
$
\end{center}

\strut

i.e., the partial isometry $a$ has the following block operator-matricial
form,

\strut

\begin{center}
$a$ $=$ $\left( 
\begin{array}{lll}
u_{a} &  &  \\ 
& s_{a} &  \\ 
&  & 0
\end{array}
\right) ,$
\end{center}

\strut

on $H_{u_{a}}$ $\oplus $ $H_{s_{a}}$ $\oplus $ $\ker a.$

\strut

Let $a$ $\in $ $PI_{H}(\varepsilon _{+},$ $\varepsilon _{-}),$ with its Wold
decomposition $a$ $=$ $u_{a}$ $+$ $s_{a},$ and assume that the shift part $%
s_{a}$ is contained in $PI_{H_{s_{a}}}(0,$ $\varepsilon ^{-}).$ Then we have
that

\strut

\begin{center}
$\varepsilon _{-}$ $=$ $\varepsilon ^{-}$ $+$ $\varepsilon _{-}^{-},$ for
some $\varepsilon _{-}^{-}$ $\in $ $\Bbb{N}_{0}^{\infty },$
\end{center}

\strut where

\begin{center}
$\varepsilon ^{-}$ $=$ $\dim \left( \ker s_{a}^{*}\right) $ and $\varepsilon
_{-}^{-}$ $=$ $\varepsilon _{-}$ $-$ $\varepsilon ^{-}.$
\end{center}

\strut

\begin{definition}
Let $a$ $\in $ $PI_{H}(\varepsilon _{+},$ $\varepsilon _{-})$ with its Wold
decomposition $u_{a}$ $+$ $s_{a},$ and assume that $s_{a}$ $\in $ $%
PI_{H_{s_{a}}}(0,$ $\varepsilon ^{-}).$ Then we will denote the collection
of such partial isometries $a$ by $PI_{H}(\varepsilon _{+},$ $\varepsilon
^{-},$ $\varepsilon _{-}^{-}).$
\end{definition}

\strut

Let $a$ $\in $ $PI_{H}(\varepsilon _{+},$ $\varepsilon ^{-},$ $\varepsilon
_{-}^{-})$ be a partial isometry. Then we can have the following properties.

\strut

\begin{proposition}
Let $a$ $\in $ $PI_{H}(\varepsilon _{+},$ $\varepsilon ^{-},$ $\varepsilon
_{-}^{-})$ be a nonzero partial isometry with its Wold decomposition $a$ $=$ 
$u$ $+$ $s,$ where $u$ and $s$ are the unitary part and the shift part of $a,
$ respectively.

\strut \strut 

(1) We can have that

\strut 

\begin{center}
$C^{*}(\{a\})$ $\overset{*\text{-isomorphic}}{=}$ $C^{*}(\{u,$ $s\})$ $=$ $%
C^{*}(\{u\})$ $\oplus $ $C^{*}(\{s\})$ $\overset{C^{*}\text{-subalgebra}}{%
\subseteq }$ $B(H).$
\end{center}

\strut 

(2) If $s$ $=$ $0,$ then

\strut 

\begin{center}
$C^{*}(\{a\})$ $=$ $C^{*}(\{u\})$ $\overset{*\text{-isomorphic}}{=}$ $(\Bbb{C%
}$ $\cdot $ $1_{H_{u}})$ $\otimes _{\Bbb{C}}$ $C\left( spec(u)\right) .$
\end{center}

\strut 

(3) If $u$ $=$ $0$ and if $s$ $\in $ $PI_{H_{s}}(0,$ $\varepsilon ^{-})$
with $\varepsilon ^{-}$ $\in $ $\Bbb{N},$ then

\strut 

\begin{center}
$C^{*}(\{a\})$ $=$ $C^{*}(\{s\})$ $\overset{*\text{-isomorphic}}{=}$ $%
\mathcal{T}(H_{s}),$
\end{center}

\strut 

where $\mathcal{T}(H_{s})$ means the classical Toeplitz algebra defined on $%
H_{s}.$

\strut \strut \strut 

(4) If $u$ $=$ $0$ and if $s$ $\in $ $PI_{H_{s}}(0,$ $\infty ),$ then

\strut 

\begin{center}
$C^{*}(\{a\})$ $=$ $C^{*}(\{s\})$ $\overset{*\text{-isomorphic}}{=}$ $(\Bbb{C%
}\cdot 1_{H_{s}})$ $\otimes _{\Bbb{C}}$ $M_{2}(\Bbb{C}).$
\end{center}
\end{proposition}

\strut

\begin{proof}
(1) Let $a$ $\in $ $PI_{H}(\varepsilon _{+},$ $\varepsilon ^{-},$ $%
\varepsilon _{-}^{-})$ be a partial isometry on $H$ having its Wold
decomposition $a$ $=$ $u$ $+$ $s,$ where $u$ is the unitary part of $a$ on $%
H_{u}$ and $s$ is the shift part of $a$ on $H_{s},$ where the subspaces $%
H_{u}$ and $H_{s}$ are the unitary part and the shift part of $H$ in terms
of $a.$ Recall the diagram

\strut

\begin{center}
$
\begin{array}{lllllll}
H & = & H_{u} & \oplus & \,\,\,\,\,\,\,\,\,\,\,\,\,\,H_{s} & \oplus & \ker a
\\ 
\downarrow _{a} &  & \downarrow _{u} &  & \,\,\,\,\,\,\,\,\,\,\,\,\,%
\downarrow _{s} &  &  \\ 
H & = & H_{u} & \oplus & \ker s^{*}\oplus H_{s} & \oplus & (\ker
a^{*}\ominus \ker s^{*}).
\end{array}
$
\end{center}

\strut \strut

Notice that $\ker s^{*}$ $\oplus $ $H_{s}$ is Hilbert-space isomorphic to $%
H_{s},$ since $s$ is a shift. So, the $C^{*}$-subalgebra $C^{*}(\{a\})$
generated by $a$ is $*$-isomorphic to the $C^{*}$-algebra $C^{*}(\{u,$ $s\})$
generated by $u$ and $s.$ Moreover, since $H_{u}$ and $H_{s}$ are orthogonal
(as subspaces) in $H,$ the $C^{*}$-subalgebra $C^{*}(\{u,$ $s\})$ is $*$%
-isomorphic to $C^{*}(\{u\})$ $\oplus $ $C^{*}(\{s\})$ as an embedded $C^{*}$%
-subalgebra of $B(H).$

\strut

(2) Suppose the shift part $s$ of a partial isometry $a$ is the zero
operator. Then $a$ $=$ $u,$ and hence the $C^{*}$-subalgebra $C^{*}(\{a\})$
is $*$-isomorphic to $C^{*}(\{u\})$, as an embedded $C^{*}$-subalgebra of $%
B(H_{u})$ $\overset{C^{*}\text{-subalgebra}}{\subseteq }$ $B(H).$ By [13],
for the given unitary part $u$, we can create a one-vertex-one-edge graph $%
G. $ i.e., $\{u\}$ is a connected $G$-family in the sense of Section 2.3.
Therefore,

\strut

\begin{center}
$C^{*}(\{u\})$ $\overset{*\text{-isomorphic}}{=}$ $(\Bbb{C}$ $\cdot $ $%
1_{H_{u}})$ $\otimes _{\Bbb{C}}$ $C(\Bbb{T}).$
\end{center}

\strut

(3) Assume now that the unitary part $a$ is the zero operator. Then $a$ $=$ $%
s,$ and hence the $C^{*}$-subalgebra $C^{*}(\{a\})$ is $*$-isomorphic to $%
C^{*}(\{s\}),$ as an embedded $C^{*}$-subalgebra of $B(H_{s}).$ Suppose that 
$\varepsilon ^{-}$ $=$ $k_{s}$ $\in $ $\Bbb{N}.$ For example, let $k_{s}$ $=$
$1,$ and hence let $s$ $\in $ $PI_{H_{s}}(0,$ $1).$ Then this operator $s$
is unitarily equivalent to the unilateral shift $U$ on $l^{2}(\Bbb{N}_{0})$
which is Hilbert-space isomorphic to $H_{s}$. It is well-known that the
unilateral shift $U$ generates the Toeplitz algebra $\mathcal{T}(l^{2}(\Bbb{N%
}_{0})).$ So, the shift part $s$ generates $C^{*}$-algebra $*$-isomorphic to
the classical Toeplitz algebra $\mathcal{T}(H_{s}).$

\strut

Suppose $k_{s}$ $>$ $1$ in $\Bbb{N}.$ Then the shift part $s$ on $H_{s}$ is
unitarily equivalent to the operator $U^{k_{s}},$ where $U$ is the
unilateral shift on $l^{2}(\Bbb{N}_{0}).$ Note that the operator $U^{k_{s}}$
satisfies that

\strut

\begin{center}
$U^{k_{s}}\left( (\xi _{0},\text{ }\xi _{1},\text{ }\xi _{2},\text{ ...}%
)\right) $ $=$ $\left( \underset{k_{s}\text{-times}}{\underbrace{0,\text{
........, }0}}\text{ },\text{ }\xi _{0},\text{ }\xi _{1},\text{ ...}\right)
, $
\end{center}

\strut

for all $(\xi _{0},$ $\xi _{1},$ $\xi _{2},$ ...$)$ $\in $ $l^{2}(\Bbb{N}%
_{0}),$ with $\dim (U^{k_{s}})^{*}$ $=$ $k_{s}.$ So, the $C^{*}$-algebra $%
C^{*}(\{U^{k_{s}}\})$ is also $*$-isomorphic to the Toeplitz algebra $%
\mathcal{T}(l(\Bbb{N}_{0})).$ Therefore, if $k_{s}$ $<$ $\infty $ in $\Bbb{N}%
,$ then the $C^{*}$-subalgebra $C^{*}(\{s\})$ is $*$-isomorphic to the
classical Toeplitz algebra $\mathcal{T}(H_{s}),$ as an embedded $C^{*}$%
-subalgebra of $B(H_{s})$ $\overset{C^{*}\text{-subalgebra}}{\subseteq }$ $%
B(H).$

\strut \strut

(4) Suppose $s$ $\in $ $PI_{H_{s}}(0,$ $\infty ).$ Then the operator $s$ is
unitarily equivalent to the block operator matrix

\strut

\begin{center}
$\left( 
\begin{array}{ll}
0_{H_{s}} & 0_{H_{s}} \\ 
1_{H_{s}} & 0_{H_{s}}
\end{array}
\right) ,$
\end{center}

\strut

and hence the $C^{*}$-algebra $C^{*}(\{s\})$ is $*$-isomorphic to $(\Bbb{G}$ 
$\cdot $ $1_{H_{s}})$ $\otimes _{\Bbb{C}}$ $M_{2}(\Bbb{C}),$ by [13]. Notice
that, on the previous case, the Hilbert space $H_{s}$ is Hilbert-space
isomorphic to $H_{1}$ $\oplus $ $H_{2},$ where both $H_{1}$ and $H_{2}$ are
Hilbert-space isomorphic to $H_{s}.$
\end{proof}

\strut \strut \strut

Motivated by the statement (2) of the previous proposition, we now define a
new number $\varepsilon _{u}$ $\in $ $\Bbb{N}_{0}^{\infty }.$

\strut

\begin{definition}
Let $a$ $\in $ $PI_{H}(\varepsilon _{+},$ $\varepsilon ^{-},$ $\varepsilon
_{-}^{-}),$ with its Wold decomposition $a$ $=$ $u$ $+$ $s,$ where $u$ and $s
$ are the unitary part and the shift part of $a,$ respectively, and let $%
H_{u}$ be the unitary part of $H,$ in terms of $a.$ Define the number $%
\varepsilon _{u}$ by

\strut 

\begin{center}
$\varepsilon _{u}$ $\overset{def}{=}$ $\dim H_{u}$ $\in $ $\Bbb{N}%
_{0}^{\infty }.$
\end{center}

\strut 

And define the subclass $PI_{H}(\varepsilon _{0},$ $\varepsilon _{+},$ $%
\varepsilon ^{-},$ $\varepsilon _{-}^{-})$ of $PI_{H}(\varepsilon _{+},$ $%
\varepsilon ^{-},$ $\varepsilon _{-}^{-})$ by

\strut 

\begin{center}
$PI_{H}(\varepsilon _{0},$ $\varepsilon _{+},$ $\varepsilon ^{-},$ $%
\varepsilon _{-}^{-})$ $\overset{def}{=}$ $\{a$ $\in $ $PI_{H}(\varepsilon
_{+},$ $\varepsilon ^{-},$ $\varepsilon _{-}^{-})$ $:$ $\dim H_{u_{a}}$ $=$ $%
\varepsilon _{0}\},$
\end{center}

\strut 

where $H_{u_{a}}$ is the unitary part of $H,$ in terms of $a.$
\end{definition}

\strut \strut \strut

Since $PI(H)$ $=$ $\underset{(\varepsilon _{+},\text{ }\varepsilon _{-})\in 
\Bbb{N}_{0}^{\infty }}{\cup }$ $PI_{H}(\varepsilon _{+},$ $\varepsilon
_{-}), $ we can get that

\strut

\begin{center}
$PI(H)$ $=$ $\underset{(\varepsilon _{0},\text{ }\varepsilon _{+},\text{ }%
\varepsilon ^{-},\text{ }\varepsilon _{-}^{-})\in (\Bbb{N}_{0}^{\infty })^{4}%
}{\cup }$ $\left( PI_{H}(\varepsilon _{0},\text{ }\varepsilon _{+},\text{ }%
\varepsilon ^{-},\text{ }\varepsilon _{-}^{-})\right) .$
\end{center}

\strut

Define an equivalence relation $\mathcal{R}_{*}$ on $PI(H)$ by

\strut

\strut (3.1)

\begin{center}
$a_{1}$ $\mathcal{R}_{*}$ $a_{2}$ $\overset{def}{\Longleftrightarrow }$ $%
C^{*}(\{a_{1}\})$ $\overset{*\text{-isomorphic}}{=}$ $C^{*}(\{a_{2}\})$ in $%
B(H).$
\end{center}

\strut \strut \strut \strut

By the previous proposition, if $s_{1}$ $\in $ $PI_{H}(0,$ $k_{1})$ and $%
s_{2}$ $\in $ $PI_{H}(0,$ $k_{2})$ are shifts (which are isometries on $H$),
where $k_{1},$ $k_{2}$ $<$ $\infty $ in $\Bbb{N},$ then we can get that

$\strut $

\begin{center}
$
\begin{array}{ll}
C^{*}(\{s_{1}\}) & \overset{*\text{-isomorphic}}{=}C^{*}(\{U^{k_{1}}\}) \\ 
&  \\ 
& \overset{*\text{-isomorphic}}{=}\mathcal{T}(H) \\ 
&  \\ 
& \overset{*\text{-isomorphic}}{=}C^{*}(\{U^{k_{2}}\})\overset{*\text{%
-isomorphic}}{=}C^{*}(\{s_{2}\}),
\end{array}
$
\end{center}

\strut

where $\mathcal{T}(H)$ is the classical Toeplitz algebra.

\strut

\begin{definition}
Let $a$ $\in $ $PI_{H}(\varepsilon _{0},$ $\varepsilon _{+},$ $\varepsilon
^{-},$ $\varepsilon _{-}^{-}).$ The quadruple $(\varepsilon _{0},$ $%
\varepsilon _{+},$ $\varepsilon ^{-},$ $\varepsilon _{-}^{-})$ $\in $ $(\Bbb{%
N}_{0}^{\infty })^{4}$ is called the $*$-isomorphic index of $a$. The $*$%
-isomorphic index of $a$ is denoted by $i_{*}(a).$ i.e., $i_{*}(a)$ $=$ $%
(\varepsilon _{0},$ $\varepsilon _{+},$ $\varepsilon ^{-},$ $\varepsilon
_{-}^{-}).$
\end{definition}

\strut

Consider the set

\strut

\begin{center}
$(\Bbb{N}_{0}^{\infty })^{4}$ $=$ $\Bbb{N}_{0}^{\infty }$ $\times $ $\Bbb{N}%
_{0}^{\infty }$ $\times $ $\Bbb{N}_{0}^{\infty }$ $\times $ $\Bbb{N}%
_{0}^{\infty },$
\end{center}

\strut

in detail. We will define a binary operation $(-)$ on $(\Bbb{N}_{0}^{\infty
})^{4}$ by

\strut

$\qquad (i_{1},$ $j_{1},$ $k_{1},$ $l_{1})$ $-$ $(i_{2},$ $j_{2},$ $k_{2},$ $%
l_{2})$

$\strut $

\begin{center}
$\overset{def}{=}$ $\left( \left| i_{1}-i_{2}\right| ,\text{ }\left|
j_{1}-j_{2}\right| ,\text{ }\left| k_{1}-k_{2}\right| ,\text{ }\left|
l_{1}-l_{2}\right| \right) $
\end{center}

\strut \strut \strut

under the additional rule:

\strut

\begin{center}
\strut $\infty -\infty \overset{def}{=}$ $0$ in $\Bbb{N}_{0}^{\infty },$
\end{center}

\strut

for all $(i_{t},$ $j_{t},$ $k_{t},$ $l_{t})$ $\in $ $(\Bbb{N}_{0}^{\infty
})^{4},$ for $t$ $=$ $1,$ $2.$ From now, if we denote $(\Bbb{N}_{0}^{\infty
})^{4},$ then it means the algebraic pair $((\Bbb{N}_{0}^{\infty })^{4},$ $%
-).$ Then, we can get the following theorem.

\strut

\begin{theorem}
Let $a_{j}$ $\in $ $PI(H)$ be partial isometries, and suppose the spectra $%
spec(u_{j})$ of the unitary parts $u_{j}$ of $a_{j}$ are identical (if they
are nonzero), for $j$ $=$ $1,$ $2.$ Then we can have that

\strut 

(3.2)$\qquad \qquad \qquad \qquad \quad i_{*}(a_{1})$ $-$ $i_{*}(a_{2})$ $=$ 
$(0,$ $k_{1},$ $k_{2},$ $0)$

\strut 

for some $k_{1},$ $k_{2}$ $<$ $\infty $ in $\Bbb{N}_{0},$ if and only if $%
a_{1}$ $\mathcal{R}_{*}$ $a_{2},$ equivalently, the $C^{*}$-algebras $%
C^{*}(\{a_{1}\})$ and $C^{*}(\{a_{2}\})$ are $*$-isomorphic, as embedded $%
C^{*}$-subalgebras of $B(H).$
\end{theorem}

\strut

\begin{proof}
Suppose the $*$-isomorphic indices $i_{*}(a_{k})$ are $(t_{1}^{(k)},$ $%
t_{2}^{(k)},$ $t_{3}^{(k)},$ $t_{4}^{(k)})$ in $(\Bbb{N}_{0}^{\infty })^{4},$
for $k$ $=$ $1,$ $2.$ Equivalently, assume that $a_{k}$ $\in $ $%
PI_{H}(t_{1}^{(k)},$ $t_{2}^{(k)},$ $t_{3}^{(k)},$ $t_{4}^{(k)}),$ for $k$ $%
= $ $1,$ $2.$ Suppose $a_{k}$ have their Wold decomposition $u_{k}$ $+$ $%
s_{k}, $ where $u_{k}$ are the unitary parts of $a_{k}$ and $s_{k}$ are the
shift parts of $a_{k},$ for $k$ $=$ $1,$ $2.$ Also, let $H_{u_{k}}$ and $%
H_{s_{k}}$ be the unitary parts of $H$ and the shift parts of $H,$ in terms
of $a_{k},$ for $k$ $=$ $1,$ $2,$ respectively.

\strut

($\Rightarrow $) Suppose the $*$-isomorphic indices $i_{*}(a_{1})$ and $%
i_{*}(a_{2})$ satisfies (3.2). Then the condition (3.2) says that

\strut

\begin{center}
$t_{1}^{(1)}$ $=$ $\dim H_{u_{1}}$ $=$ $t_{1}$ $=$ $\dim H_{u_{2}}$ $=$ $%
t_{1}^{(2)},$
\end{center}

\strut

\begin{center}
$t_{2}^{(1)}$ $=$ $\dim (\ker a_{1})$ $=$ $t_{2}$ $=$ $\dim (\ker a_{2})$ $=$
$t_{2}^{(2)},$
\end{center}

\strut

\begin{center}
$t_{4}^{(1)}$ $=$ $\dim (\ker a_{1}^{*}\ominus \ker s_{1}^{*})$ $=$ $t_{4}$ $%
=$ $\dim (\ker a_{2}^{*}\ominus \ker s_{2}^{*})$ $=$ $t_{4}^{(2)}$
\end{center}

and

\begin{center}
$\left| t_{2}^{(1)}-t_{2}^{(2)}\right| =k_{1},$ $\left|
t_{3}^{(1)}-t_{3}^{(2)}\right| =k_{2}$ in $\Bbb{N}_{0}.$
\end{center}

\strut

Assume now that $k_{1}$ $<$ $\infty $ and $k_{2}$ $<$ $\infty .$ Then

\strut

\begin{center}
$a_{1}$ $\in $ $PI_{H}(t_{1},$ $t_{2}^{(1)},$ $t_{3}^{(1)},$ $t_{4})$ and $%
a_{2}$ $\in $ $(t_{1},$ $t_{2}^{(1)},$ $t_{3}^{(2)},$ $t_{4}).$
\end{center}

\strut

Since $C^{*}(\{a_{j}\})$ $\overset{*\text{-isomorphic}}{=}$ $%
C^{*}(\{u_{j}\}) $ $\oplus $ $C^{*}(\{s_{j}\}),$ for $j$ $=$ $1,$ $2,$ it
suffices to show that $C^{*}(\{s_{1}\})$ and $C^{*}(\{s_{2}\})$ are $*$%
-isomorphic. By the previous proposition, $C^{*}(\{s_{k}\})$ are $*$%
-isomorphic either $(\Bbb{C}$ $\cdot $ $1_{H_{k}})$ $\otimes _{\Bbb{C}}$ $%
M_{2}(\Bbb{C})$ or $\mathcal{T}(H_{s_{k}}),$ where $\mathcal{T}(H_{s_{k}})$
is the Toeplitz algebra on $H_{s_{k}},$ for $k$ $=$ $1,$ $2.$ By the
assumption that

\strut

\begin{center}
$k_{1}$ $=$ $\left| t_{2}^{(1)}-t_{2}^{(2)}\right| $ $<$ $\infty $ and $%
k_{2} $ $=$ $\left| t_{3}^{(1)}\text{ }-\text{ }t_{3}^{(2)}\right| $ $<$ $0,$
\end{center}

\strut

either [$t_{3}^{(1)}$ $<$ $\infty $ and $t_{3}^{(2)}$ $<$ $\infty $] or [$%
t_{3}^{(1)}$ $=$ $\infty $ $=$ $t_{3}^{(2)}$]. If $t_{3}^{(1)}$ $<$ $\infty $
and $t_{3}^{(2)}$ $<$ $\infty ,$ then

\strut

\begin{center}
$C^{*}(\{s_{1}\})$ $\overset{*\text{-isomorphic}}{=}$ $\mathcal{T}(l^{2}(%
\Bbb{N}_{0}))$ $\overset{*\text{-isomorphic}}{=}$ $C^{*}(\{s_{2}\}),$
\end{center}

\strut

by the previous proposition. If $t_{3}^{(1)}$ $=$ $t_{3}^{(2)},$ then

\strut

\begin{center}
$C^{*}(\{s_{1}\})$ $\overset{*\text{-isomorphic}}{=}$ $(\Bbb{C}\cdot
1_{H_{0}})$ $\otimes _{\Bbb{C}}$ $M_{2}(\Bbb{C})$ $\overset{*\text{%
-isomorphic}}{=}$ $C^{*}(\{s_{2}\}),$
\end{center}

\strut

again by the previous proposition. Therefore, if (3.2) holds true, then $%
C^{*}(\{s_{1}\})$ and $C^{*}(\{s_{2}\})$ are $*$-isomorphic, and hence, $%
C^{*}(\{a_{1}\})$ and $C^{*}(\{a_{2}\})$ are $*$-isomorphic, as embedded $%
C^{*}$-subalgebras of $B(H).$ i.e., $a_{1}$ $\mathcal{R}_{*}$ $a_{2}.$

\strut

($\Leftarrow $) Assume now that two partial isometries $a_{1},$ $a_{2}$ $\in 
$ $PI(H)$ satisfy $a_{1}$ $\mathcal{R}_{*}$ $a_{2},$ equivalently, the $%
C^{*} $-subalgebras $C^{*}(\{a_{1}\})$ and $C^{*}(\{a_{2}\})$ are $*$%
-isomorphic in $B(H).$ This means that $C^{*}(\{u_{1}\})$ and $%
C^{*}(\{u_{2}\})$ (resp., $C^{*}(\{s_{1}\})$ and $C^{*}(\{s_{2}\})$) are $*$%
-isomorphic, by the previous proposition. It is clear that $C^{*}(\{u_{1}\})$
and $C^{*}(\{u_{2}\})$ are $*$-isomorphic, as embedded $C^{*}$-subalgebras
of $B(H),$ if and only if $H_{u_{1}}$ is Hilbert-space isomorphic to $%
H_{u_{2}}$ if and only if $t_{1}^{(1)}$ $=$ $t_{1}^{(2)},$ whenever $%
spec(u_{1})$ $=$ $spec(u_{2})$ in $\Bbb{T}$ $\subset $ $\Bbb{C}.$ Clearly,
as partial isometries, if $a_{1}$ $\mathcal{R}_{*}$ $a_{2},$ then $%
t_{4}^{(1)}$ $=$ $t_{4}^{(2)}.$

\strut

Assume that both $C^{*}(\{s_{1}\})$ and $C^{*}(\{s_{2}\})$ are $*$%
-isomorphic to the classical Toeplitz algebra $\mathcal{T}(l^{2}(\Bbb{N}%
_{0})).$ Then $t_{2}^{(1)},$ $t_{2}^{(2)},$ $t_{3}^{(1)},$ $t_{3}^{(2)}$ $<$ 
$\infty $ in $\Bbb{N}_{0}^{\infty }.$ Thus

\strut

\begin{center}
$\left| t_{2}^{(1)}-t_{2}^{(2)}\right| =k_{1}$ $<$ $\infty $ and $\left|
t_{3}^{(1)}-t_{3}^{(2)}\right| $ $=$ $k_{2}$ $<$ $\infty .$
\end{center}

\strut

Assume now that both $C^{*}(\{s_{1}\})$ and $C^{*}(\{s_{2}\})$ are $*$%
-isomorphic to the $C^{*}$-algebra $(\Bbb{C}$ $\cdot $ $1_{K})$ $\otimes _{%
\Bbb{C}}$ $M_{2}(\Bbb{C}).$ Then $t_{3}^{(1)}$ $=$ $\infty $ $=$ $%
t_{3}^{(2)} $ in $\Bbb{N}_{0}^{\infty }.$ Thus

\strut

\begin{center}
$\left| t_{3}^{(1)}-t_{3}^{(2)}\right| $ $=$ $0$ $=$ $\left|
t_{2}^{(1)}-t_{2}^{(2)}\right| ,$
\end{center}

\strut

by the rule $\infty $ $-$ $\infty $ $\overset{def}{=}$ $0$ in $\Bbb{N}%
_{0}^{\infty }.$ This shows that if $a_{1}$ $\mathcal{R}_{*}$ $a_{2}$ and if 
$i_{*}(a_{k})$ $=$ $(t_{1}^{(k)},$ $t_{2}^{(k)},$ $t_{3}^{(k)},$ $%
t_{4}^{(k)})$ in $(\Bbb{N}_{0}^{\infty })^{4},$ then

\strut

\begin{center}
$i_{*}(a_{1})$ $-$ $i_{*}(a_{2})$ $=$ $(0,$ $k_{1},$ $k_{2},$ $0),$ with $%
k_{1},$ $k_{2}$ $<$ $\infty .$
\end{center}

\strut
\end{proof}

\strut

The above theorem shows that the $*$-isomorphic indexing

$\strut $

\begin{center}
$i_{*}(\cdot )$ $:$ $PI(H)$ $\rightarrow $ $(\Bbb{N}_{0}^{\infty })^{4}$
\end{center}

\strut

is \emph{a kind of} invariant on $*$-isomorphic $C^{*}$-subalgebras
generated by a single partial isometry in $B(H),$ up to the spectrums of
unitary parts and the ``finite'' difference of the second and the third
entries of $i_{*}(\cdot ).$

\strut

\begin{example}
Let $H$ $=$ $l^{2}(\Bbb{N}_{0})$, and let $a$ $\in $ $B(H)$ be an operator
which is unitarily equivalent to the following matrix

\strut 

\begin{center}
$\left( 
\begin{array}{lll}
\left( 
\begin{array}{lll}
\theta  &  &  \\ 
& \ddots  &  \\ 
&  & \theta 
\end{array}
\right)  &  & \qquad \qquad 0 \\ 
& \left( 
\begin{array}{lllll}
0 &  &  &  &  \\ 
1 & 0 &  &  &  \\ 
& 1 & 0 &  &  \\ 
&  & 1 & 0 &  \\ 
&  &  & \ddots  & \ddots 
\end{array}
\right)  &  \\ 
0 &  & \left( 
\begin{array}{lll}
0 &  &  \\ 
& 0 &  \\ 
&  & \ddots 
\end{array}
\right) 
\end{array}
\right) ,$
\end{center}

\strut 

where $\theta $ $\in $ $\Bbb{T}$ in $\Bbb{C}$. Without loss of generality,
let $a$ be the operator having above matrix form on $H.$ Then it is a
partial isometry on $H.$ Assume that the left cornered block matrix

\strut 

\begin{center}
$\left( 
\begin{array}{lll}
\theta  &  &  \\ 
& \ddots  &  \\ 
&  & \theta 
\end{array}
\right) $ $\overset{denote}{=}$ $a_{u}$
\end{center}

\strut 

of $a$ is contained in the matricial algebra $M_{n}(\Bbb{C}),$ with $%
spec(a_{u})$ $=$ $\{\theta \},$ and the middle block matrix

\strut \newline

\begin{center}
$\left( 
\begin{array}{llll}
0 &  &  &  \\ 
1 & 0 &  &  \\ 
& 1 & 0 &  \\ 
&  & 1 & \ddots  \\ 
&  &  & \ddots 
\end{array}
\right) $ $\overset{denote}{=}$ $a_{s}$
\end{center}

\strut 

is the unilateral shift. Then this partial isometry $a$ has its Wold
decomposition $a$ $=$ $a_{u}$ $+$ $a_{s},$ and the unitary part $H_{u_{a}}$
of $H$ is Hilbert-space isomorphic to $\Bbb{C}^{\oplus \,n}$ and the shift
part $H_{a_{s}}$ of $H$ is Hilbert-space isomorphic to $l^{2}(\Bbb{N}_{0}).$
So, we can get that this partial isometry $a$ is contained in $PI_{H}(n,$ $%
\infty ,$ $1,$ $\infty ),$ and hence the $*$-isomorphic index $i_{*}(a)$ is

\strut 

\begin{center}
$i_{*}(a)$ $=$ $(n,$ $\infty ,$ $1,$ $\infty ).$
\end{center}

\strut 

Also, the $C^{*}$-subalgebra $C^{*}(\{a\})$ generated by $a$ is $*$%
-isomorphic to

\strut 

\begin{center}
$
\begin{array}{ll}
\Bbb{C}^{\oplus \,n} & \oplus \text{ }\mathcal{T}\left( l^{2}(\Bbb{N}%
_{0})\right)  \\ 
&  \\ 
& \overset{C^{*}\text{-subalgebra}}{\subseteq }B(H_{a_{u}})\oplus
B(H_{a_{s}}) \\ 
&  \\ 
& \overset{C^{*}\text{-subalgebra}}{\subseteq }B(H).
\end{array}
$
\end{center}
\end{example}

\strut

\strut \strut \strut \strut

\strut

\subsection{Graph Groupoids Induced by Partial Isometries}

\strut

\strut

\strut

Throughout this section, we will use the same notations we used in the
previous sections. In Section 3.1, we characterized $C^{*}$-subalgebras of $%
B(H)$ generated by a single partial isometry. In the rest of this Section,
we will observe $C^{*}$-subalgebras of $B(H)$ generated by multi-partial
isometries on $H.$ So, let $\mathcal{G}$ $=$ $\{a_{1},$ ..., $a_{N}\}$ be a
family of partial isometries $a_{j}$'s on $H,$ for all $j$ $=$ $1,$ ..., $N.$
And let's assume that $N$ $>$ $1.$ Also, let

\strut

\begin{center}
$i_{*}(a_{j})$ $=$ $(k_{1}^{(j)},$ $k_{2}^{(j)},$ $k_{3}^{(j)},$ $%
k_{4}^{(j)}),$ for all $j$ $=$ $1,$ ..., $N.$
\end{center}

\strut i.e.,

\begin{center}
$k_{1}^{(j)}$ $=$ $\dim $ $H_{u_{j}},$ where $u_{j}$ is the unitary part of $%
a_{j},$
\end{center}

\strut

\begin{center}
$k_{2}^{(j)}$ $=$ $\dim \left( \ker a_{j}^{*}\right) ,$\quad $k_{3}^{(j)}$ $%
= $ $\dim \left( \ker s_{j}^{*}\right) $
\end{center}

and

\begin{center}
$k_{4}^{(j)}$ $=$ $\dim \left( \ker a_{j}^{*}\ominus \ker s_{j}^{*}\right) ,$
\end{center}

\strut

for all $j$ $=$ $1,$ ..., $N,$ where $a_{j}$ have their Wold decomposition $%
a_{j}$ $=$ $u_{j}$ $+$ $s_{j},$ for $j$ $=$ $1,$ ..., $N.$ Then we can
construct a family $\mathcal{G}_{W}$ induced by $\mathcal{G}$:

\strut

\begin{center}
$\mathcal{G}_{W}$ $=$ $\mathcal{G}^{(u)}$ $\cup $ $\mathcal{G}^{(s)},$
\end{center}

where

\begin{center}
$\mathcal{G}^{(u)}$ $\overset{def}{=}$ $\{u_{j}$ $:$ $j$ $=$ $1,$ ..., $N\}$
\end{center}

and

\begin{center}
$\mathcal{G}^{(s)}$ $\overset{def}{=}$ $\{s_{j}$ $:$ $j$ $=$ $1,$ ..., $N\}.$
\end{center}

\strut

Also, for our purpose, we decompose $\mathcal{G}^{(s)}$ by

\strut

\begin{center}
$\mathcal{G}^{(s)}$ $=$ $\mathcal{G}_{f}^{(s)}$ $\cup $ $\mathcal{G}_{\infty
}^{(s)},$
\end{center}

where

\begin{center}
\strut $\mathcal{G}_{f}^{(s)}$ $\overset{def}{=}$ $\left\{ x\in \mathcal{G}%
^{(s)}\left| 
\begin{array}{c}
i_{*}(x)=(0,\varepsilon _{+},\varepsilon ^{-},\varepsilon _{-}^{-}), \\ 
\varepsilon _{+},\text{ }\varepsilon _{-}^{-}\in \Bbb{N}_{0}^{\infty },\text{
and} \\ 
\varepsilon ^{-}<\infty \text{ in }\Bbb{N}_{0}
\end{array}
\right. \right\} $
\end{center}

and

\begin{center}
$\mathcal{G}_{\infty }^{(s)}$ $\overset{def}{=}$ $\left\{ x\in \mathcal{G}%
^{(s)}\left| 
\begin{array}{c}
i_{*}(x)=(0,\text{ }\varepsilon _{+},\text{ }\varepsilon ^{-},\text{ }%
\varepsilon _{-}^{-}) \\ 
\varepsilon _{+},\text{ }\varepsilon _{-}^{-}\text{ }\in \Bbb{N}_{0}^{\infty
}\text{, and} \\ 
\varepsilon ^{-}=\infty
\end{array}
\right. \right\} .$
\end{center}

\strut

Define a subset $\widehat{\mathcal{G}_{W}}$ of $B(H)$ by $\mathcal{G}_{W}$ $%
\cup $ $\mathcal{G}_{W}^{*}.$ Also, define the set $\mathcal{G}_{pro}$
induced by $\mathcal{G}_{W}$ by

\strut

\begin{center}
$\mathcal{G}_{pro}$ $\overset{def}{=}$ $\{x^{*}$ $x,$ $x$ $x^{*}$ $:$ $x$ $%
\in $ $\mathcal{G}_{W}\}$ $\subset $ $B(H)_{pro},$
\end{center}

\strut

where $B(H)_{pro}$ is the subset of $B(H)$ consisting of all projections on $%
H.$

\strut

Consider now the subset $\Bbb{G}_{0}^{0}$ $\subset $ $B(H),$ consisting of
words in $\{x^{k\,*}$ $x^{k}$ $:$ $x$ $\in $ $\widehat{\mathcal{G}_{W}},$ $k$
$\in $ $\Bbb{N}\}.$ i.e.,

\strut

\begin{center}
$
\begin{array}{ll}
\Bbb{G}_{0}^{0} & \overset{def}{=}\{0\}\cup \mathcal{G}_{pro}\cup \left( 
\underset{k=2}{\overset{\infty }{\cup }}\text{ }\{(x^{k})(x^{k})^{*},\text{ }%
(x^{k})^{*}(x^{k}):x\in \mathcal{G}_{W}\}\right) \\ 
&  \\ 
& =\{0\}\cup \left( \underset{k=1}{\overset{\infty }{\cup }}\left( \left\{
(x^{k})(x^{k})^{*},\text{ }(x^{k})^{*}(x^{k}):x\in \mathcal{G}_{W}\right\}
\right) \right)
\end{array}
$
\end{center}

\strut

Notice that, if $x$ $\in $ $\widehat{\mathcal{G}^{(u)}}$ $\subset $ $%
\widehat{\mathcal{G}_{W}},$ where $\widehat{\mathcal{G}^{(u)}}$ $=$ $%
\mathcal{G}^{(u)}$ $\cup $ $\mathcal{G}^{(u)\,*},$ then

\strut

\begin{center}
$x^{*}$ $x$ $=$ $x$ $x^{*}$ $=$ $1_{H_{u}}$ $=$ $(x^{n})$ $(x^{n})^{*}$ $=$ $%
(x^{n})^{*}$ $(x^{n}),$ for all $n$ $\in $ $\Bbb{N}.$
\end{center}

\strut

If $x$ $\in $ $\widehat{\mathcal{G}_{\infty }^{(s)}}$ $=$ $\mathcal{G}%
_{\infty }^{(s)}$ $\cup $ $\mathcal{G}_{\infty }^{(s)\,*},$ then

\strut

\begin{center}
$x^{n}$ $=$ $0$ $=$ $(x^{n})^{*}$ $=$ $(x^{*})^{n},$ for all $n$ $\in $ $%
\Bbb{N}$ $\setminus $ $\{1\}.$
\end{center}

\strut

Suppose now that $x$ $\in $ $\widehat{\mathcal{G}_{f}^{(s)}}$ $=$ $\mathcal{G%
}_{f}^{(s)}$ $\cup $ $\mathcal{G}_{f}^{(s)\,*}.$ Then we can have that

\strut

\begin{center}
$(x^{n_{1}})$ $(x^{n_{1}})^{*}$ $\neq $ $(x^{n_{2}})(x^{n_{2}})^{*}$ in $%
B(H)_{pro},$ for all $n_{1}$ $\neq $ $n_{2}$ $\in $ $\Bbb{N}.$
\end{center}

\strut

\begin{remark}
Suppose $x$ $\in $ $\mathcal{G}_{f}^{(s)}.$ Then, $x^{n}$ $\in $ $B(H)$ is a
shift, too. Indeed, if $i_{*}(x)$ $=$ $(0,$ $\varepsilon _{+},$ $k,$ $%
\varepsilon _{-}^{-}),$ with $k$ $<$ $\infty $ in $\Bbb{N}_{0},$ then $%
i_{*}(x^{n})$ $=$ $(0,$ $\varepsilon _{+},$ $nk,$ $\varepsilon _{-}^{-})$ in 
$(\Bbb{N}_{0}^{\infty })^{4}.$ This shows that the operators $(x^{n})$ $%
(x^{n})^{*}$ are projections on $H,$ too, since $x^{n}$ are partial
isometries on $H,$ for all $n$ $\in $ $\Bbb{N}.$
\end{remark}

\strut

By the previous remark, we can conclude that:

\strut

\begin{lemma}
Let $p$ $\in $ $\Bbb{G}_{0}^{0}.$ Then $p$ is a projection on $H.$
\end{lemma}

\strut

\begin{proof}
Let $p$ $\in $ $\Bbb{G}_{0}^{0}.$ Then there exists $x$ $\in $ $\widehat{%
\mathcal{G}_{W}}$ and $m$ $\in $ $\Bbb{N},$ such that $p$ $=$ $(x^{m})$ $%
(x^{m})^{*}.$ Assume that $x$ $\in $ $\widehat{\mathcal{G}^{(u)}}.$ Then

\strut

\begin{center}
$(x^{m})$ $(x^{m})^{*}$ $=$ $1_{H_{x}}$ $=$ $x$ $x^{*}$ $=$ $x^{*}$ $x,$
\end{center}

\strut

for all $m$ $\in $ $\Bbb{N}$, where $H_{x}$ $=$ $(x$ $x^{*})$ $H$ is the
subspace of $H,$ which is both the initial and the final spaces of $x.$
Suppose that $x$ $\in $ $\widehat{\mathcal{G}_{\infty }^{(s)}}.$ Then

\strut

\begin{center}
$(x^{m})(x^{m})^{*}$ $=$ $\left\{ 
\begin{array}{lll}
x\text{ }x^{*} &  & \text{if }m=1 \\ 
0 &  & \text{otherwise,}
\end{array}
\right. $
\end{center}

\strut

since $x^{m}$ $=$ $0$ $=$ $x^{*\,m},$ whenever $m$ $>$ $1.$ Therefore, it is
a projection on $H.$ Assume now that $x$ $\in $ $\widehat{\mathcal{G}%
_{f}^{(s)}}.$ Then, by the previous remark, the operator $(x^{m})$ $%
(x^{m})^{*}$ is a projection on $H,$ too. This shows that $\Bbb{G}_{0}^{0}$ $%
\subset $ $B(H)_{pro}.$
\end{proof}

\strut

Now, define the partial ordering $\leq $ on $\Bbb{G}_{0}^{0}$ by the rule:

\strut

\begin{center}
$p$ $\leq $ $q$ $\overset{def}{\Longleftrightarrow }$ $p$ $H$ $\overset{%
\text{Subspace}}{\subseteq }$ $q$ $H,$ for all $p,$ $q$ $\in $ $\Bbb{G}%
_{0}^{0}.$
\end{center}

\strut

The above partial ordering $\leq $ is well-defined, since all elements of $%
\Bbb{G}_{0}^{0}$ are projections on $H.$

\strut

\textbf{Notation} Without loss of generality, we will denote the partially
ordered set $(\Bbb{G}_{0}^{0},$ $\leq )$ in $B(H)_{pro}$ simply by $\Bbb{G}%
_{0}^{0}.$ $\square $

\strut

\strut

\subsubsection{Construction of Corresponding Graphs of Partial Isometries}

\strut

\strut \strut \strut \strut

Let $x$ $\in $ $\mathcal{G}_{W}.$ Then we can construct a corresponding
directed graph $G_{x}$ as follows:

\strut

(Case I)\ \ \ Suppose $u$ $\in $ $\mathcal{G}^{(u)}.$ Then construct the
one-vertex-one-loop-edge graph $G_{u}$ by a directed graph with

\strut

\begin{center}
$V(G_{u})$ $=$ $\{u^{*}$ $u$ $=$ $u$ $u^{*}\}$ and $E(G_{u})$ $=$ $\{u\}.$
\end{center}

\strut

i.e., we regard $u$ as a loop-edge and we regard the projections $u^{*}$ $u$
and $u$ $u^{*}$ as the initial and the terminal vertices of $u$:

\strut

\begin{center}
$G_{u}$ $=\qquad \overset{\overset{u^{*}u=uu^{*}}{\bullet }}{\underset{u}{%
\circlearrowleft }}$
\end{center}

\strut

(Case II)\ \ Suppose $s$ $\in $ $\mathcal{G}_{\infty }^{(s)}.$ Then
construct the one-edge graph $G_{s}$ by a directed graph with

\strut

\begin{center}
$V(G_{s})$ $=$ $\{s^{*}$ $s,$ $s$ $s^{*}\}$ and $E(G_{s})$ $=$ $\{s\}.$
\end{center}

\strut

i.e., we regard $s,$ $s^{*}$ $s$ and $s$ $s^{*}$ as the non-loop edge and
the corresponding initial and terminal vertices of the edge, respectively:

\strut

\begin{center}
$G_{s}$ $=$ \quad $\underset{s*s}{\bullet }\overset{s}{\longrightarrow }%
\underset{ss^{*}}{\bullet }$
\end{center}

\strut

(Case III) Assume now that $x$ $\in $ $\mathcal{G}_{f}^{(s)}.$ Then
construct the directed graph $G_{x}$ by the countable directed linear graph
with

\strut

\begin{center}
$V(G_{x})$ $=$ $\{x^{*}$ $x\}$ $\cup $ $\{(x^{n})$ $(x^{n})^{*}$ $:$ $n$ $%
\in $ $\Bbb{N}\}$
\end{center}

and

\begin{center}
$E(G_{x})$ $=$ $\{x^{(1)}$ $=$ $x\}$ $\cup $ $\{x^{(k)}$ $:$ $k$ $\in $ $%
\Bbb{N}$ $\setminus $ $\{1\}\}.$
\end{center}

\strut

i.e., we can have the following graph $G_{x}$:

\strut

\begin{center}
$G_{x}$ $=\quad \underset{x^{*}x}{\bullet }\overset{x}{\longrightarrow }%
\underset{xx^{*}}{\bullet }\overset{x^{(2)}}{\longrightarrow }\underset{%
x^{2}x^{2\,*}}{\bullet }\overset{x^{(3)}}{\longrightarrow }\underset{%
x^{3}x^{3\,*}}{\bullet }\overset{x^{(4)}}{\longrightarrow }\cdots $
\end{center}

\strut

Notice that $x^{(n)}$'s satisfy that

$\strut $

\begin{center}
$x^{(1)}$ $=$ $x$ and $x^{(n+1)}$ $=$ $x$ $=$ $x\mid _{(x^{n}x^{n\,*})H},$
for all $n$ $\in $ $\Bbb{N}$
\end{center}

\strut and

\begin{center}
$x^{(1)}$ $x^{(2)}$ ... $x^{(m)}$ $=$ $x^{m}$ $\in $ $B(H),$ for all $m$ $%
\in $ $\Bbb{N}$ $\setminus $ $\{1\},$
\end{center}

\strut

which are shifts on $H.$ Moreover, each vertex $(x^{n})$ $(x^{n})^{*}$ are
the projection onto the final space of the shift $x^{n},$ acting as the
terminal vertex of the edge $x^{(n)},$ for all $n$ $\in $ $\Bbb{N}$ $%
\setminus $ $\{1\}.$

\strut

\begin{definition}
Let $G_{x}$ be the directed graphs induced by $x$ $\in $ $\mathcal{G}_{W},$
by the above rules. Denote the collection of countable directed graphs $G_{x}
$'s induced by $x$'s in $\mathcal{G}_{W}$ by $\mathcal{G}_{G}.$ The
collection $\mathcal{G}_{G}$ is called the (collection of) corresponding
graphs induced by $\mathcal{G}_{W}.$
\end{definition}

\strut

\strut

\strut

\subsubsection{The Glued Graph $G_{1}$ $^{v_{1}}\#^{v_{2}}$ $G_{2}$ of $%
G_{1} $ and $G_{2}$}

\strut

\strut

Let $G_{k}$ be countable directed graphs and let $v_{k}$ $\in $ $V(G_{k})$
be the fixed vertices, for $k$ $=$ $1,$ $2.$ Then, by identifying the fixed
two vertices $v_{1}$ and $v_{2},$ we can construct a new countable directed
graph $G$ $=$ $G_{1}$ $^{v_{1}}\#^{v_{2}}$ $G_{2},$ by a countable directed
graph

\strut

\begin{center}
$V(G)$ $=$ $\left( V(G_{1})\text{ }\setminus \text{ }\{v_{1}\}\right) \cup $ 
$\{v_{12}\}$ $\cup $ $\left( V(G_{2})\text{ }\setminus \text{ }%
\{v_{2}\}\right) $
\end{center}

and

\begin{center}
$E(G)$ $=$ $E(G_{1})$ $\cup $ $E(G_{2}),$
\end{center}

\strut

where $v_{12}$ is the identified vertex of the vertices $v_{1}$ and $v_{2}$
(after identifying $v_{1}$ and $v_{2}$). If $e_{k}$ $\in $ $E(G_{k})$ and if 
$e_{k}$ $=$ $v_{k}$ $e_{k}$ or $e_{k}$ $=$ $e_{k}$ $v_{k},$ then this edge $%
e_{k}$ is regarded as $e_{k}$ $=$ $v_{12}$ $e_{k}$, respectively, $e_{k}$ $=$
$e_{k}$ $v_{12}$ in $E(G).$ For instance, let

\strut

\begin{center}
$G_{1}$ $=$\quad $\bullet \longrightarrow $ $\overset{v_{1}}{\bullet }$ $%
\longrightarrow \bullet $
\end{center}

and

\begin{center}
$G_{2}$ $=$ \quad $\bullet \longrightarrow $ $\overset{v_{2}}{\bullet }.$
\end{center}

\strut

Then the graph $G$ $=$ $G_{1}$ $^{v_{1}}\#^{v_{2}}$ $G_{2}$ is the graph

\strut

\begin{center}
$G$ $=\qquad 
\begin{array}{lllll}
\bullet & \longrightarrow & \overset{v_{12}}{\bullet } & \longrightarrow & 
\bullet \\ 
&  & \,\uparrow &  &  \\ 
&  & \,\,\bullet &  & 
\end{array}
.$
\end{center}

\strut

\begin{definition}
The graph $G$ $=$ $G_{1}$ $^{v_{1}}\#^{v_{2}}$ $G_{2}$ is called the glued
graph of $G_{1}$ and $G_{2}$ by gluing (or identifying) $v_{1}$ and $v_{2}$.
The vertex $v_{12}$ of $v_{1}$ and $v_{2}$ in $V(G)$ is said to be the
identified (or glued) vertex.
\end{definition}

\strut

Inductively, we can construct the iterated glued graph $G$ of the countable
directed graphs $K$ and $T$ by

\strut

\begin{center}
$G$ $=$ $\left( K\text{ }^{v_{1}}\#^{v_{2}}\text{ }T\right) $ $%
^{v_{3}}\#^{v_{4}}$ $T$ $^{v_{5}}\#^{v_{6}}$ $T$ ...
\end{center}

\strut etc.

\strut

As usual, we can create the shadowed graph $\widehat{G}$ and the
corresponding graph groupoid $\Bbb{G}$ of the glued graph $G.$ Notice that

\strut

\begin{center}
$G_{1}$ $^{v_{1}}\#^{v_{2}}$ $G_{2}$ $=$ $G_{2}$ $^{v_{2}}\#^{v_{1}}$ $%
G_{1}. $
\end{center}

\strut

\strut

\subsubsection{Conditional Iterated Gluing on $\mathcal{G}_{G}$}

\strut

\strut

Let $\mathcal{G}_{G}$ be the corresponding graphs induced by the Wold
decomposed family $\mathcal{G}_{W}$ of $\mathcal{G}$ in $PI(H).$ Define now
the subset $\Bbb{G}_{0}$ of $B(H)$ by

\strut

\begin{center}
$\Bbb{G}_{0}$ $\overset{def}{=}$ $\{$all reduced words in $\Bbb{G}_{0}^{0}\}$
$\subset $ $B(H).$
\end{center}

\strut

i.e., if $y$ $\in $ $\Bbb{G}_{0}$ is nonzero, then there exist $n$ $\in $ $%
\Bbb{N}$ and $p_{1},$ ..., $p_{n}$ $\in $ $\Bbb{G}_{0}^{0}$ such that

\strut

\begin{center}
$y$ $=$ the operator formed by the product $p_{1}$ ... $p_{n}.$
\end{center}

\strut \strut

Remark that, even though $p_{1},$ ..., $p_{n}$ are projections on $H,$ we
cannot guarantee that the operator $y$ is a projection on $H.$ Recall that
the operator product $p$ $q$ of projections $p$ and $q$ is a projection if
and only if $p$ $q$ $=$ $q$ $p$ on $H.$

\strut

First, we will determine the connecting relation on $\widehat{\mathcal{G}_{W}%
}$ $=$ $\mathcal{G}_{W}$ $\cup $ $\mathcal{G}_{W}^{*}$ on $H,$ by defining
the map

$\strut $

\begin{center}
$\pi $ $:$ $\widehat{\mathcal{G}_{W}}$ $\times $ $\widehat{\mathcal{G}_{W}}$ 
$\rightarrow $ $\Bbb{G}_{0}$
\end{center}

\strut by

\begin{center}
$\pi (x,$ $y)$ $\overset{def}{=}$ $(x^{*}$ $x)$ $(y$ $y^{*})$,
\end{center}

\strut

for all $x,$ $y$ $\in $ $\widehat{\mathcal{G}_{W}}.$ Clearly, we can get that

\strut

\begin{center}
$\pi (x,$ $y)$ $=$ $\left\{ 
\begin{array}{ll}
x^{*}\text{ }x & \text{if }x^{*}\text{ }x\leq y\text{ }y^{*}\text{ in }\Bbb{G%
}_{0}^{0} \\ 
y\text{ }y^{*} & \text{if }x^{*}\text{ }x\geq y\text{ }y^{*}\text{ in }\Bbb{G%
}_{0}^{0} \\ 
(x^{*}\text{ }x)(y\text{ }y^{*}) & \text{if }(x^{*}\text{ }x)\text{ }(y\text{
}y^{*})\neq 0\text{ on }H \\ 
0 & \text{otherwise,}
\end{array}
\right. $
\end{center}

\strut

for all $x,$ $y$ $\in $ $\widehat{\mathcal{G}_{W}}.$ Consider that if $\pi
(x,$ $y)$ $\neq $ $0$ in $\Bbb{G}_{0},$ then

\strut

\begin{center}
$x$ $y$ $=$ $x$ $(x^{*}$ $x)$ $(y$ $y^{*})$ $y$ $=$ $x$ $\pi (x,$ $y)$ $y$ $%
\neq $ $0$ on $H.$
\end{center}

\strut

\begin{definition}
The map $\pi $ $:$ $\widehat{\mathcal{G}_{W}}$ $\times $ $\widehat{\mathcal{G%
}_{W}}$ $\rightarrow $ $\Bbb{G}_{0}$ is said to be the $\mathcal{G}_{W}$%
-admissibility map.
\end{definition}

\strut

Now, define a new subset $\Bbb{G}$ of $B(H)$ by the collection of all
reduced words in $\widehat{\mathcal{G}_{W}}.$ i.e.,

\strut

\begin{center}
$\Bbb{G}$ $\overset{def}{=}$ $\{$all reduced words in $\widehat{\mathcal{G}%
_{W}}\}$ $\subset $ $B(H).$
\end{center}

\strut

The reduction on $\Bbb{G}$ is determined by the operator multiplication on $%
B(H),$ and this reduction is totally explained by the $\mathcal{G}_{W}$%
-admissibility map $\pi $:

\strut

$\qquad x_{1},$ ..., $x_{n}$ $\in $ $\widehat{\mathcal{G}_{W}}$ and $w$ $=$ $%
x_{1}$ ... $x_{n}$ $\in $ $\Bbb{G}$ $\setminus $ $\{0\}$

\strut

\begin{center}
$\Longleftrightarrow $ $\pi (x_{i},$ $x_{i+1})$ $\neq $ $0$ in $\Bbb{G}_{0},$
for all $i$ $=$ $1,$ ..., $n$ $-$ $1.$
\end{center}

\strut

Notice that every element $y$ $\in $ $\Bbb{G}_{0}$ is contained in $\Bbb{G},$
too. i.e.,

\strut

\begin{center}
$\Bbb{G}_{0}$ $\subseteq $ $\Bbb{G}.$
\end{center}

\strut

Thus we can extend the $\mathcal{G}_{W}$-admissibility map $\pi $ on $%
\widehat{\mathcal{G}_{W}}$ to that on $\Bbb{G}$, also denoted by $\pi $:

\strut

\begin{center}
$\pi $ $:$ $\Bbb{G}$ $\times $ $\Bbb{G}$ $\rightarrow $ $\Bbb{G}_{0}$
\end{center}

defined by

\begin{center}
$\pi (w_{1},$ $w_{2})$ $\overset{def}{=}$ $\left\{ 
\begin{array}{ll}
\begin{array}{l}
\pi (w_{1},\text{ }w_{2})
\end{array}
& 
\begin{array}{l}
\text{if }(w_{1},\text{ }w_{2})\in \widehat{\mathcal{G}_{W}}\text{ }\times 
\text{ }\widehat{\mathcal{G}_{W}}
\end{array}
\\ 
\begin{array}{l}
w_{1}\text{ }w_{2}
\end{array}
& 
\begin{array}{l}
\text{if }(w_{1},\text{ }w_{2})\in \Bbb{G}_{0}\text{ }\times \text{ }\Bbb{G}%
_{0}
\end{array}
\\ 
\ w_{1}\text{ }(x_{1}\text{ }x_{1}^{*}) & 
\begin{array}{l}
\text{if }w_{1}\in \Bbb{G}_{0}\text{ and }w_{2}\notin \Bbb{G}_{0} \\ 
\text{and }w_{2}=x_{1}\text{ ... }x_{n},\text{ }x_{j}\in \widehat{\mathcal{G}%
_{W}}
\end{array}
\\ 
&  \\ 
\ (x_{n}^{*}\text{ }x_{n})\text{ }w_{2} & 
\begin{array}{l}
\text{if }w_{1}\notin \Bbb{G}_{0}\text{ and }w_{2}\in \Bbb{G}_{0} \\ 
\text{and }w_{1}=x_{1}\text{ ... }x_{n},\text{ }x_{j}\in \widehat{\mathcal{G}%
_{W}}
\end{array}
\\ 
&  \\ 
\ \pi (x_{n},\text{ }y_{1}) & 
\begin{array}{l}
\text{if }w_{1},\text{ }w_{2}\notin \Bbb{G}_{0}\text{ in }\Bbb{G},\text{ and}
\\ 
w_{1}=x_{1}\text{ ... }x_{n}\text{ and }w_{2}=y_{1}\text{ ... }y_{m}, \\ 
\text{where }x_{j},\text{ }y_{i}\text{ }\in \text{ }\widehat{\mathcal{G}_{W}}%
,
\end{array}
\end{array}
\right. $
\end{center}

\strut

for all $w_{1},$ $w_{2}$ $\in $ $\Bbb{G}.$

\strut

Denote now the collection of the shadowed graphs $\widehat{G_{x}}$'s of $%
G_{x}$ $\in $ $\mathcal{G}_{G}$ by $\widehat{\mathcal{G}_{G}}.$ i.e.,

\strut

\begin{center}
$\widehat{\mathcal{G}_{G}}$ $\overset{def}{=}$ $\{\widehat{G_{x}}$ $:$ $%
G_{x} $ $\in $ $\mathcal{G}_{G}\}.$
\end{center}

\strut

Take $G_{1}$ and $G_{2}$ in $\mathcal{G}_{G}$ and choose $v_{1}$ $\in $ $%
V(G_{1})$ and $v_{2}$ $\in $ $V(G_{2}).$ For the chosen vertices $v_{1}$ and 
$v_{2},$ we can compute

\strut

\begin{center}
$\pi (v_{1},$ $v_{2})$ $=$ $v_{1}$ $v_{2}$ $\in $ $\Bbb{G}_{0}$ in $B(H),$
\end{center}

\strut

since $V(G_{x})$ $\in $ $\Bbb{G}_{0}^{0}$ $\subset $ $\Bbb{G}_{0},$ for all $%
G_{x}$ $\in $ $\mathcal{G}_{G}.$

\strut

Observe now the partial ordering on $\mathcal{G}_{G}.$ Define

\strut

\begin{center}
$G_{1}$ $\leq $ $G_{2}$ in $\mathcal{G}_{G}$ $\overset{def}{%
\Longleftrightarrow }$ $G_{1}$ is a full-subgraph of $G_{2}.$
\end{center}

\strut

Recall that we say that $K_{1}$ is a full-subgraph of $K_{2},$ where $K_{1}$
and $K_{2}$ are countable directed graphs, if

\strut

\begin{center}
$E(K_{1})$ $\subseteq $ $E(K_{2}),$
\end{center}

and

\begin{center}
$V(K_{1})$ $=$ $\{v,$ $v^{\prime }$ $\in $ $V(K_{2})$ $:$ $e$ $=$ $v$ $e$ $%
v^{\prime },$ $\forall $ $e$ $\in $ $E(K_{1})\}.$
\end{center}

\strut

Remark the difference between subgraphs and full-subgraphs. We say that $%
K_{1}$ is a subgraph of $K_{2}$ if

\strut

\begin{center}
$V(K_{1})$ $\subseteq $ $V(K_{2})$
\end{center}

and

\begin{center}
$E(K_{1})$ $=$ $\{e$ $\in $ $E(K_{2})$ $:$ $e$ $=$ $v$ $e$ $v^{\prime },$ $%
\forall $ $v,$ $v^{\prime }$ $\in $ $V(K_{1})\}.$
\end{center}

\strut

Our partial ordering $\leq $ on $\mathcal{G}_{G}$ is determined by the
concept, ``full''-subgraphs.

\strut

\textbf{Notation} From now, if we denote $\mathcal{G}_{G},$ then it means
the partial ordered set $(\mathcal{G}_{G},$ $\leq ).$ $\square $

\strut

We now can construct the ($\pi $-depending) conditional glued graph $G_{1}$ $%
^{v_{1}}\#_{\pi }^{v_{2}}$ $G_{2},$ by a directed graph,

\strut

\begin{center}
$G_{1}$ $^{v_{1}}\#_{\pi }^{v_{2}}$ $G_{2}$ $\overset{def}{=}$ $\left\{ 
\begin{array}{ll}
G_{1}\text{ }^{v_{1}}\#^{v_{2}}\text{ }G_{2} & \,\,\,\text{if }\pi (v_{1},%
\text{ }v_{2})\neq 0 \\ 
\begin{array}{l}
G_{2} \\ 
G_{1}
\end{array}
& 
\begin{array}{l}
\text{if }G_{1}\leq G_{2} \\ 
\text{if }G_{1}\geq G_{2}
\end{array}
\\ 
G_{1}\text{ }\cup \text{ }G_{2} & \,\,\,\,\text{otherwise,}
\end{array}
\right. $
\end{center}

\strut \strut

\strut where $G_{1}$ $\cup $ $G_{2}$ is the directed graph with

\strut

\begin{center}
$V\left( G_{1}\cup G_{2}\right) $ $=$ $V(G_{1})$ $\sqcup $ $V(G_{2})$
\end{center}

and

\begin{center}
$E\left( G_{1}\cup G_{2}\right) $ $=$ $E(G_{1})$ $\sqcup $ $E(G_{2}),$
\end{center}

\strut

where the symbol ``$\sqcup $'' on the right-hand sides means the disjoint
union.

\strut \strut

Now, consider the conditional gluing $\#_{\pi }$ on $\mathcal{G}_{G}$ as
follows: if $G_{1},$ $G_{2}$ $\in $ $\mathcal{G}_{G}$ and if

\strut

\begin{center}
$V(G_{1})$ $=$ $\{p_{1},$ $p_{2},$ $p_{3},$ ...$\}$ and $V(G_{2})$ $=$ $%
\{q_{1},$ $q_{2},$ $q_{3},$ ...$\},$
\end{center}

\strut \strut

in $\Bbb{G}_{0}^{0},$ then the conditional glued graph $G_{1}$ $\#_{\pi }$ $%
G_{2}$ of $G_{1}$ and $G_{2}$ is defined by

\strut

\strut (3.3)

$\qquad \qquad G_{1}$ $^{p_{1}}\#_{\pi }^{q_{1}}$ $G_{2}$ $^{p_{1}}\#_{\pi
}^{q_{2}}$ $G_{2}$ $^{p_{1}}\#_{\pi }^{q_{3}}$ $G_{2}$ ...

\strut

\begin{center}
$G_{1}$ $^{p_{2}}\#_{\pi }^{q_{1}}$ $G_{2}$ $^{p_{2}}\#_{\pi }^{q_{2}}$ $%
G_{2}$ $^{p_{2}}\#_{\pi }^{q_{3}}$ $G_{2}$ ...
\end{center}

\qquad ...

\begin{center}
$G_{1}$ $^{p_{n}}\#_{\pi }^{q_{1}}$ $G_{2}$ $^{p_{n}}\#_{\pi }^{q_{2}}$ $%
G_{2}$ $^{p_{n}}\#_{\pi }^{q_{3}}$ $G_{2}$ ...
\end{center}

\qquad ....

\strut

\begin{definition}
The graph (3.3) of $G_{1}$ and $G_{2}$ in $\mathcal{G}_{G}$ is denoted by $%
G_{1}$ $\#_{\pi }$ $G_{2}.$ And we call it the conditional (or the $\pi $%
-dependent) glued graph of $G_{1}$ and $G_{2}.$ And the symbol ``$\#_{\pi }$%
'' is called the conditional (or the $\pi $-dependent) gluing.
\end{definition}

\strut

By the conditional \emph{iterated} gluing on $\mathcal{G}_{G},$ we can
decide the countable directed graph $G_{\mathcal{G}_{W}}$ as a directed
graph with

\strut

\strut (3.4)

\begin{center}
$G_{\mathcal{G}_{W}}$ $\overset{def}{=}$ $\underset{x\in \mathcal{G}_{W}}{%
\#_{\pi }}$ $G_{x}.$
\end{center}

\strut

\begin{definition}
The conditional iterated glued graph $G_{\mathcal{G}_{W}}$ defined in (3.4)
is called the $\mathcal{G}$-graph. And the corresponding graph groupoid $%
\Bbb{G}_{\mathcal{G}_{W}}$ of $G_{\mathcal{G}_{W}}$ is called the $\mathcal{G%
}$-groupoid.
\end{definition}

\strut

Notice that every element in $\Bbb{G}_{\mathcal{G}_{W}}$ is the reduced
words in $\widehat{\mathcal{G}_{W}},$ under the operator multiplication on $%
B(H).$ Thus we can have the following lemma.

\strut

\begin{lemma}
$\Bbb{G}_{\mathcal{G}_{W}}$ $=$ $\Bbb{G},$ where $\Bbb{G}$ $\overset{def}{=}$
$\{$all reduced words in $\widehat{\mathcal{G}_{W}}\}$ defined at the
beginning of this subsection. $\square $
\end{lemma}

\strut

\begin{example}
Suppose that $\mathcal{G}_{W}$ $=$ $\mathcal{G}^{(u)}$ $\cup $ $\mathcal{G}%
_{\infty }^{(s)}$ and assume that

\strut 

\begin{center}
$\pi (x,$ $y)$ $\neq $ $0,$ for $(x,$ $y)$ $\in $ $\mathcal{G}_{W}$ $\times $
$\mathcal{G}_{W}$ $\Longleftrightarrow $ $x^{*}$ $x$ $=$ $y$ $y^{*}$ in $%
\Bbb{G}_{0}^{0}.$
\end{center}

\strut 

Then this finite family $\mathcal{G}_{W}$ of partial isometries in $B(H)$ is
a $G_{\mathcal{G}_{W}}$-family in the sense of [13] (Also See Section 2.3).
\end{example}

\strut

\begin{example}
Let $\mathcal{G}$ $=$ $\{a\}$ and assume that $\mathcal{G}_{W}$ $=$ $\{u,$ $%
s\},$ where $a$ has its Wold decomposition $a$ $=$ $u$ $+$ $s,$ with $u$ $%
\neq $ $0$ and $s$ $\neq $ $0$ on $H$, equivalently, the $*$-isomorphic
index $i_{*}(a)$ of $a$ is $(k_{1},$ $k_{2},$ $k_{3},$ $k_{4})$ in $\Bbb{N}%
_{0}^{\infty },$ satisfying that $k_{1}$ $\neq $ $0,$ $k_{3}$ $=$ $\infty $
in $\Bbb{N}_{0}^{\infty }.$ Then we can create a family $\mathcal{G}_{G}$ $=$
$\{G_{u},$ $G_{s}\}$ of directed graphs. In particular,

\strut 

\begin{center}
$V(G_{u})$ $=$ $\{u^{*}$ $u$ $=$ $u$ $u^{*}\}$ and $E(G_{u})$ $=$ $\{u\}$
\end{center}

and

\begin{center}
$V(G_{s})$ $=$ $\{s^{*}$ $s,$ $s$ $s^{*}\}$ and $E(G_{s})$ $=$ $\{s\}.$
\end{center}

\strut \strut 

Assume that

\strut \strut 

\begin{center}
$\pi (u,$ $s)$ $=$ $(u^{*}$ $u)$ $(s$ $s^{*})$ $=$ $(u$ $u^{*})$ $(s$ $s^{*})
$ $=$ $\pi (u^{*},$ $s)$ $\neq $ $0,$
\end{center}

\strut and

\begin{center}
$\pi (u,$ $s^{*})$ $=$ $(u^{*}$ $u)$ $(s^{*}$ $s)$ $=$ $(u$ $u^{*})$ $(s^{*}$
$s)$ $=$ $\pi (u^{*},$ $s^{*})$ $=$ $0.$
\end{center}

\strut 

Then we can construct the conditional iterated glued graph $G$ $=$ $G_{u}$ $%
\#_{\pi }$ $G_{s},$

\strut 

\begin{center}
$G$ $=$ $G_{u}$ $\#_{\pi }$ $G_{s}$ $=$ $\left( G_{u}\text{ }%
^{uu^{*}}\#^{ss^{*}}G_{s}\right) ,$
\end{center}

\strut \strut 

Then this graph $G$ is graph-isomorphic to the graph $K$,

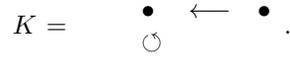
\strut 

\begin{center}
$K$ $=\qquad 
\begin{array}{lll}
\bullet  & \longleftarrow  & \bullet  \\ 
\circlearrowleft  &  & 
\end{array}
.$\strut 
\end{center}
\end{example}

\strut

\begin{example}
Let $\mathcal{G}$ $=$ $\{u_{1},$ $u_{2}\}$ $=$ $\mathcal{G}^{(u)}$ be a
finite family of partial isometries in $B(H).$ Then we can have $\mathcal{G}%
_{G}$ $=$ $\{G_{u_{1}},$ $G_{u_{2}}\},$ where $G_{k}$ are the graph with

\strut 

\begin{center}
$V(G_{k})$ $=$ $\{u_{k}^{*}$ $u_{k}$ $=$ $u_{k}$ $u_{k}^{*}\}$ and $E(G_{k})$
$=$ $\{u_{k}\},$
\end{center}

\strut 

for $k$ $=$ $1,$ $2.$ Suppose $\pi (u_{1},$ $u_{2})$ $=$ $\pi (u_{1}^{*},$ $%
u_{2})$ $=$ $\pi (u_{1},$ $u_{2}^{*})$ $=$ $\pi (u_{1}^{*},$ $u_{2}^{*})$ $%
\neq $ $0$ in $\Bbb{G}_{0}^{0}.$ Then we can construct the iterated glued
graph $G$ $=$ $G_{u_{1}}$ $\#_{\pi }$ $G_{u_{2}}$ which is identified with $%
G_{u_{1}}$ $^{u_{1}^{*}u_{1}}\#_{\pi }^{u_{2}^{*}u_{2}}$ $G_{u_{2}}.$ Then
this graph $G$ is graph-isomorphic to the graph $K$ with

\strut 

\begin{center}
$V(K)$ $=$ $\{v\}$ and $E(K)$ $=$ $\{e_{1}$ $=$ $v$ $e_{1}$ $v,$ $e_{2}$ $=$ 
$v$ $e_{2}$ $v\}.$
\end{center}

\strut 

Assume now that $\pi (u_{1},$ $u_{2})$ $=$ $0$ in $\Bbb{G}_{0}^{0}.$ Then
the iterated glued graph $G$ $=$ $G_{u_{1}}$ $\#_{\pi }$ $G_{u_{2}}$ is
identified with the graph $G_{u_{1}}$ $\sqcup $ $G_{u_{2}}.$ This graph is
graph-isomorphic to the graph $H$ with

\strut 

\begin{center}
$V(H)$ $=$ $\{v_{1},$ $v_{2}\}$ and $E(H)$ $=$ $\{e_{1}$ $=$ $v_{1}$ $e_{1}$ 
$v_{1},$ $e_{2}$ $=$ $v_{2}$ $e_{2}$ $v_{2}\}.$
\end{center}
\end{example}

\strut

\begin{example}
Let $\mathcal{G}$ $=$ $\{s_{1},$ $s_{2}\}$ $=$ $\mathcal{G}_{\infty }^{(s)}.$
Then we have the corresponding family $\mathcal{G}_{G}$ $=$ $\{G_{s_{1}},$ $%
G_{s_{2}}\}$ of directed graphs. Suppose $\pi (s_{1},$ $s_{2})$ $\neq $ $0,$
in $\Bbb{G}.$ Then the conditional iterated glued graph $\Delta $ $=$ $%
G_{s_{1}}$ $\#_{\pi }$ $G_{s_{2}}$ is graph-isomorphic to the graph $\Delta $%
:

\strut 

\begin{center}
$\Delta $ $=$ $\quad 
\begin{array}{lll}
\bullet  & \longrightarrow  & \bullet  \\ 
&  & \downarrow  \\ 
&  & \bullet 
\end{array}
.$
\end{center}

\strut 

Here, we can understand the horizontal part

\strut 

\begin{center}
$\bullet \longrightarrow \bullet $
\end{center}

\strut 

of $\Delta $ are graph-isomorphic to $G_{s_{2}},$ and the vertical part

\strut 

\begin{center}
$
\begin{array}{l}
\bullet  \\ 
\downarrow  \\ 
\bullet 
\end{array}
$
\end{center}

\strut 

of $\Delta $ are graph-isomorphic to $G_{s_{1}}.$
\end{example}

\strut

The more examples would be provided in Section 5. In the following section,
we will observe how the $\mathcal{G}$-groupoid $\Bbb{G}_{\mathcal{G}_{W}}$
induced by the $\mathcal{G}$-graph $G_{\mathcal{G}_{W}}$ works on the fixed
Hilbert space $H.$

\strut \strut

\strut

\strut

\subsubsection{The Vertex Set of the $\mathcal{G}$-Graph}

\strut

\strut

In this subsection, we observe the vertex set $V(G)$ of the $\mathcal{G}$%
-graph $G,$ which is the conditional iterated glued graph of $\mathcal{G}%
_{G} $ $=$ $\{G_{x}$ $:$ $x$ $\in $ $\mathcal{G}_{W}\}.$ Suppose $x,$ $y$ $%
\in $ $\mathcal{G}_{W},$ and assume that $\pi (x,$ $y)$ $\neq $ $0$ in $\Bbb{%
G}_{0}. $ Then we can identify the vertices $x^{*}$ $x$ $\in $ $V(G_{x})$
and $y$ $y^{*}$ $\in $ $G_{y},$ and we can create the glued graph,

$\strut $

\begin{center}
$G_{x,y}$ $=$ $G_{x}$ $^{x^{*}x}\#_{\pi }^{yy^{*}}$ $G_{y},$
\end{center}

\strut

which is a full-subgraph of the $\mathcal{G}$-graph $G.$ Then, inside $%
G_{x,y},$ the vertices $x^{*}x$ $\in $ $V(G_{x})$ and $y$ $y^{*}$ $\in $ $%
V(G_{y})$ are identified. Denote the identified vertex by $p_{0}.$ Then, how
we can understand this vertex $p_{0}$? Combinatorially (in the sense of
Subsection 3.2.1), we can simply understand $p_{0}$ is the identified vertex
which is a pure combinatorial object. However, operator-theoretically, it is
not clear to see how this vertex $p_{0}$ works as an operator on $H$. In
fact, this vertex $p_{0}$ is not determined uniquely as an operator.
Operator-theoretically, we can understand $p_{0}$ as:

$\strut $

\begin{center}
$p_{0}$ $=$ $x^{*}$ $x$ or $p_{0}$ $=$ $y$ $y^{*}$ or $p_{0}$ $=$ $\pi (x,$ $%
y)$
\end{center}

or

\begin{center}
$p_{0}$ $=$ $\pi (x,$ $y)^{*}$ $=$ $\pi (y^{*},$ $x^{*}),$
\end{center}

\strut \strut

case by case. Notice that

$\strut $

\begin{center}
$\pi (x,$ $y)$ $\neq $ $0$ in $\Bbb{G}_{0}$ if and only if $\pi (y^{*},$ $%
x^{*})$ $\neq $ $0$ in $\Bbb{G}_{0}$,
\end{center}

\strut

since $\pi (x,$ $y)^{*}$ $=$ $\pi (y^{*},$ $x^{*}).$ Anyway, the vertex $%
p_{0}$ just represents the connection on $G_{x,\,y}$ (or on $G$), in terms
of the $\mathcal{G}_{W}$-admissibility (map $\pi $).

\strut

The best way to choose $p_{0}$ $\in $ $V(G_{x,y})$ $\subset $ $V(G),$ as an
operator (\emph{if we have to: but we have not to}), is that: $p_{0}$ may be
the projection from $H$ onto the subspace $H_{init}^{x}$ $\cap $ $%
H_{fin}^{y}.$ So, by gluing $x^{*}$ $x$ and $y$ $y^{*}$ to $p_{0}$ as the
projection onto $H_{init}^{x}$ $\cap $ $H_{fin}^{y}$
(operator-theoretically), we may construct the conditional glued graph $%
G_{x,y}.$ However, we will not fix $p_{0},$ as an operator. We regard $p_{0}$
as a pure combinatorial object representing the connection of the operators $%
x$ and $y$. In other words, the identified vertex $p_{0}$ can be explained
case by case, differently in Operator Theory point of view.

\strut

\textbf{Observation} Combinatorially, the identified vertices of the $%
\mathcal{G}$-graph represent the $\mathcal{G}_{W}$-admissibility.
Operator-theoretically, these vertices are understood differently case by
case. For instance, the identified vertex $p_{0}$ in the text represents $%
x^{*}$ $x$ or $y$ $y^{*}$ or $\pi (x,$ $y)$ or $\pi (y^{*},$ $x^{*}),$ or
more depending on the $\mathcal{G}_{W}$-admissibility,
operator-theoretically. $\square $

\strut

\strut \strut \strut

\strut \strut

\subsection{A Representation of the $\mathcal{G}$-Groupoid}

\strut

\strut

\strut

Throughout this section, we will use the same notations we used in the
previous sections. We will consider a certain representation of the $%
\mathcal{G}$-groupoid $\Bbb{G},$ which is the graph groupoid induced by the
conditional iterated glued graph, the $\mathcal{G}$-graph $G$. Recall that,
the operator $w$ $\in $ $B(H)$ is contained in $\Bbb{G}$ if and only if $w$ $%
=$ $x_{1}$ ... $x_{n}$, for $n$ $\in $\thinspace $\Bbb{N}$ and $x_{1},$ ..., 
$x_{n}$ $\in $ $E(\widehat{G})$ $=$ $\widehat{\mathcal{G}_{W}}.$

\strut \strut \strut \strut

\strut Let $x$ $\in $ $E(\widehat{G}).$ Then we can construct the subspace $%
H_{x}$ $=$ $(x^{*}$ $x)$ $H$ of $H$, where $H_{x}$ is a Hilbert space where $%
x$ $\in $ $\widehat{\mathcal{G}_{W}}$ is acting on, as a unitary or a shift.
Then we can have the class $\mathcal{H}_{\widehat{\mathcal{G}_{W}}}$ of
subspaces of $H$,

\strut

\begin{center}
$\mathcal{H}_{\widehat{\mathcal{G}_{W}}}$ $\overset{def}{=}$ $\{H_{x}$ $=$ $%
(x^{*}$ $x)$ $H$ $:$ $x$ $\in E(\widehat{G})\}.$
\end{center}

\strut \strut

Let $H_{x_{1}},$ $H_{x_{2}}$ $\in $ $\mathcal{H}_{\widehat{\mathcal{G}_{W}}%
}. $ Then we can define subspaces of $H,$

\strut

\begin{center}
$H_{x_{1}}$ $\wedge _{\pi }$ $H_{x_{2}}$ and $H_{x_{1}}$ $\vee _{\pi }$ $%
H_{x_{2}},$
\end{center}

\strut by

\begin{center}
$H_{x_{1}}$ $\wedge _{\pi }$ $H_{x_{2}}$ $\overset{def}{=}$ $\left\{ 
\begin{array}{ll}
H_{x_{1}}\cap H_{x_{2}} & \text{if }\pi (x_{1},\text{ }x_{2})\neq 0 \\ 
&  \\ 
\{0_{H}\} & \text{otherwise.}
\end{array}
\right. $
\end{center}

and

\begin{center}
$H_{x_{1}}$ $\vee _{\pi }$ $H_{x_{2}}$ $\overset{def}{=}$ $\left\{ 
\begin{array}{ll}
\overline{span(H_{x_{1}}\cup H_{x_{2}})}^{H} & \text{if }\pi (x_{1},\text{ }%
x_{2})\neq 0 \\ 
&  \\ 
H_{x_{1}}\text{ }\oplus \text{ }H_{x_{2}} & \text{otherwise.}
\end{array}
\right. $
\end{center}

\strut

Then we can define the Hilbert space $H_{\mathcal{G}}$ by

\strut

\strut (3.5)

\begin{center}
$H_{\mathcal{G}}$ $\overset{def}{=}$ $\underset{w\in FP_{r}(\widehat{G})}{%
\vee _{\pi }}$ $H_{w}$
\end{center}

with

\begin{center}
$H_{w}$ $\overset{def}{=}$ $\overset{n}{\underset{j=1}{\,\wedge _{\pi }}}$ $%
H_{x_{j}},$
\end{center}

whenever

\begin{center}
$w$ $=$ $x_{1}$ ... $x_{n}$ $\in $ $FP_{r}(\widehat{G}),$ with $x_{1},$ ..., 
$x_{n}$ $\in $ $E(\widehat{G}),$
\end{center}

\strut where

\begin{center}
$FP_{r}(\widehat{G})$ $\overset{def}{=}$ $\Bbb{G}$ $\setminus $ $\left( V(%
\widehat{G})\text{ }\cup \text{ }\{0\}\right) .$
\end{center}

\strut

Clearly, the Hilbert space $H_{\mathcal{G}}$ is a well-determined subspace
of $H,$ and every operator $w$ $\in $ $\Bbb{G}$ $\subset $ $B(H)$ acts on $%
H_{\mathcal{G}}.$ Suppose $v$ $\in $ $V(\widehat{G}).$ Then there always
exists $e$ $\in $ $E(\widehat{G})$ such that $(e^{*}$ $e)$ $H$ $\overset{%
\text{Subspace}}{\supseteq }$ $v$ $H$ or $(e$ $e^{*})$ $H$ $\overset{\text{%
Subspace}}{\supseteq }$ $v$ $H$. So, all operators in $\Bbb{G}$ acts on the
subspace $H_{\mathcal{G}}$ of $H.$

\strut \strut

We now define the groupoid action $\alpha $ of $\Bbb{G}$ acting on $H_{%
\mathcal{G}},$

$\strut $

\begin{center}
$\alpha $ $:$ $\Bbb{G}$ $\rightarrow $ $B(H_{\mathcal{G}})$
\end{center}

\strut \strut

sending $w$ to $\alpha _{w}$ $\in $ $B(H_{\mathcal{G}}),$ by

\strut

\strut (3.6)$\qquad \qquad \qquad \qquad \qquad \alpha _{w}$ $=$ $w,$ for
all $w$ $\in $ $\Bbb{G}.$

\strut

Then, clearly, the action $\alpha $ is a groupoid action of $\Bbb{G}$ on $H_{%
\mathcal{G}}.$ Indeed,

\strut

\begin{center}
$\alpha _{w_{1}w_{2}}$ $=$ $w_{1}$ $w_{2}$ $=$ $\alpha _{w_{1}}$ $\alpha
_{w_{2}},$ for all $w_{1},$ $w_{2}$ $\in $ $\Bbb{G}.$
\end{center}

\strut \strut

We can easily check that $\alpha _{w_{1}w_{2}}$ $\neq $ $0$ if and only if $%
\pi (w_{1},$ $w_{2})$ $\neq $ $0,$ for all $w_{1},$ $w_{2}$ $\in $ $\Bbb{G}.$
Now, we can determine a representation of the $\mathcal{G}$-groupoid $\Bbb{G}
$ by the pair $(H_{\mathcal{G}},$ $\alpha ).$

\strut

\begin{definition}
The representation $(H_{\mathcal{G}},$ $\alpha )$ of the $\mathcal{G}$%
-groupoid $\Bbb{G}$ is called the $\mathcal{G}$-representation in $B(H).$
\end{definition}

\strut

In the following theorem, we observe the equivalence of $\mathcal{G}$%
-representations.

\strut

\begin{theorem}
Let $\mathcal{G}_{k}$ $=$ $\{a_{1}^{(k)},$ ..., $a_{N_{k}}^{(k)}\}$ $\subset 
$ $PI(H)$ be the finite families of partial isometries in $B(H),$ and let $%
G_{k}$ be the $\mathcal{G}_{k}$-graphs having the $\mathcal{G}_{k}$%
-groupoids $\Bbb{G}_{k},$ respectively, for $k$ $=$ $1,$ $2.$ Assume that
(i) the shadowed graphs $\widehat{G_{k}}$ of $G_{k}$ are graph-isomorphic,
via the graph-isomorphism $g$ $:$ $\widehat{G_{1}}$ $\rightarrow $ $\widehat{%
G_{2}},$ and (ii) $i_{*}(x)$ $-$ $i_{*}(g(x))$ $=$ $(0,$ $0,$ $0,$ $0)$ in $%
\Bbb{N}_{0}^{\infty },$ for all $x$ $\in $ $E(\widehat{G_{1}}).$ Then the $%
\mathcal{G}_{k}$-representations $(H_{\mathcal{G}_{k}},$ $\alpha ^{(k)})$
are equivalent, for $k$ $=$ $1,$ $2.$
\end{theorem}

\strut

\begin{proof}
Suppose $\mathcal{G}_{1}$ and $\mathcal{G}_{2}$ are the finite families of
partial isometries in $B(H),$ and let $\mathcal{G}_{k:W}$ $=$ $\mathcal{G}%
_{k}^{(u)}$ $\cup $ $\mathcal{G}_{k}^{(s)},$ where $\mathcal{G}_{k}^{(u)}$
and $\mathcal{G}_{k}^{(s)}$ are the collection of all unitary parts and all
shift parts of $\mathcal{G}_{k},$ respectively, for $k$ $=$ $1,$ $2.$

\strut

Also, by the condition (i), we can conclude that the $\mathcal{G}_{k}$%
-groupoids $\Bbb{G}_{k}$ are groupoid-isomorphic, for $k$ $=$ $1,$ $2.$
Indeed, we can have the groupoid-isomorphism $\varphi $ $:$ $\Bbb{G}_{1}$ $%
\rightarrow $ $\Bbb{G}_{2}$:

\strut

\begin{center}
$\varphi (w)$ $\overset{def}{=}$ $\left\{ 
\begin{array}{ll}
0 & \text{if }w=0 \\ 
g(w) & \text{if }w\in V(\widehat{G})\cup E(\widehat{G}) \\ 
g(x_{1})\text{ ... }g(x_{n}) & \text{if }w=e_{1}\text{ ... }e_{n}\in FP_{r}(%
\widehat{G}),\text{ }n>1,
\end{array}
\right. $
\end{center}

\strut

in $\Bbb{G}_{2},$ for all $w$ $\in $ $\Bbb{G}_{1}$ (Also, see [9] and [10]).
This means that algebraically the groupoidal structures of $\Bbb{G}_{1}$ and 
$\Bbb{G}_{2}$ are same.

\strut

By the condition (ii), the subspaces $H_{x}$ $=$ $(x^{*}$ $x)$ $H$ and $%
H_{g(x)}$ $=$ $(g(x)^{*}$ $g(x))$ $H$ are Hilbert-space isomorphic in $H,$
for all $x$ $\in $ $E(G_{1}).$ Therefore, the Hilbert spaces $H_{\mathcal{G}%
_{1}}$ and $H_{\mathcal{G}_{2}}$ are Hilbert-space isomorphic, as embedded
subspaces in $H.$ i.e., we can have the Hilbert-space isomorphism $\Phi $ $:$
$H_{\mathcal{G}_{1}}$ $\rightarrow $ $H_{\mathcal{G}_{2}}$ induced by the
groupoid-isomorphism $\varphi $ $:$ $\Bbb{G}_{1}$ $\rightarrow $ $\Bbb{G}%
_{2} $:

\strut

\begin{center}
$\Phi \left( H_{w}\right) $ $=$ $H_{\varphi (w)},$ in $H_{\mathcal{G}_{2}},$
for all $H_{w}$ $\subset $ $H_{\mathcal{G}_{1}},$ for all $w$ $\in $ $FP_{r}(%
\widehat{G}).$
\end{center}

\strut

Therefore, the Hilbert spaces $H_{\mathcal{G}_{1}}$ and $H_{\mathcal{G}_{2}}$
are Hilbert-space isomorphic, via $\Phi $.

\strut

Since we have the following commuting diagram,

\strut

\begin{center}
$
\begin{array}{lll}
H_{\mathcal{G}_{1}} & \overset{\Phi }{\longrightarrow } & H_{\mathcal{G}_{2}}
\\ 
\downarrow _{\alpha _{w}^{(1)}} &  & \downarrow _{\alpha _{\varphi
(w)}^{(2)}} \\ 
H_{\mathcal{G}_{1}} & \underset{\Phi }{\longrightarrow } & H_{\mathcal{G}%
_{2}}
\end{array}
$
\end{center}

\strut

for all $w$ $\in $ $\Bbb{G}_{1},$ the groupoid actions $\alpha ^{(k)}$
satisfy that

\strut

\begin{center}
$\alpha _{\varphi (w)}^{(2)}$ $=$ $\varphi (w)=$ $\Phi $ $w$ $\Phi ^{-1}$ $=$
$\Phi $ $\alpha _{w}$ $\Phi ^{-1}$ $=$ $\alpha _{\Phi \,w\,\Phi
^{-1}}^{(1)}, $
\end{center}

\strut

for all $w$ $\in $ $\Bbb{G}_{1}.$ Therefore, the actions $\alpha ^{(1)}$ and 
$\alpha ^{(2)}$ are equivalent. Thus the $\mathcal{G}_{k}$-representations $%
(H_{\mathcal{G}_{k}},$ $\alpha ^{(k)})$ are equivalent, for $k$ $=$ $1,$ $2.$
\end{proof}

\strut \strut

The above theorem shows that the equivalence of the $\mathcal{G}_{k}$%
-representations $(H_{\mathcal{G}_{k}},$ $\alpha ^{(k)}),$ for $k$ $=$ $1,$ $%
2,$ is determined by the combinatorial conditions of $\mathcal{G}_{k:W},$
for $k$ $=$ $1,$ $2,$ where $\mathcal{G}_{k:W}$ $=$ $\mathcal{G}_{k}^{(u)}$ $%
\cup $ $\mathcal{G}_{k}^{(s)}.$

\strut \strut \strut \strut \strut

\strut \strut \strut

\strut \strut \strut \strut

\subsection{$C^{*}$-Subalgebras Generated by Partial Isometries}

\strut

\strut

As usual, throughout this section, we will use the same notations we used in
the previous sections. Let $\mathcal{G}$ be a finite family of partial
isometries in $B(H)$, with its Wold decomposed family $\mathcal{G}_{W}.$ Let 
$x$ $\in $ $\mathcal{G}_{W}.$ We already observed the $C^{*}$-subalgebra $%
\mathcal{A}_{x}$ $\overset{denote}{=}$ $C^{*}(\{x\})$ of $B(H),$ in Section
3.1. Notice that, by the very definition, $\mathcal{A}_{x}$ $=$ $\mathcal{A}%
_{x^{*}},$ for all $x$ $\in $ $\mathcal{G}_{W}.$ Here, the symbol ``$=$''
means ``being identically same in $B(H)$''.

\strut

\begin{proposition}
(See Section 3.1) Let $x$ $\in $ $\widehat{\mathcal{G}_{W}}$ and $\mathcal{A}%
_{x}$ $=$ $C^{*}(\{x\}).$ Then

\strut 

\strut (3.7)

\begin{center}
$\mathcal{A}_{x}$ $\overset{*\text{-isomorphic}}{=}$ $\left\{ 
\begin{array}{ll}
(\Bbb{C}\cdot 1_{H_{x}})\otimes _{\Bbb{C}}\text{ }C\left( spec(x)\right)  & 
\begin{array}{l}
\text{if }i^{*}(x)=(k_{1},\text{ }k_{2},\text{ }0,\text{ }k_{4}), \\ 
\quad \text{for }k_{1},\text{ }k_{2},\text{ }k_{4}\in \Bbb{N}_{0}^{\infty }
\end{array}
\\ 
&  \\ 
\mathcal{T}(H_{x}) & 
\begin{array}{l}
\text{if }i^{*}(x)=(0,\text{ }k_{2},\text{ }k_{3},\text{ }k_{4}) \\ 
\text{for }k_{2},\text{ }k_{4}\in \Bbb{N}_{0}^{\infty }\text{ and }k_{3}\in 
\Bbb{N}
\end{array}
\\ 
&  \\ 
(\Bbb{C}\cdot 1_{H_{x}})\otimes _{\Bbb{C}}\text{ }M_{2}(\Bbb{C}) & 
\begin{array}{l}
\text{if }i_{*}(x)=(0,\text{ }k_{2},\text{ }\infty ,\text{ }k_{4}) \\ 
\qquad \text{for }k_{2},\text{ }k_{4}\in \Bbb{N}_{0}^{\infty },
\end{array}
\end{array}
\right. $
\end{center}

\strut 

where $i_{*}(\cdot )$ means the $*$-isomorphic index. $\square $
\end{proposition}

\strut \strut

The following theorem is the main result of this paper. This shows that a $%
C^{*}$-subalgebra of $B(H)$ generated by finitely many partial isometries is 
$*$-isomorphic to a groupoid $C^{*}$-algebra.

\strut

\begin{theorem}
Let $\mathcal{G}$ be a finite family of partial isometries in $B(H)$ and
assume that it has its $\mathcal{G}$-groupoid $\Bbb{G}.$ Then the $C^{*}$%
-algebra $C^{*}(\mathcal{G})$ generated by $\mathcal{G}$ is $*$-isomorphic
to the groupoid $C^{*}$-algebra $C_{\alpha }^{*}(\Bbb{G})$, as embedded $%
C^{*}$-subalgebras of $B(H_{\mathcal{G}})$ $\subseteq $ $B(H),$ where $(H_{%
\mathcal{G}},$ $\alpha )$ is the $\mathcal{G}$-representation of $\Bbb{G}.$
\end{theorem}

\strut

\begin{proof}
Observe that

\strut

\begin{center}
$
\begin{array}{ll}
C^{*}(\mathcal{G}) & =C^{*}(\mathcal{G}_{W})=C^{*}\left( \widehat{\mathcal{G}%
_{W}}\right) \overset{*\text{-isomorphic}}{=}\overline{\Bbb{C}[\alpha (\Bbb{G%
})]} \\ 
&  \\ 
& =C_{\alpha }^{*}(\Bbb{G})\text{ }\overset{C^{*}\text{-subalgebra}}{%
\subseteq }B(H_{\mathcal{G}})\overset{C^{*}\text{-subalgebra}}{\subseteq }%
B(H).
\end{array}
$
\end{center}

\strut

The first identity of the above formulae holds true because of the Wold
decomposition. And the $*$-isomorphic relation and the first $*$-subalgebra
inclusion hold true, by the definition of the $\mathcal{G}$-representation $%
(H_{\mathcal{G}},$ $\alpha )$ of $\Bbb{G}.$
\end{proof}

\strut

The above theorem shows that $C^{*}$-subalgebras of $B(H)$ generated by
finitely many partial isometries on $H$ are $*$-isomorphic to certain
groupoid $C^{*}$-algebras.

\strut

\begin{corollary}
Let $\mathcal{G}_{1}$ and $\mathcal{G}_{2}$ be the collection of finitely
many partial isometries in $B(H)$ having their $\mathcal{G}$-graphs $G_{1}$
and $G_{2}.$ Assume that (i) the shadowed graphs $\widehat{G_{1}}$ and $%
\widehat{G_{2}}$ are graph-isomorphic, via a graph-isomorphism $g$ $:$ $%
\widehat{G_{1}}$ $\rightarrow $ $\widehat{G_{2}},$ (ii) $i_{*}(x)$ $-$ $%
i_{*}(g(x))$ $=$ $(0,$ $0,$ $0,$ $0)$ in $\Bbb{N}_{0}^{\infty },$ for all $x$
$\in $ $E(\widehat{G_{1}})$, and (iii) For any unitary parts $u$ $\in $ $%
E(G_{1})$, $spec(u)$ $=$ $spec\left( g(u)\right) $ in $\Bbb{T}.$ Then the $%
C^{*}$-algebras $C^{*}(\mathcal{G}_{1})$ and $C^{*}(\mathcal{G}_{2})$ are $*$%
-isomorphic, as embedded $C^{*}$-subalgebras of $B(H).$
\end{corollary}

\strut \strut

\begin{proof}
Observe that

\strut

$\qquad C^{*}(\mathcal{G}_{1})$ $\overset{*\text{-isomorphic}}{=}$ $%
C_{\alpha ^{(1)}}^{*}(\Bbb{G}_{1})$

\strut

in $B(H_{\mathcal{G}_{1}}),$ by the previous theorem, where $\Bbb{G}_{k}$
are the $\mathcal{G}_{k}$-groupoids and where $(H_{\mathcal{G}_{k}},$ $%
\alpha ^{(k)})$ are the $\mathcal{G}_{k}$-representations of $\Bbb{G}_{k},$
for $k$ $=$ $1,$ $2$

\strut

$\qquad \qquad \quad \;\overset{*\text{-isomorphic}}{=}$ $C_{\alpha
^{(2)}}^{*}(\Bbb{G}_{2})$ $\overset{*\text{-isomorphic}}{=}$ $C^{*}(\mathcal{%
G}_{2}),$

\strut

of $B(H_{\mathcal{G}_{2}}),$ by (ii) and (iii), since $(H_{\mathcal{G}_{k}},$
$\alpha ^{(k)})$'s are equivalent.
\end{proof}

\strut \strut \strut

\strut

\strut

\section{Block Structures of $C^{*}(\mathcal{G})$}

\strut

\strut

\strut

The main result in this section is a two-part structure theorem (Theorem 4.3
and 4.13) for the $C^{*}$-algebras introduced above. In proving them, we
will need some facts and lemmas from topological (reduced) free product
algebras, introduced in Section 4.1.

\strut 

Throughout this Section, we will keep using the same notations we used
before. Let $\mathcal{G}$ be the given family of finitely many partial
isometries in $B(H)$ and let $\mathcal{G}_{W}$ be the corresponding Wold
decomposed family of $\mathcal{G}.$ We showed that the $C^{*}$-algebra $%
C^{*}(\mathcal{G})$ is $*$-isomorphic to the groupoid $C^{*}$-algebra $%
C_{\alpha }^{*}(\Bbb{G}),$ where $\Bbb{G}$ is the $\mathcal{G}$-groupoid
induced by the $\mathcal{G}$-graph $G$ in the operator algebra $B(H_{%
\mathcal{G}})$ $\overset{C^{*}\text{-subalgebra}}{\subseteq }$ $B(H),$ where 
$(H_{\mathcal{G}},$ $\alpha )$ is the $\mathcal{G}$-representation.

\strut

\strut

\strut

\subsection{Topological Reduced Free Product Algebras}

\strut

\strut

Let $X_{1}$ and $X_{2}$ be arbitrary sets. Then we can define the \emph{%
algebraic free product} $X_{1}$ $*$ $X_{2}$ of $X_{1}$ and $X_{2}$ by the
set of all words in $X_{1}$ $\cup $ $X_{2}.$ i.e.,

\strut

(4.1)$\qquad \qquad \qquad \qquad X_{1}$ $*$ $X_{2}$ $\overset{def}{=}$ all
words in $X_{1}$ $\cup $ $X_{2}.$

\strut

i.e., $x$ $\in $ $X_{1}$ $*$ $X_{2}$ if and only if there exist $n$ $\in $ $%
\Bbb{N}$ and $(i_{1},$ ..., $i_{n})$ $\in $ $\{1,$ $2\}^{n}$ and $x_{i_{j}}$ 
$\in $ $X_{i_{j}},$ for $j$ $=$ $1,$ ..., $n,$ such that

\strut

\begin{center}
$i_{1}$ $\neq $ $i_{2},$ $i_{2}$ $\neq $ $i_{3},$ ..., $i_{n-1}$ $\neq $ $%
i_{n}$
\end{center}

\strut and

\begin{center}
$x$ $=$ $x_{i_{1}}$ $x_{i_{2}}$ ... $x_{i_{n}}.$
\end{center}

\strut \strut

Recall now that we say a topological space $A$ is a topological algebra over
a scalar field $\Bbb{F}$, if it is an algebra over $\Bbb{F}$ and all
operations on $A$ are continuous under the topology for $A$. i.e., $A$ is a
topological algebra over a scalar field $\Bbb{F},$ if a topological space $%
A, $ containing $\Bbb{F},$ satisfies

\strut

(i)\ \ \ $A$ is a vector space over $\Bbb{F}$ with its vector addition $(+)$ 
$:$ $(a_{1},$ $a_{2})$ $\mapsto $ $a_{1}$ $+$ $a_{2}$ and the scalar
multiplication $(\times )$ $:$ $(t,$ $a_{1})$ $\mapsto $ $ta_{1},$ for all $%
t $ $\in $ $\Bbb{F}$ and $a_{1},$ $a_{2}$ $\in $ $A$,

\strut

(ii)\ \ $A$ has a vector multiplication $(\cdot )$ $:$ $(a_{1},$ $a_{2})$ $%
\mapsto $ $a_{1}$ $a_{2},$ for all $a_{1},$ $a_{2}$ $\in $ $A,$ and it is
associative,

\strut

(iii) the vector addition $(+)$ and the vector multiplication $(\cdot )$ on $%
A$ are left and right distributive, and

\strut

(iv) the operations $(+),$ $(\cdot )$ and $(\times )$ are continuous, under
the topology for $A.$

\strut

Recall also that we say $A$ is an (algebraic) algebra over $\Bbb{F},$ if $A$
satisfies (i), (ii) and (iii). For example, $C^{*}$-algebras and von Neumann
algebras are topological algebras.\strut

\strut

\begin{definition}
Let $A_{1}$ and $A_{2}$ be algebraic algebras over the same scalar field $%
\Bbb{F}.$ Then we can define their algebraic free product algebra $A_{1}$ $%
*_{a\lg }$ $A_{2}$ of $A_{1}$ and $A_{2}$ by the algebra generated by all
words in $A_{1}$ $\cup $ $A_{2},$ in the algebra $A\lg _{\Bbb{F}}(A_{1},$ $%
A_{2})$ $=$ $\Bbb{F}[A_{1}$ $\cup $ $A_{2}],$ generated by $A_{1}$ and $%
A_{2}.$ Notice that $A_{1}$ $*_{a\lg }$ $A_{2}$ is not a topological space
(even though $A_{1}$ and $A_{2}$ are topological spaces).
\end{definition}

\strut

By the definition of algebraic free product and algebraic free product
algebras, we can get the following proposition.

\strut

\begin{proposition}
Let $A_{1}$ and $A_{2}$ be algebraic algebras over a scalar field $\Bbb{F}$,
and let $A_{1}$ $*_{a\lg }$ $A_{2}$ be the algebraic free product algebra of 
$A_{1}$ and $A_{2}.$ Assume now that $A_{k}$ $=$ $\Bbb{F}[X_{k}],$ where $%
X_{k}$ are the generator set of $A_{k},$ for $k$ $=$ $1,$ $2.$ Then

\strut 

\begin{center}
$A_{1}$ $*_{a\lg }$ $A_{2}$ $\overset{\text{Algebra}}{=}$ $\Bbb{F}[X_{1}$ $*$
$X_{2}]$,
\end{center}

\strut 

where ``$\overset{\text{Algebra}}{=}$'' means ``being algebra-isomorphic''.
\end{proposition}

\strut

\begin{proof}
Suppose $A_{k}$ $=$ $\Bbb{F}[X_{k}]$ are the algebraic algebras generated by 
$X_{k},$ for $k$ $=$ $1,$ $2.$ Then we have that

\strut

$\qquad A_{1}$ $*_{a\lg }$ $A_{2}$ $=$ $\Bbb{F}[A_{1}$ $*$ $A_{2}]$

\strut

by the definition of ``$*_{a\lg }$''

\strut

$\qquad \qquad =$ $\Bbb{F}\left[ \Bbb{F}[X_{1}]\text{ }*\text{ }\Bbb{F}%
[X_{2}]\right] $ $\overset{\text{Algebra}}{=}$ $\Bbb{F}\left[ \Bbb{F}[X_{1}%
\text{ }*\text{ }X_{2}]\right] $

\strut

by the definition of ``$*$''

\strut

$\qquad \qquad =$ $\Bbb{F}[X_{1}*X_{2}].$

\strut
\end{proof}

\strut

i.e., the above proposition shows that

\strut

(4.2)$\qquad \qquad \qquad \quad \Bbb{F}[X_{1}*X_{2}]$ $\overset{\text{%
Algebra}}{=}$ $\Bbb{F}[X_{1}]$ $*_{A\lg }$ $\Bbb{F}[X_{2}],$

\strut

where $X_{1}$ and $X_{2}$ are arbitrary sets.

\strut

Similar to the previous construction, we can determine algebraic \emph{%
reduced} free structures. Let $(X,$ $\cdot )$ be an arbitrary algebraic
pair. i.e., $(\cdot )$ $:$ $X$ $\times $ $X$ $\rightarrow $ $X$ is a binary
operation on $X.$ (Notice that $(X,$ $\cdot )$ is not necessarily be a
well-known algebraic structures, for instance, a semigroup or a group, or a
groupoid, etc.) Let $\mathcal{X}$ be the collection of all reduced words in $%
X,$ where the reduction is totally depending on the binary operation $(\cdot
)$ on $X.$

\strut

\begin{definition}
Let $(X,$ $\cdot )$ and $\mathcal{X}$ be given as in the previous paragraph.
Let $\mathcal{X}_{1}$ and $\mathcal{X}_{2}$ be subsets of $\mathcal{X}.$
Define the algebraic reduced free product set $\mathcal{X}_{1}$ $*^{r}$ $%
\mathcal{X}_{2}$ of $\mathcal{X}_{1}$ and $\mathcal{X}_{2},$ by the subset
of $\mathcal{X}$ consisting of all reduced words in $\mathcal{X}_{1}$ $\cup $
$\mathcal{X}_{2}.$
\end{definition}

\strut

Similarly, we can determine the algebraic reduced free product algebra of
algebraic algebras $A_{1}$ and $A_{2}$:

\strut

\begin{definition}
Let $A$ be an algebraic algebra over its scalar field $\Bbb{F},$ and let $%
A_{1}$ and $A_{2}$ be the algebraic subalgebras of $A.$ Define the algebraic
reduced free product algebra $A_{1}$ $*_{a\lg }^{r}$ $A_{2}$ of $A_{1}$ and $%
A_{2}$ in $A,$ by the algebra generated by all reduced words in $A_{1}$ $%
\cup $ $A_{2}.$ i.e.,

\strut 

\begin{center}
$A_{1}$ $*_{a\lg }^{r}$ $A_{2}$ $\overset{def}{=}$ $\Bbb{F}\left[ \{\text{%
all reduced words in }A_{1}\text{ }\cup \text{ }A_{2}\}\right] .$
\end{center}

\strut 

Remark that the reduction is dependent upon the vector multiplication on $A.$
\end{definition}

\strut

Then, similar to the previous proposition or (4.2), we can get the following
proposition.

\strut

\begin{proposition}
Let $A$ be an algebraic algebra over a scalar field $\Bbb{F}$ and assume
that $A$ $=$ $\Bbb{F}[X]$. i.e., $A$ is generated by a set $X.$ Let $A_{1}$
and $A_{2}$ be algebraic subalgebras of $A$ and suppose $A_{k}$ $=$ $\Bbb{F}%
[X_{k}],$ where $X_{k}$'s are the subset of $X.$ Then the algebraic reduced
free product algebra $A_{1}$ $*_{a\lg }^{r}$ $A_{2}$ in $A$ is $*$%
-isomorphic to the algebraic algebra $\Bbb{F}[X_{1}$ $*^{r}$ $X_{2}],$ where 
$X_{1}$ $*^{r}$ $X_{2}$ is the reduced free product of $X_{1}$ and $X_{2},$
under the vector multiplication on $A.$
\end{proposition}

\strut

\begin{proof}
First, notice that we can construct an algebraic pair $(X,$ $\cdot ),$ where 
$(\cdot )$ is the restricted vector multiplication on $A,$ and we can
understand $X$ as $(X,$ $\cdot ).$ (So, the algebraic reduced free product
set $X_{1}$ $*^{r}$ $X_{2}$ is well-defined.) Similar to the proof of the
previous proposition, we can get that

\strut

$\qquad A_{1}$ $*_{a\lg }^{r}$ $A_{2}$ $=$ $\Bbb{F}[A_{1}$ $*^{r}$ $A_{2}]$

\strut

by the definition of ``$*_{a\lg }^{r}$''

\strut \strut

$\qquad \qquad =$ $\Bbb{F}\left[ \Bbb{F}[X_{1}]\text{ }*^{r}\text{ }\Bbb{F}%
[X_{2}]\right] $ $\overset{\text{Algebra}}{=}$ $\Bbb{F}\left[ \Bbb{F}%
[X_{1}*^{r}X_{2}]\right] $

\strut

by the definition of ``$*^{r}$'' and by the vector addition and
multiplication on $A$

\strut

$\qquad \qquad =$ $\Bbb{F}[X_{1}$ $*^{r}$ $X_{2}].$

\strut
\end{proof}

\strut

So, the above proposition shows the relation between ``$*^{r}$'' and ``$%
*_{a\lg }^{r}$'':

\strut

(4.3)$\qquad \qquad \Bbb{F}\left[ X_{1}*^{r}X_{2}\right] $ $\overset{\text{%
Algebra}}{=}$ $\Bbb{F}[X_{1}]$ $*_{a\lg }^{r}$ $\Bbb{F}[X_{2}],$ inside $%
\Bbb{F}[X],$

\strut

whenever $X_{1}$ and $X_{2}$ are subsets of $X.$

\strut Now, we define the topological reduced free product algebras.

\strut

\begin{definition}
Let $A$ be a topological algebra over its scalar field $\Bbb{F},$ and let $%
A_{1}$ and $A_{2}$ be topological subalgebras of $A.$ Define the topological
reduced free product algebra $A_{1}$ $*_{top}^{r}$ $A_{2}$ of $A_{1}$ and $%
A_{2}$ in $A$ by the $\tau $-closure of $A_{1}$ $*_{a\lg }^{r}$ $A_{2},$
where $\tau $ is the topology for $A.$ i.e.,

\strut 

(4.4) $\qquad \qquad \quad A_{1}$ $*_{top}^{r}$ $A_{2}$ $\overset{def}{=}$ $%
\overline{A_{1}*_{a\lg }^{r}\text{ }A_{2}}^{\tau }$ inside $A.$
\end{definition}

\strut

By definition, we can get the following theorem.

\strut

\begin{theorem}
Let $A$ be a topological algebra over its scalar field $\Bbb{F},$ and let $%
A_{1}$ and $A_{2}$ be topological subalgebras of $A.$ Assume that $A$ $=$ $%
\overline{\Bbb{F}[X]}^{\tau }$ and $A_{k}$ $=$ $\overline{\Bbb{F}[X_{k}]}%
^{\tau },$ where $X$ is the generator set of $A$ and the $X_{k}$'s are
subsets of $X,$ for $k$ $=$ $1,$ $2.$ Then

\strut 

\begin{center}
$A_{1}$ $*_{top}^{r}$ $A_{2}$ $\overset{\text{Top-Algebra}}{=}$ $\overline{%
\Bbb{F}[X_{1}*^{r}X_{2}]}^{\tau },$
\end{center}

\strut 

where ``$\overset{\text{Top-Algebra}}{=}$'' means ``being topological
algebra isomorphic''.
\end{theorem}

\strut

\begin{proof}
Observe that

\strut

\begin{center}
$
\begin{array}{ll}
A_{1}*_{top}^{r}A_{2} & =\overline{\Bbb{F}[A_{1}\text{ }*_{a\lg }^{r}\text{ }%
A_{2}]}^{\tau }=\overline{\Bbb{F}\left[ \Bbb{F}[X_{1}]\text{ }*_{a\lg }^{r}%
\text{ }\Bbb{F}[X_{2}]\right] }^{\tau } \\ 
&  \\ 
& =\overline{\Bbb{F}\left[ \Bbb{F}[X_{1}\text{ }*^{r}\text{ }X_{2}]\right] }%
^{\tau }=\overline{\Bbb{F}\left[ X_{1}\text{ }*^{r}\text{ }X_{2}\right] }%
^{\tau },
\end{array}
$
\end{center}

\strut \strut

by (4.3).
\end{proof}

\strut

i.e., by the previous theorem, we have

\strut

(4.5)$\qquad \qquad \quad \overline{\Bbb{F}[X_{1}]}^{\tau }$ $*_{top}^{r}$ $%
\overline{\Bbb{F}[X_{2}]}^{\tau }$ $=$ $\overline{\Bbb{F}\left[ X_{1}\text{ }%
*^{r}\text{ }X_{2}\right] }^{\tau },$ in $\overline{\Bbb{F}[X]},$

\strut

whenever $X_{1},$ $X_{2}$ are subsets of $X.$

\strut

The following corollary is the direct consequence of the previous theorem.

\strut

\begin{corollary}
Let $\mathcal{G}$ be a finite family of partial isometries in $B(H)$ and let 
$\Bbb{G}$ be the $\mathcal{G}$-groupoid. Let $\{\Bbb{G}_{x}$ $:$ $x$ $\in $ $%
\mathcal{G}_{W}\}$ be the subgroupoids of the $\mathcal{G}$-groupoid $\Bbb{G}%
,$ where $\Bbb{G}_{x}$ are the graph groupoids induced by the corresponding
graphs $G_{x}$ of $x$ in $\mathcal{G}_{G},$ for all $x$ $\in $ $\mathcal{G}%
_{W}.$ Then

\strut 

(4.6)$\qquad \qquad \qquad C^{*}(\mathcal{G})$ $\overset{*\text{-isomorphic}%
}{=}$ $\underset{x\in \mathcal{G}_{W}}{\,\,*_{top}^{r}}$ $\left( C_{\alpha
}^{*}(\Bbb{G}_{x})\right) $ in $B(H_{\mathcal{G}})$,

\strut 

as embedded $C^{*}$-subalgebras of $B(H),$ where $(H_{\mathcal{G}},$ $\alpha
)$ is the $\mathcal{G}$-representation of $\Bbb{G}$.
\end{corollary}

\strut

\begin{proof}
Denote $C_{\alpha }^{*}(\Bbb{G})$ and $C_{\alpha }^{*}(\Bbb{G}_{x})$'s by $%
\mathcal{A}$ and $\mathcal{A}_{x}$'s, respectively.. By Section 4, we know
that the $C^{*}$-subalgebra $C^{*}(\mathcal{G})$ of $B(H)$ is $*$-isomorphic
to the groupoid $C^{*}$-algebra $\mathcal{A}$ in $B(H_{\mathcal{G}})$ $%
\subseteq $ $B(H),$ where $(H_{\mathcal{G}},$ $\alpha )$ is the $\mathcal{G}$%
-representation of $\Bbb{G}.$ It is easily checked that

\strut

\begin{center}
$\Bbb{G}$ $=$ $\underset{x\in \mathcal{G}_{W}}{\,*^{r}}$ $\Bbb{G}_{x},$
\end{center}

\strut

Thus, we can have that

\strut

$\qquad C^{*}(\mathcal{G})$ $\overset{*\text{-isomorphic}}{=}$ $\mathcal{A}$ 
$=$ $\overline{\Bbb{C}\left[ \alpha (\Bbb{G})\right] }$ $=$ $\overline{\Bbb{C%
}[\Bbb{G}]}$

\strut \strut

$\qquad \qquad \qquad =$ $\overline{\Bbb{C}\left[ \underset{x\in \mathcal{G}%
}{\,\,*^{r}}\text{ }\Bbb{G}_{x}\right] }$ $\overset{*\text{-isomorphic}}{=}$ 
$\underset{x\in \mathcal{G}_{W}}{\,\,*_{top}^{r}}$ $\left( \overline{\Bbb{C}[%
\Bbb{G}_{x}]}\right) $

by (4.5)

\strut $\qquad \qquad \qquad =$ $\underset{x\in \mathcal{G}_{W}}{%
\,\,*_{top}^{r}}$ $\mathcal{A}_{x}.$

\strut
\end{proof}

\strut

The above corollary shows that the $C^{*}$-algebra $C^{*}(\mathcal{G})$
generated by a finite family $\mathcal{G}$ of partial isometries in $B(H)$
has its block structure determined by the topological reduced free product,
and the blocks are $\mathcal{A}_{x}$ $=$ $C_{\alpha }^{*}(\Bbb{G}_{x})$'s,
for all $x$ $\in $ $\mathcal{G}_{W}.$ Since $\mathcal{A}_{x}$'s are
characterized in Section 3, we have the topological reduced free block
structures and characterized blocks.

\strut

\begin{remark}
Our topological reduced free (product) structure is basically different from
the (amalgamated reduced) free structures observed in [9], [10], [11] and
[13]. (Amalgamated) Reduced freeness on those papers are determined by the
free probabilistic settings of Voiculescu (See [16]). However, there are
some connections between them (See [13]).
\end{remark}

\strut

In the following section, we will consider this topological reduced free
block structures of $C^{*}(\mathcal{G})$ more in detail.

\strut \strut

\strut \strut

\strut

\subsection{Topological Free Block Structures on $C^{*}(\mathcal{G})$}

\strut

\strut

Let $\mathcal{G}$ be the given finite family of partial isometries in $B(H),$
and let $\mathcal{G}_{G}$ $=$ $\{G_{x}$ $:$ $x$ $\in $ $\mathcal{G}_{W}\}$
be the family of corresponding graphs induced by $\mathcal{G}_{W}.$ Also,
let $G$ be the $\mathcal{G}$-graph induced by $\mathcal{G}_{W},$ which is
the conditional iterated glued graph of $\mathcal{G}_{G},$ and let $\Bbb{G}$
be the $\mathcal{G}$-grouopoid, the graph groupoid of $G.$ In the previous
section, we showed that $\Bbb{G}$ is the topological reduced free product of 
$\Bbb{G}_{x}$'s, where $\Bbb{G}_{x}$ are the graph groupoids of $G_{x},$
which are the subgroupoid of $\Bbb{G},$ for all $x$ $\in $ $\mathcal{G}_{W},$
and

\strut

\begin{center}
$C^{*}\left( \mathcal{G}\right) $ $\overset{*\text{-isomorphic}}{=}C_{\alpha
}^{*}(\Bbb{G})$ $\overset{*\text{-isomorphic}}{=}$ $\underset{x\in \mathcal{G%
}_{W}}{\,\,*_{top}^{r}}$ $C_{\alpha }^{*}(\Bbb{G}_{x}).$
\end{center}

\strut

This means that the $C^{*}$-subalgebra $C^{*}(\mathcal{G})$ of $B(H)$ has
rough but characterized block structures $\{\mathcal{A}_{x}$ $:$ $x$ $\in $ $%
\mathcal{G}_{W}\}$, by Section 3.1. In this section, we will consider this
block structure more in detail, inside $C^{*}(\mathcal{G}).$

\strut

Suppose $u$ $\in $ $\mathcal{G}^{(u)}$ and $G_{u}$ $\in $ $\mathcal{G}_{G}.$
Then, by Subsection 3.2.1, the graph $G_{u}$ is the one-vertex-one-loop-edge
graph. Therefore, the graph groupoid $\Bbb{G}_{u}$ of $G_{u},$ which is a
subgroupoid of $\mathcal{G}$-groupoid $\Bbb{G},$ is a group. Also, if $s$ $%
\in $ $\mathcal{G}_{\infty }^{(s)},$ then the corresponding graph $G_{s}$ $%
\in $ $\mathcal{G}_{G}$ is a one-edge graph having distinct initial and
terminal vertices. So, the graph groupoid $\Bbb{G}_{s}$ of $G_{s}$ is a
finite (sub)groupoid consisting of the elements $0,$ $s,$ $s^{*}$, $s^{*}$ $%
s $ and $s$ $s^{*}.$ Therefore, if $x$ $\in $ $\mathcal{G}^{(u)}$ $\cup $ $%
\mathcal{G}_{\infty }^{(s)},$ then we can handle the corresponding
subgroupoid induced by $x$ relatively easy.

\strut

\begin{lemma}
Let $u_{1},$ $u_{2},$ $u$ $\in $ $\mathcal{G}^{(u)}$ and $s_{1},$ $s_{2},$ $s
$ $\in $ $\mathcal{G}_{\infty }^{(s)}.$

\strut 

(1) If $\pi (u_{1},$ $u_{2})$ $\neq $ $0$ in $\Bbb{G}_{0},$ then

\strut 

\begin{center}
$C^{*}(\{u_{1},$ $u_{2}\})$ $\overset{*\text{-isomorphic}}{=}$ $\overline{%
C^{*}(\{u_{1}\})*_{a\lg }C^{*}(\{u_{2}\})}.$
\end{center}

\strut 

Thus, $C^{*}(\{u_{1},$ $u_{2}\})$ is $*$-isomorphic to

\strut 

\begin{center}
$\overline{(\Bbb{C}\cdot 1_{H_{\{u_{1},\,u_{2}\}}})\otimes _{\Bbb{C}}\left(
C\left( spec(u_{1})\right) *_{a\lg }C\left( spec(u_{2})\right) \right) },$
\end{center}

\strut 

where $H_{\{u_{1},\text{ }u_{2}\}}$ is the $\{u_{1},$ $u_{2}\}$-Hilbert
space in the sense of Section 3.

\strut 

(2) if $\pi (s_{1},$ $s_{2})$ $\neq $ $0$ in $\Bbb{G}_{0},$ then

\strut 

\begin{center}
$C^{*}(\{s_{1},$ $s_{2}\})$ $\overset{*\text{-isomorphic}}{=}$ $\overline{%
C^{*}(\{s_{1}\})*_{a\lg }C^{*}(\{s_{2}\})}.$
\end{center}

\strut 

Thus, $C^{*}(\{s_{1},$ $s_{2}\})$ is $*$-isomorphic to

\strut 

\begin{center}
$(\Bbb{C}$ $\cdot $ $1_{H_{\{s_{1},\,s_{2}\}}})$ $\otimes _{\Bbb{C}}$ $%
\left( M_{2}(\Bbb{C})\text{ }\otimes _{\Bbb{C}}\text{ }M_{2}(\Bbb{C})\right) 
$ $=$ $(\Bbb{C}$ $\cdot $ $1_{H_{\{s_{1},\,s_{2}\}}})\otimes _{\Bbb{C}}$ $%
M_{4}(\Bbb{C}),$
\end{center}

\strut 

where $H_{\{s_{1},\,s_{2}\}}$ is the $\{s_{1},$ $s_{2}\}$-Hilbert space in
the sense of Section 3.

\strut 

(3) if $\pi (u,$ $s)$ $\neq $ $0$ in $\Bbb{G}_{0},$ then

\strut 

\begin{center}
$C^{*}(\{u,$ $s\})$ $\overset{*\text{-isomorphic}}{=}$ $\overline{%
C^{*}(\{u\})*_{a\lg }C^{*}(\{s\})}.$
\end{center}

\strut 

So, $C^{*}(\{u,$ $s\})$ is $*$-isomorphic to

\strut 

\begin{center}
$\overline{(\Bbb{C}\cdot 1_{H_{\{u,\,s\}}})\otimes _{\Bbb{C}}\left( C\left(
spec(u)\right) *_{a\lg }M_{2}(\Bbb{C})\right) }.$
\end{center}
\end{lemma}

\strut

\begin{proof}
By Section 4.1, we have that if $\mathcal{G}_{W}$ $=$ $\{x,$ $y\},$ then

\strut

\begin{center}
$
\begin{array}{ll}
C^{*}\left( \mathcal{G}_{W}\right) & \,\overset{*\text{-isomorphic}}{=}%
C_{\alpha }^{*}(\Bbb{G})=C_{\alpha }^{*}\left( \Bbb{G}_{x}*^{r}\Bbb{G}%
_{y}\right) \\ 
&  \\ 
& 
\begin{array}{l}
\overset{\ast \text{-isomorphic}}{=}C_{\alpha }^{*}(\Bbb{G}%
_{x})*_{top}^{r}C_{\alpha }^{*}(\Bbb{G}_{y}) \\ 
\\ 
\,\,\,\,\,\,\,\overset{def}{=}\text{ }\overline{C_{\alpha }^{*}(\Bbb{G}%
_{x})*_{a\lg }^{r}\text{ }C_{\alpha }^{*}(\Bbb{G}_{y})} \\ 
\\ 
\overset{\ast \text{-isomorphic}}{=}\overline{C^{*}(\{x\})*_{a\lg
}^{r}C^{*}(\{y\})},
\end{array}
\end{array}
$\strut
\end{center}

\strut \strut

where ``$*_{a\lg }^{r}$'' means the ``algebraic \textbf{reduced} algebra
free product''.

\strut

(1) Suppose $u_{1}$ and $u_{2}$ are unitary parts of certain partial
isometries in $B(H),$ and let $G_{u_{1}}$ and $G_{u_{2}}$ be the
corresponding one-vertex-one-loop-edge graphs in $\mathcal{G}_{G}.$ Then we
can construct the conditional glued graph $G$ $=$ $G_{u_{1}}$ $\#_{\pi }$ $%
G_{u_{2}}$ of $G_{u_{1}}$ and $G_{u_{2}},$ which is the $\mathcal{G}$-graph.
It is graph-isomorphic to the one-vertex-two-loop-edge graph. Then $\mathcal{%
G}$-groupoid $\Bbb{G}$ satisfies that $\Bbb{G}$ $=$ $\Bbb{G}_{u_{1}}$ $*^{r}$
$\Bbb{G}_{u_{2}},$ where $\Bbb{G}_{u_{k}}$ are the graph groupoids of $%
G_{u_{k}},$ for $k$ $=$ $1,$ $2.$ However, by the conditional gluing on the $%
\mathcal{G}$-graph $G,$ we can realize that $\Bbb{G}_{u_{1}}$ $*^{r}$ $\Bbb{G%
}_{u_{2}}$ is identified with $\Bbb{G}_{u_{1}}$ $*$ $\Bbb{G}_{u_{2}},$ where
``$*$'' means the ``algebraic (\textbf{non-reduced}) free product''.
Therefore,

\strut

\begin{center}
$
\begin{array}{ll}
C^{*}\left( \{u_{1},\text{ }u_{2}\}\right) & \overset{*\text{-isomorphic}}{=}%
\overline{C^{*}(\{u_{1}\})*_{a\lg }^{r}\text{ }C^{*}(\{u_{2}\})} \\ 
&  \\ 
& \overset{*\text{-isomorphic}}{=}\overline{C^{*}(\{u_{1}\})\text{ }*_{a\lg }%
\text{ }C^{*}(\{u_{2}\})},
\end{array}
$
\end{center}

\strut

where ``$*_{a\lg }$'' means the ``algebraic \emph{non-reduced} algebra free
product''.

\strut

(2) Assume now that $s_{1}$ and $s_{2}$ are shift parts in $\mathcal{G}%
_{\infty }^{(s)}$ of certain partial isometries in $B(H),$ and let $%
G_{s_{k}} $ be the corresponding graphs with their graph groupoids $\Bbb{G}%
_{s_{k}},$ for $k$ $=$ $1,$ $2.$ Suppose $\pi (s_{1},$ $s_{2})$ $\neq $ $0.$
Then we can create the conditional glued graph $G$ $=$ $G_{s_{1}}$ $\#_{\pi
} $ $G_{s_{2}}$ which is graph-isomorphic to

\strut

\begin{center}
$\bullet \longrightarrow \bullet \longrightarrow \bullet .$
\end{center}

\strut

Then we can easily check that the graph groupoid $\Bbb{G}$ of $G$ is $\{0,$ $%
s_{1},$ $s_{2},$ $s_{1}^{*},$ $s_{2}^{*},$ $s_{1}$ $s_{2},$ $s_{2}^{*}$ $%
s_{1}^{*}\},$ set-theoretically. So,

\strut

\begin{center}
$\Bbb{G}$ $=$ $\Bbb{G}_{s_{1}}$ $*^{r}$ $\Bbb{G}_{s_{2}}$ $=$ $\Bbb{G}%
_{s_{1}}$ $*$ $\Bbb{G}_{s_{2}}.$
\end{center}

\strut

Therefore, we have that

\strut

\begin{center}
$
\begin{array}{ll}
C^{*}(\{s_{1},\text{ }s_{2}\}) & \overset{*\text{-isomorphic}}{=}\text{ }%
\overline{C^{*}(\{s_{1}\})*_{a\lg }^{r}\text{ }C^{*}(\{s_{2}\})} \\ 
&  \\ 
& \overset{*\text{-isomorphic}}{=}\text{ }\overline{C^{*}(\{s_{1}\})\text{ }%
*_{a\lg }C^{*}(\{s_{2}\}).}
\end{array}
$
\end{center}

\strut

(3) Let $u$ $\in $ $\mathcal{G}^{(u)}$ and $s$ $\in $ $\mathcal{G}_{\infty
}^{(s)}$ with their corresponding graphs $G_{u},$ $G_{s}$ $\in $ $\mathcal{G}%
_{G},$ respectively. Then similar to the previous two cases, we can get the
desired result.
\end{proof}

\strut

The above lemma shows that it is relatively easy to deal with the groupoid $%
C^{*}$-algebraic structures in terms of their topological reduced free
blocks of them, since the topological reduced free product of blocks is
determined by the algebraic ``non-reduced'' free product of generator sets
of blocks.\ This means that we do not need to consider the reduction on
blocks induced by the elements in $\mathcal{G}^{(u)}$ $\cup $ $\mathcal{G}%
_{\infty }^{(s)},$ inside $C^{*}(\mathcal{G}).$

\strut \strut

However, if $x$ is contained in $\mathcal{G}_{f}^{(s)},$ then the
corresponding graph $G_{x}$ is graph-isomorphic to the infinite linear graph,

\strut

\begin{center}
$\underset{x^{*}x}{\bullet }\overset{x}{\longrightarrow }\underset{xx^{*}}{%
\bullet }\overset{x^{(2)}}{\longrightarrow }\underset{x^{2}x^{2\,*}}{\bullet 
}\overset{x^{(3)}}{\longrightarrow }\underset{x^{3}x^{3\,*}}{\bullet }%
\overset{x^{(4)}}{\longrightarrow }\cdots ,$
\end{center}

\strut \strut

and hence it is hard to handle the graph groupoid $\Bbb{G}_{x}$ of $G_{x},$
where $x^{(k)}$ $=$ $x$ on $H,$ for all $k$ $\in $ $\Bbb{N},$ with $x^{(1)}$ 
$=$ $x,$ satisfying

\strut

\begin{center}
$x$ $x^{(2)}$ ... $x^{(m)}$ $=$ $x^{m}$ on $H,$ for all $m$ $\in $ $\Bbb{N}.$
\end{center}

\strut

We will concentrate on the case where we have a fixed shift part $x$ in $%
\mathcal{G}_{f}^{(s)}.$ Let $x$ $\in $ $\mathcal{G}_{f}^{(s)}$ and $G_{x}$ $%
\in $ $\mathcal{G}_{G},$ and let $y$ $\in $ $\mathcal{G}_{W}$ with its
corresponding graph $G_{y}$ $\in $ $\mathcal{G}_{G}.$ First, assume that $y$ 
$\in $ $\mathcal{G}^{(u)}.$ If $\pi (x,$ $y)$ $=0$ in $\Bbb{G}_{0},$ then we
can have the disconnected graph $G_{x,y}$ $=$ $G_{x}$ $\#_{\pi }$ $G_{y}$ $=$
$G_{x}$ $\cup $ $G_{y},$

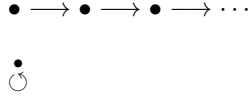
\strut

\begin{center}
$
\begin{array}{l}
\bullet \longrightarrow \bullet \longrightarrow \bullet \longrightarrow
\cdots \\ 
\\ 
\overset{\bullet }{\circlearrowleft }
\end{array}
$
\end{center}

\strut \strut \strut

If $\pi (x,$ $y)$ $\neq $ $0$ in $\Bbb{G}_{0},$ then $\pi (x^{n},$ $y)$ $%
\neq $ $0,$ for all $n$ $\in $ $\Bbb{N}$, since

\strut

\begin{center}
$x^{*}$ $x$ $\geq $ $x$ $x^{*}$ $\geq $ $x^{2}$ $x^{2\,*}\geq $ $x^{3}$ $%
x^{3\,*}\geq $ $x^{4}$ $x^{4\,*}\geq $ ...$.$
\end{center}

\strut

Thus, the conditional iterated glued graph $G_{x,y}$ is

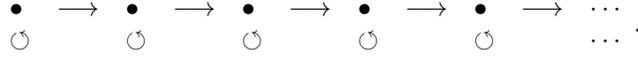
\strut

\begin{center}
$
\begin{array}{lllllllllll}
\bullet & \longrightarrow & \bullet & \longrightarrow & \bullet & 
\longrightarrow & \bullet & \longrightarrow & \bullet & \longrightarrow & 
\cdots \\ 
\circlearrowleft &  & \circlearrowleft &  & \circlearrowleft &  & 
\circlearrowleft &  & \circlearrowleft &  & \cdots
\end{array}
.$
\end{center}

\strut

If we consider the graph groupoid $\Bbb{G}_{x,y}$ of the conditional glued
graph $G_{x,y}$ of $G_{x}$ and $G_{y},$ then we can conclude that

$\strut $

(4.7)$\qquad \qquad \qquad \qquad \Bbb{G}_{x,y}$ $=$ $\Bbb{G}_{x}$ $*^{r}$ $%
\Bbb{G}_{y}$ $=$ $\Bbb{G}_{x}$ $*$ $\Bbb{G}_{y},$

\strut

where ``$*^{r}$'' means the reduced algebraic free product and ``$*$'' means
the non-reduced algebraic free product.

\strut

\begin{lemma}
Let $x$ $\in $ $\mathcal{G}_{f}^{(s)}$ and $u$ $\in $ $\mathcal{G}^{(u)},$
assume that $\pi (x,$ $u)$ $\neq $ $0$ in $\Bbb{G}_{0}.$ Then the $C^{*}$%
-subalgebra $C^{*}(\{x,$ $u\})$ of $B(H)$ is $*$-isomorphic to $\overline{%
C^{*}(\{x\})*_{a\lg }C^{*}(\{y\})}.$ Thus, $C^{*}(\{x,$ $u\})$ is $*$%
-isomorphic to

\strut 

\begin{center}
$\overline{(\Bbb{C}\cdot 1_{H_{\{x,\text{ }u\}}})\otimes _{\Bbb{C}}\left( 
\mathcal{T}(l^{2}(\Bbb{N}_{0}))\text{ }*_{a\lg }\text{ }C(spec(u))\right) }.$
\end{center}
\end{lemma}

\strut

\begin{proof}
It suffices to show that the graph groupoid $\Bbb{G}$ of the conditional
glued graph $G_{x,\,y}$ of $G_{x}$ and $G_{y}$ satisfies $\Bbb{G}$ $=$ $\Bbb{%
G}_{x}$ $*$ $\Bbb{G}_{y},$ where $\Bbb{G}_{x}$ and $\Bbb{G}_{y}$ are the
graph groupoids of the corresponding graphs $G_{x}$ and $G_{y}$ of $x$ and $%
y.$ By (4.7), it is done.
\end{proof}

\strut

Assume now that $x$ $\in $ $\mathcal{G}_{f}^{(s)}$ and $s$ $\in $ $\mathcal{G%
}_{\infty }^{(s)},$ and let $\pi (x,$ $s)$ $\neq $ $0$ in $\Bbb{G}_{0}.$
Then, $\pi (x^{n},$ $s)$ $\neq $ $0$ in $\Bbb{G}_{0},$ too, for all $n$ $\in 
$ $\Bbb{N},$ because

\strut

\begin{center}
$x^{*}$ $x$ $\geq $ $x$ $x^{*}$ $\geq $ $x^{2}$ $x^{2\,*}$ $\geq $ $x^{3}$ $%
x^{3\,*}$ $\geq $ ...,
\end{center}

\strut \strut

as projections on $H.$ Since $\pi (x,$ $s)^{*}$ $=$ $\pi (s^{*},$ $x^{*})$,
for any $x,$ $s$ $\in $ $\mathcal{G}_{W},$ Therefore, we can get the
conditional glued graph $G$ $=$ $G_{x}$ $\#_{\pi }$ $G_{s}$,
graph-isomorphic to the following graph,

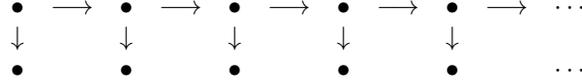
\strut

\begin{center}
$
\begin{array}{lllllllllll}
\bullet & \longrightarrow & \bullet & \longrightarrow & \bullet & 
\longrightarrow & \bullet & \longrightarrow & \bullet & \longrightarrow & 
\cdots \\ 
\downarrow &  & \downarrow &  & \downarrow &  & \downarrow &  & \downarrow & 
&  \\ 
\bullet &  & \bullet &  & \bullet &  & \bullet &  & \bullet &  & \cdots
\end{array}
$
\end{center}

\strut

The upper row of the above graph can be regarded as the graph $G_{x},$ and
the columns are regarded as the graphs $G_{s}$'s. This shows that, similar
to the previous lemma, we can conclude that

\strut

\begin{center}
$C^{*}(\{x,$ $s\})$ $\overset{*\text{-isomorphic}}{=}$ $\overline{%
C^{*}(\{x\})*_{a\lg }C^{*}(\{s\})}.$
\end{center}

\strut

\begin{lemma}
Let $x$ $\in $ $\mathcal{G}_{f}^{(s)},$ and $s$ $\in $ $\mathcal{G}_{\infty
}^{(s)}.$ Then the $C^{*}$-algebra $C^{*}(\{x,$ $s\})$ is $*$-isomorphic to
the $C^{*}$-closure $\overline{C^{*}(\{x\})*_{a\lg }C^{*}(\{s\})}$ of the
algebraic free product algebra of $C^{*}(\{x\})$ and $C^{*}(\{s\}).$ Thus, $%
C^{*}(\{x,$ $s\})$ is $*$-isomorphic to

\strut 

\begin{center}
$\overline{(\Bbb{C}\cdot 1_{H_{\{x,\,s\}}})\otimes _{\Bbb{C}}\left( \mathcal{%
T}(H_{x})\otimes _{\Bbb{C}}M_{2}(\Bbb{G})\right) }.$
\end{center}

$\square $
\end{lemma}

\strut

Finally, let's take $x,$ $y$ $\in $ $\mathcal{G}_{f}^{(s)},$ and assume that 
$\pi (x,$ $y)$ $\neq $ $0,$ in $\Bbb{G}_{0}.$ Let $H_{\{x,\,y\}}$ be the $%
\{x,$ $y\}$-Hilbert space in the sense of Section 3, induced by the
subspaces $H_{x}$ and $H_{y}$ of $H,$ and let $\alpha $ $\overset{denote}{=}$
$\alpha \mid _{\Bbb{G}_{x,\,y}}$ be the restriction of the $\mathcal{G}$%
-groupoid action $\alpha ,$ where $\Bbb{G}_{x,\,y}$ is the graph groupoid,
which is a subgroupoid of the $\mathcal{G}$-groupoid $\Bbb{G},$ induced by
the conditional glued graph $G_{x,\,y}$ $=$ $G_{x}$ $\#_{\pi }$ $G_{y}$ of $%
G_{x}$ and $G_{y}$ of $\mathcal{G}_{G}.$ We know that

\strut

$\qquad C^{*}(\{x,$ $y\})$ $\overset{*\text{-isomorphic}}{=}$ $C_{\alpha
}^{*}(\Bbb{G}_{x,\,y})$,

\strut

$\qquad \quad \qquad \overset{*\text{-isomorphic}}{=}$ $C_{\alpha }^{*}(\Bbb{%
G}_{x})$ $*_{top}^{r}$ $C_{\alpha }^{*}(\Bbb{G}_{y})$

\strut

by the definition of $\Bbb{G}_{x,\,y}$ (or $G_{x,\text{ }y}$) and by Section
4.2

\strut

$\qquad \qquad \quad \overset{*\text{-isomorphic}}{=}$ $\overline{C_{\alpha
}^{*}\left( \Bbb{G}_{x}\text{ }*^{r}\text{ }\Bbb{G}_{y}\right) }$

\strut

$\qquad \quad \qquad \overset{*\text{-isomorphic}}{=}$ $\overline{%
C^{*}(\{x\})*_{a\lg }^{r}\text{ }C^{*}(\{y\})}$ $\quad \left( \overset{*%
\text{-isomorphic}}{\neq }\overline{C^{*}(\{x\})*_{a\lg }C^{*}(y)}\right) $

\strut (4.8)

$\qquad \qquad \quad \overset{*\text{-isomorphic}}{=}$ $\mathcal{T}(H_{x})$ $%
*_{top}^{r}$ $\mathcal{T}(H_{y}),$

\strut

as a $C^{*}$-subalgebra of $B(H).$ By the reduction on ``$*_{top}^{r}$''
(and on ``$*_{a\lg }^{r}$'' and on ``$*^{r}$''), it is somewhat hard to deal
with it, compared with the previous cases. In the previous cases, we only
need to observe the non-reduced free products, moreover simple algebraic
(non-reduced) free products. However, in the above case, we need to observe
the reduction on the free product.

\strut

In fact, it is also interesting to observe the $*$-isomorphic $C^{*}$%
-subalgebra of $\mathcal{T}(H_{x})$ $*_{top}^{r}$ $\mathcal{T}(H_{y}),$ in
(4.8), inside $B(H).$ Clearly, if the corresponding graphs $G_{x}$ and $%
G_{y} $ in $\mathcal{G}_{G}$ satisfies either $G_{x}$ $>$ $G_{y}$ or $G_{x}$ 
$<$ $G_{y},$ or equivalently, if either

$\strut $

\begin{center}
$\dim \left( \ker y^{*}\right) $ $\mid $ $\dim (\ker x^{*})$ or $\dim (\ker
x^{*})$ $\mid $ $\dim (\ker y^{*}),$
\end{center}

\strut

then $G_{x,\,y}$ $=$ $G_{y},$ respectively, $G_{x,\,y}$ $=$ $G_{x}.$
Therefore, we can get that:

\strut

\begin{lemma}
Let $x,$ $y$ $\in $ $\mathcal{G}_{f}^{(s)},$ and assume that the
corresponding graphs $G_{x}$ and $G_{y}$ in $\mathcal{G}_{G}$ satisfy $G_{x}$
$>$ $G_{y}.$ Then $C^{*}$-algebra $C^{*}(\{x,$ $y\})$ generated by $\{x,$ $%
y\}$ is $*$-isomorphic to the Toepilitz algebra $\mathcal{T}(H_{y}).$ i.e.,

\strut 

\begin{center}
$C^{*}(\{x,$ $y\})\overset{*\text{-isomorphic}}{=}C_{\alpha }^{*}(\Bbb{G}%
_{y})\overset{*\text{-isomorphic}}{=}\mathcal{T}(H_{y}).$
\end{center}

$\square $
\end{lemma}

\strut

Assume now that $\pi (x,$ $y)$ $\neq $ $0$ in $\Bbb{G}_{0}$, and also assume
that neither $G_{x}$ $>$ $G_{y}$ nor $G_{x}$ $<$ $G_{y}.$ Then the
conditional glued graph $G_{x,\text{\thinspace }y}$ of $G_{x}$ and $G_{y}$
is graph-isomorphic to the following graph,

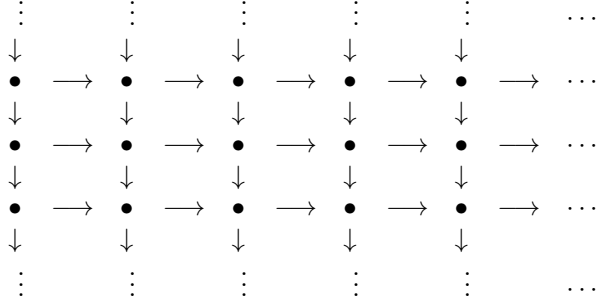
\strut

\begin{center}
$
\begin{array}{lllllllllll}
\,\,\vdots &  & \,\,\vdots &  & \,\,\vdots &  & \,\,\vdots &  & \,\,\vdots & 
& \cdots \\ 
\downarrow &  & \downarrow &  & \downarrow &  & \downarrow &  & \downarrow & 
&  \\ 
\bullet & \longrightarrow & \bullet & \longrightarrow & \bullet & 
\longrightarrow & \bullet & \longrightarrow & \bullet & \longrightarrow & 
\cdots \\ 
\downarrow &  & \downarrow &  & \downarrow &  & \downarrow &  & \downarrow & 
&  \\ 
\bullet & \longrightarrow & \bullet & \longrightarrow & \bullet & 
\longrightarrow & \bullet & \longrightarrow & \bullet & \longrightarrow & 
\cdots \\ 
\downarrow &  & \downarrow &  & \downarrow &  & \downarrow &  & \downarrow & 
&  \\ 
\bullet & \longrightarrow & \bullet & \longrightarrow & \bullet & 
\longrightarrow & \bullet & \longrightarrow & \bullet & \longrightarrow & 
\cdots \\ 
\downarrow &  & \downarrow &  & \downarrow &  & \downarrow &  & \downarrow & 
&  \\ 
\,\,\vdots &  & \,\,\vdots &  & \,\,\vdots &  & \,\,\vdots &  & \,\,\vdots & 
& \cdots
\end{array}
$
\end{center}

\strut

We can regard all edges of the previous graph as shifts $x$'s or $y$'s. Thus
all finite paths of the above graph are regarded as shifts, too. And hence,
all reduced finite paths in the graph groupoid $\Bbb{G}_{x,\,y}$ of $%
G_{x,\,y}$ are regarded as finite dimensional shifts on $H.$ So, it is easy
to check that

\strut

(4.9)$\qquad \qquad \qquad \mathcal{T}(H_{x})$ $*_{top}^{r}$ $\mathcal{T}%
(H_{y})$ $\overset{C^{*}\text{-subalgebra}}{\subseteq }$ $\mathcal{T}%
(H_{\{x,\,y\}}),$

\strut

in $B(H).$ Conversely, assume that $T$ $\in $ $\mathcal{T}(H_{\{x,\,y\}}).$
Since $\mathcal{T}(H_{\{x,\,y\}})$ is $*$-isomorphic to

$\strut $

\begin{center}
$C_{\alpha }^{*}(\Bbb{G}_{x,\,y})$ $\overset{*\text{-isomorphic}}{=}$ $%
C_{\alpha }^{*}(\Bbb{G}_{x})$ $*_{top}^{r}$ $C_{\alpha }^{*}(\Bbb{G}_{y})$ $%
\overset{*\text{-isomorphic}}{=}$ $\overline{\Bbb{C}[\Bbb{G}_{x}*^{r}\Bbb{G}%
_{y}]},$
\end{center}

\strut

the element $T$ is represented by

\strut

\begin{center}
$T$ $=$ $\underset{w\in \Bbb{G}_{x,\,y}}{\sum }$ $t_{w}$ $w,$ for $t_{w}$ $%
\in $ $\Bbb{C}$.
\end{center}

\strut

Since $\Bbb{G}_{x,\,y}$ $=$ $\Bbb{G}_{x}$ $*^{r}$ $\Bbb{G}_{y},$ the
operator $T$ is contained in $\overline{\Bbb{C}[\Bbb{G}_{x,\,y}]}$ in $%
B(H_{\{x,\,y\}})$ $\subseteq $ $B(H).$ This shows that

\strut

(4.10)$\qquad \qquad \qquad \mathcal{T}(H_{\{x,\,y\}})$ $\overset{C^{*}\text{%
-subalgebra}}{\subseteq }$ $\mathcal{T}(H_{x})$ $*_{top}^{r}$ $\mathcal{T}%
(H_{y}),$

\strut

inside $B(H),$ whenever $x$ and $y$ are finite dimensional shifts satisfying
that $\pi (x,$ $y)$ $\neq $ $0$ and neither $G_{x}$ $<$ $G_{y}$ nor $G_{x}$ $%
<$ $G_{y}.$

\strut

By (4.9) and (4.10), we can get the following lemma.

\strut

\begin{lemma}
Let $x,$ $y$ $\in $ $\mathcal{G}_{f}^{(s)}$ and assume that $\pi (x,$ $y)$ $%
\neq $ $0$ in $\Bbb{G}_{0},$ and neither $G_{x}$ $<$ $G_{y}$ nor $G_{y}$ $<$ 
$G_{x}$ in $\mathcal{G}_{G}.$ Then the $C^{*}$-algebra $C^{*}(\{x,$ $y\})$
is $*$-isomorphic to the Toeplitz algebra $\mathcal{T}(H_{\{x,\,y\}}),$
where $H_{\{x,\,y\}}$ is the $\{x,$ $y\}$-Hilbert space in the sense of
Section 3. $\square $
\end{lemma}

\strut

Readers can realize the difference between the above lemma and the previous
lemmas. The above lemma shows that, inductively, if we have a subfamily $%
\mathcal{G}_{f}$ of $\mathcal{G}_{f}^{(s)},$ satisfying that (i) $\pi (x_{1},
$ $x_{2})$ $\neq $ $0$ in $\Bbb{G}_{0},$ and (ii) neither $G_{x_{1}}$ $<$ $%
G_{x_{2}}$ nor $G_{x_{2}}$ $<$ $G_{x_{1}},$ in $\mathcal{G}_{G},$ for all
pair $(x_{1},$ $x_{2})$ $\in $ $\mathcal{G}_{f}$ $\times $ $\mathcal{G}_{f},$
then the $C^{*}$-subalgebra $C^{*}(\mathcal{G}_{f})$ of $B(H)$, generated by 
$\mathcal{G}_{f},$ is $*$-isomorphic to the Toeplitz algebra $\mathcal{T}(H_{%
\mathcal{G}_{f}}),$ where $H_{\mathcal{G}_{f}}$ is the $\mathcal{G}_{f}$%
-Hilbert space in the sense of Section 3, which is a subspace of $H_{%
\mathcal{G}}$ $\subseteq $ $H.$

\strut

\begin{proposition}
Let $\mathcal{G}_{f}$ be a subfamily of $\mathcal{G}_{f}^{(s)}$ $\subseteq $ 
$\mathcal{G}_{W}$ satisfying (i) $\pi (x_{1},$ $x_{2})$ $\neq $ $0$ in $\Bbb{%
G}_{0},$ and (ii) neither $G_{x_{1}}$ $<$ $G_{x_{2}}$ nor $G_{x_{1}}$ $>$ $%
G_{x_{2}}$ in $\mathcal{G}_{G},$ for all $(x_{1},$ $x_{2})$ $\in $ $\mathcal{%
G}_{f}\Bbb{\ \times }$ $\mathcal{G}_{f},$ then the $C^{*}$-subalgebra $C^{*}(%
\mathcal{G}_{f})$ of $B(H)$ is $*$-isomorphic to $\mathcal{T}(H_{\mathcal{G}%
_{f}}),$ where $H_{\mathcal{G}_{f}}$ is the $\mathcal{G}_{f}$-Hilbert space,
which is a subspace of $H_{\mathcal{G}}.$ $\square $
\end{proposition}

\strut

Let $\mathcal{G}_{f}^{(s)}$ $\subset $ $\mathcal{G}_{W}$ be nonempty. Define
the corresponding subset $\mathcal{G}_{G}^{(s:f)}$ of $\mathcal{G}_{G}$ by

\strut

\begin{center}
$\mathcal{G}_{G}^{(s:f)}$ $\overset{def}{=}$ $\{G_{x}$ $\in $ $\mathcal{G}%
_{G}$ $:$ $x$ $\in $ $\mathcal{G}_{f}^{(s)}\}$ $\subseteq $ $\mathcal{G}%
_{G}. $
\end{center}

\strut

Recall that the family $\mathcal{G}_{G}$ is in fact a partially ordered set $%
(\mathcal{G}_{G},$ $\leq ),$ where ``$\leq $'' is the full-subgraph
inclusion on directed graphs. Then, without loss of generality, we can
regard $\mathcal{G}_{G}^{(s:f)}$ as a partially ordered set $(\mathcal{G}%
_{G}^{(s:f)},$ $\leq ),$ under the inherited partial ordering $\leq $. By
this partial ordering, we can take the chains inside $\mathcal{G}%
_{G}^{(s:f)},$ and, for each chain, we can take the minimal elements in $%
\mathcal{G}_{G}^{(s:f)}.$

\strut

\begin{definition}
Let $\mathcal{G}_{G}^{(s:f)}$ $\subseteq $ $\mathcal{G}_{G}$ be given as
before. Define the subset $\mathcal{G}_{G}^{(s:\min )}$ of $\mathcal{G}%
_{G}^{(s:f)}$ $\subseteq $ $\mathcal{G}_{G}$ by

\strut 

\begin{center}
$\mathcal{G}_{G}^{(s:\min )}$ $\overset{def}{=}$ $\left\{ G\in \mathcal{G}%
_{G}^{(s:f)}:G\text{ is minimal in }\mathcal{G}_{G}^{(s:f)}\right\} .$
\end{center}

\strut 

By the collection $\mathcal{G}_{G}^{(s:\min )},$ we can take the
corresponding subset $\mathcal{G}_{f}^{(s:\min )}$ of $\mathcal{G}_{f}^{(s)}$
by

\strut 

\begin{center}
$\mathcal{G}_{f}^{(s:\min )}$ $\overset{def}{=}$ $\left\{ x\in \mathcal{G}%
_{f}^{(s)}:G_{x}\in \mathcal{G}_{G}^{(s:\min )}\right\} .$
\end{center}
\end{definition}

\strut

By the previous lemmas, we can realize that:

\strut

\begin{proposition}
The $C^{*}$-subalgebra $C^{*}\left( \mathcal{G}_{f}^{(s)}\right) $ of $C^{*}(%
\mathcal{G})$ is $*$-isomorphic to

\strut 

(4.11)$\qquad \qquad \qquad \qquad \quad \left( \Bbb{C}\cdot 1_{H_{\mathcal{G%
}_{f}^{(s)}}}\right) $ $\otimes _{\Bbb{C}}$ $C_{\alpha }^{*}\left( \Bbb{G}%
_{f}^{(s:\min )}\right) ,$

\strut 

where $\Bbb{G}_{f}^{(s:\min )}$ is the subgroupoid of the $\mathcal{G}$%
-groupoid $\Bbb{G}$, generated by the graph $G_{f}^{(s:\min )}$ $\overset{def%
}{=}$ $\underset{x\in \mathcal{G}_{f}^{(s:\max )}}{\#_{\pi }}$ $G_{x}.$ $%
\square $
\end{proposition}

\strut

Observe the graph $G_{f}^{(s:\min )}$ $\overset{def}{=}$ $\underset{x\in 
\mathcal{G}_{f}^{(s:\max )}}{\#_{\pi }}$ $G_{x}.$ By the previous lemmas, we
can conclude that if $x,$ $y$ $\in $ $\mathcal{G}_{f}^{(s:\min )}$ having
their corresponding graphs $G_{x},$ $G_{y}$ $\in $ $\mathcal{G}_{G}^{(s:\min
)},$ then there are two cases where $\pi (x,$ $y)$ $\neq $ $0$ or $\pi (x,$ $%
y)$ $=$ $0,$ in $\Bbb{G}_{0}.$ In particular, if $\pi (x,$ $y)$ $\neq $ $0,$
the full-subgraph $G_{x,\,y}$ $=$ $G_{x}$ $\#_{\pi }$ $G_{y}$ of $%
G_{f}^{(s:\min )}$ induces the graph groupoid $\Bbb{G}_{x,\,y}$ $=$ $\Bbb{G}%
_{x}$ $*^{r}$ $\Bbb{G}_{y},$ and $C_{\alpha }^{*}(\Bbb{G}_{x,\,y})$ is $*$%
-isomorphic to $\mathcal{T}(H_{\{x,\,y\}}).$ Clearly, if $\pi (x,$ $y)$ $=$ $%
0$, then

$\strut $

\begin{center}
$C_{\alpha }^{*}(\Bbb{G}_{x,\,y})$ $\overset{*\text{-isomorphic}}{=}$ $%
C_{\alpha }^{*}(\Bbb{G}_{x})$ $\oplus $ $C_{\alpha }^{*}(\Bbb{G}_{y})$ $%
\overset{*\text{-isomorphic}}{=}$ $\mathcal{T}(H_{x})$ $\oplus $ $\mathcal{T}%
(H_{y}).$
\end{center}

\strut

Therefore, we can get that:

\strut

\begin{proposition}
Let $\mathcal{G}_{f}^{(s:\min )}$ of $\mathcal{G}_{W}$ be given as above.
Then there exists $k$ $\in $ $\Bbb{N}$ such that (i) $k$ $\leq $ $\left| 
\mathcal{G}_{f}^{(s)}\right| ,$ and (ii) we can decompose $\mathcal{G}%
_{f}^{(s:\min )}$ by $\mathcal{G}_{f}^{(s:\min )}$ $=$ $\underset{\,m=1}{%
\overset{k}{\sqcup }}$ $\mathcal{G}_{f,\,m}^{(s:\min )},$ where $\sqcup $
means the disjoint union and where

$\strut $

\begin{center}
$x_{1}\neq $ $x_{2}$ $\in $ $\mathcal{G}_{f,\,m}^{(s:\min )}$ $%
\Longleftrightarrow $ $\pi (x_{1},$ $x_{2})$ $\neq $ $0\,$ or $\pi (x_{2},$ $%
x_{1})$ $\neq $ $0$ in $\Bbb{G}_{0}.$
\end{center}

\strut 

So, the $C^{*}$-algebra $C^{*}\left( \mathcal{G}_{f}^{(s:\min )}\right) $ is 
$*$-isomorphic to $\underset{m=1}{\overset{k}{\oplus }}$ $\mathcal{T}\left(
H_{\mathcal{G}_{f,m}^{(s:\min )}}\right) ,$ as $C^{*}$-subalgebras of $B(H_{%
\mathcal{G}})$ $\subseteq $ $B(H).$ $\square $
\end{proposition}

\strut

The proof of the previous proposition is done by the previous lemmas and
proposition, by induction.

\strut

\textbf{Notation} In the rest of this paper, we will keep using the same
notations used in the previous lemmas and propositions. $\square $

\strut

Now, we will combine all previous observations.\ The following theorem
provides the topological free block structures in $C^{*}(\mathcal{G}).$

\strut

\begin{theorem}
Let $\mathcal{G}$ be a finite family of partial isometries in $B(H).$ Define
the set $\mathcal{G}_{W}^{\min }$ of partial isometries induced by $\mathcal{%
G}_{W}$ of $\mathcal{G}$ by

\strut 

\begin{center}
$\mathcal{G}_{W}^{\min }$ $=$ $\mathcal{G}^{(u)}$ $\cup $ $\mathcal{G}%
_{\infty }^{(s)}$ $\cup $ $\mathcal{G}_{f}^{(s:\min )}$ $=$ $\mathcal{G}%
^{(u)}$ $\cup $ $\mathcal{G}_{\infty }^{(s)}$ $\cup $ $\left( \underset{m=1}{%
\overset{k}{\sqcup }}\text{ }\mathcal{G}_{f:m}^{(s:\min )}\right) .$
\end{center}

Then the $C^{*}$-subalgebra $C^{*}(\mathcal{G})$ of $B(H)$ is $*$-isomorphic
to

\strut 

\begin{center}
$\left( \Bbb{C}\cdot 1_{H_{\mathcal{G}}}\right) $ $\otimes _{\Bbb{C}}$ $%
C^{*}\left( \mathcal{G}_{W}^{\min }\right) ,$
\end{center}

\strut 

and hence it is $*$-isomorphic to

\strut 

$\qquad \overline{\left( \underset{u\in \mathcal{G}^{(u)}}{*_{a\lg }}\left( (%
\Bbb{C}\cdot 1_{H_{u}})\otimes _{\Bbb{C}}C(spec(u))\right) \right) }$

\begin{center}
\strut 

$\overline{*_{a\lg }\left( \underset{s\in \mathcal{G}_{\infty }^{(s)}}{%
*_{a\lg }}\left( (\Bbb{C}\cdot 1_{H_{s}})\otimes _{\Bbb{C}}M_{2}(\Bbb{C}%
)\right) \right) *_{a\lg }\left( \underset{m=1}{\overset{k}{\oplus }}\text{ }%
\mathcal{T}_{m}\right) },$
\end{center}

\strut 

where $\mathcal{T}_{m}$ $=$ $\mathcal{T}\left( H_{\mathcal{G}%
_{f,\,m}^{(s:\min )}}\right) $ are the Toeplitz algebras, for all $m$ $=$ $1,
$ ..., $k.$
\end{theorem}

\strut

\begin{proof}
Let $\mathcal{G}$ be a finite family of partial isometries in $B(H)$ and let 
$\mathcal{G}_{W}^{\min }$ be given as above. Then

\strut

$\quad C^{*}(\mathcal{G})$ $\overset{*\text{-isomorphic}}{=}$ $C_{\alpha
}^{*}(\Bbb{G})$ $\overset{*\text{-isomorphic}}{=}$ $\left( \Bbb{C}\cdot
1_{H_{\mathcal{G}}}\right) $ $\otimes _{\Bbb{C}}$ $C^{*}\left( \mathcal{G}%
_{W}^{\min }\right) $

\strut

by the previous lemmas and propositions

\strut

$\qquad \overset{*\text{-isomorphic}}{=}$ $\left( \Bbb{C}\cdot 1_{H_{%
\mathcal{G}}}\right) $ $\otimes _{\Bbb{C}}\left( C^{*}\left( \mathcal{G}%
^{(u)}\right) *_{top}^{r}C^{*}\left( \mathcal{G}_{\infty }^{(s)}\right)
*_{top}^{r}C^{*}\left( \underset{m=1}{\overset{k}{\sqcup }}\text{ }\mathcal{G%
}_{f,m}^{(s:\min )}\right) \right) $

\strut

$\qquad \overset{*\text{-isomorphic}}{=}$ $\left( \underset{u\in \mathcal{G}%
^{(u)}}{*_{top}^{r}}\text{ }\mathcal{A}_{u}\right) *_{top}^{r}\left( 
\underset{s\in \mathcal{G}_{\infty }^{(s)}}{*_{top}^{r}}\text{ }\mathcal{A}%
_{s}\right) *_{top}^{r}\left( \underset{m=1}{\overset{k}{\oplus }}\text{ }%
\mathcal{T}_{m}\right) $

\strut

where $\mathcal{A}_{x}$ $=$ $C^{*}(\{x\}),$ for all $x$ $\in $ $\mathcal{G}%
_{W}$, and $\mathcal{T}_{m}$ $=$ $\mathcal{T}\left( H_{\mathcal{G}%
_{f,\,m}^{(s:\min )}}\right) ,$ by the previous propositions, for all $m$ $=$
$1,$ ..., $k$

\strut

$\qquad \overset{*\text{-isomorphic}}{=}$ $\left( \overline{\underset{u\in 
\mathcal{G}^{(u)}}{*_{a\lg }}\mathcal{A}_{u}}\right) *_{top}^{r}$ $\left( 
\overline{\underset{s\in \mathcal{G}_{\infty }^{(s)}}{*_{a\lg }}\mathcal{A}%
_{s}}\right) $ $*_{top}^{r}$ $\left( \underset{m=1}{\overset{k}{\oplus }}%
\text{ }\mathcal{T}_{m}\right) $

\strut \strut \strut \strut

$\qquad \overset{*\text{-isomorphic}}{=}$ $\overline{\left( \underset{u\in 
\mathcal{G}^{(u)}}{*_{a\lg }}\mathcal{A}_{u}\right) *_{a\lg }\left( 
\underset{s\in \mathcal{G}_{\infty }^{(s)}}{*_{a\lg }}\mathcal{A}_{s}\right)
*_{a\lg }\left( \underset{m=1}{\overset{k}{\oplus }}\text{ }\mathcal{T}%
_{m}\right) },$

\strut

by the previous lemmas. Recall that

\strut

\begin{center}
$\mathcal{A}_{u}$ $\overset{*\text{-isomorphic}}{=}$ $(\Bbb{C}$ $\cdot $ $%
1_{H_{u}})$ $\otimes _{\Bbb{C}}$ $C\left( spec(u)\right) ,$ for all $u$ $\in 
$ $\mathcal{G}^{(u)}$,
\end{center}

and

\begin{center}
$\mathcal{A}_{s}$ $\overset{*\text{-isomorphic}}{=}$ $(\Bbb{C}$ $\cdot $ $%
1_{H_{s}})$ $\otimes _{\Bbb{C}}$ $M_{2}(\Bbb{C}),$ for all $s$ $\in $ $%
\mathcal{G}_{\infty }^{(s)}.$
\end{center}

\strut
\end{proof}

\strut \strut \strut \strut

The above theorem provides the characterization of the block structures of $%
C^{*}$-subalgebras generated by finitely many partial isometries in terms of
the algebraic (non-reduced) free product and the characterization of $C^{*}$%
-subalgebras generated by a single partial isometries.

\strut

\strut \strut

\subsection{Examples}

\strut

\strut \strut

In this section, we will consider some examples.

\strut

\begin{example}
In this first example, we will consider the graph-families we observed in
[13]. Let $\mathcal{G}$ $=$ $\{a_{1},$ ..., $a_{N}\}$ be a finite family of
partial isometries in $B(H).$ And assume that $\mathcal{G}$ constructs a
finite directed graph $G$ in the sense that

\strut 

(i)\ \ $\left| \mathcal{G}\right| $ $=$ $\left| E(G)\right| $ and $\left| 
\mathcal{G}_{p}\right| $ $=$ $\left| V(G)\right| ,$ where

\strut 

\begin{center}
$\mathcal{G}_{p}$ $\overset{def}{=}$ $\{a_{j}^{*}$ $a_{j}$ $:$ $j$ $=$ $1,$
..., $N\}$ $\cup $ $\{a_{j}$ $a_{j}^{*}$ $:$ $j$ $=$ $1,$ ..., $N\},$
\end{center}

\strut 

\qquad with the corresponding bijection $h_{E}$ $:$ $E(G)$ $\rightarrow $ $%
\mathcal{G}.$

\strut 

(ii) the edges $e_{1}$ and $e_{2}$ generate nonempty finite path $e_{1}$ $%
e_{2}$ if and only if $H_{init}^{h_{E}(e_{1})}$ $=$ $H_{fin}^{h_{E}(e_{2})},$
where the symbol ``$H_{1}$ $=$ $H_{2}$'' means that $H_{1}$ and $H_{2}$ are
identically same as subspaces of $H.$ So, ``$H_{1}$ $\neq $ $H_{2}$'' means
that $H_{1}$ $\cap $ $H_{2}$ $=$ $\{0_{H}\}.$

\strut 

The above graph-family-setting makes us understand the elements $a_{j}$
satisfies that either $a_{j}$ $=$ $u_{j}$ or $a_{j}$ $=$ $s_{j},$ where $%
u_{j}$ are the unitary parts of $a_{j}$ and $s_{j}$ are the shift parts of $%
a_{j}.$ Moreover, by the condition (ii), if $a_{j}$ $=$ $s_{j},$ for some $j$
$\in $ $\{1,$ ..., $N\},$ then its $*$-isomorphic index $i_{*}(a_{j})$ $=$ $%
(0,$ $\varepsilon ,$ $\infty ,$ $0),$ for some $\varepsilon $ $\in $ $\Bbb{N}%
_{0}^{\infty }.$ Thus we can get that the $C^{*}$-subalgebra $%
C^{*}(\{a_{j}\})$ generated by a single partial isometry $a_{j}$ $\in $ $%
\mathcal{G}$ is $*$-isomorphic to either

\strut 

\begin{center}
$(\Bbb{C}\cdot 1_{H_{a_{j}}})$ $\otimes _{\Bbb{C}}$ $C(\Bbb{T})$
\end{center}

\strut 

(if $h_{E}^{-1}(a_{j})$ is a loop-edge in $E(G)$) or

\strut 

\begin{center}
$(\Bbb{C}\cdot 1_{H_{a_{j}}})$ $\otimes _{\Bbb{C}}$ $M_{2}(\Bbb{C})$
\end{center}

\strut 

(if $h_{E}^{-1}(a_{j})$ is a non-loop edge in $E(G)$). i.e., the $\mathcal{G}%
_{W}$-admissibility is completely determined by (ii), and hence, we have
either,

\strut 

\begin{center}
$\pi (a_{i},$ $a_{j})$ $=$ $0$
\end{center}

or

\begin{center}
$\pi (a_{i},$ $a_{j})$ $=$ $(a_{i}^{*}$ $a_{i})$ $(a_{j}$ $a_{j}^{*})$ $=$ $%
a_{i}^{*}$ $a_{i}$ $=$ $a_{j}$ $a_{j}^{*},$
\end{center}

\strut 

for all $i,$ $j$ $\in $ $\{1,$ ..., $N\}.$ In other words, in this case, we
have that

\strut 

\begin{center}
$\mathcal{G}$ $=$ $\mathcal{G}_{W}$ $=$ $\mathcal{G}^{(u)}$ $\cup $ $%
\mathcal{G}_{\infty }^{(s)},$ with $\mathcal{G}_{f}^{(s)}$ $=$ $\varnothing .
$
\end{center}

\strut \strut 

Therefore, indeed, the $C^{*}$-algebras $C^{*}(\{a_{j}\})$ are $*$%
-isomorphic to either

\strut 

\begin{center}
$(\Bbb{C}\cdot 1_{H_{0}})$ $\otimes _{\Bbb{C}}$ $C(\Bbb{T})$ or $(\Bbb{C}$ $%
\cdot $ $1_{H_{0}})$ $\otimes _{\Bbb{C}}$ $M_{2}(\Bbb{C}),$
\end{center}

\strut 

by Section 3.1. Moreover, if $G$ is connected in the sense that, for any $%
(v_{1},$ $v_{2})$ $\in $ $V(\widehat{G})$ $\times $ $V(\widehat{G})$, there
always exists at least one reduced finite path $w$ $\in $ $FP_{r}(\widehat{G}%
)$ such that $w$ $=$ $v_{1}$ $w$ $v_{2}$ and $w^{-1}$ $=$ $v_{2}$ $w^{-1}$%
\thinspace $v_{1},$ then the $C^{*}$-subalgebra $C^{*}(\mathcal{G})$
generated by $\mathcal{G}$ is $*$-isomorphic to the affiliated matricial
graph $C^{*}$-algebra $\mathcal{M}_{G}(H_{0})$ which is $*$-isomorphic to

\strut 

\begin{center}
$(\Bbb{C}\cdot 1_{H_{0}})$ $\otimes _{\Bbb{C}}$ $\mathcal{M}_{G},$
\end{center}

\strut 

where $\mathcal{M}_{G}$ is the matricial graph $C^{*}$-algebra which is the $%
C^{*}$-subalgebra of $M_{n}(\Bbb{C})$ (See [13]), and where $H_{0}$ is the
Hilbert space that is Hilbert-space isomorphic to the initial and the final
subspaces of all elements in $\mathcal{G}.$ Notice that, indeed, $H_{0}$ $=$ 
$H_{\mathcal{G}}$ in $H$.

\strut 

So, in general, if the family $\mathcal{G}$ constructs a finite directed
graph $G$ in the above sense and if $G$ has its connected components $G_{1},$
..., $G_{t},$ for some $t$ $\in $ $\Bbb{N},$ where $G_{1},$ ..., $G_{t}$ are
the connected full-subgraphs of $G.$ Then we can decompose $\mathcal{G}$
into the disjoint union of $\mathcal{G}_{1},$ ..., $\mathcal{G}_{t},$
constructing $G_{1},$ ..., $G_{t},$ respectively. Then the $C^{*}$-algebra $%
C^{*}(\mathcal{G})$ is $\oplus _{j=1}^{t}$ $C^{*}(\mathcal{G}_{j})$ which is 
$*$-isomorphic to the direct sum of the affiliated matricial graph $C^{*}$%
-algebras:

\strut 

\begin{center}
$C^{*}(\mathcal{G})$ $=$ $\oplus _{j=1}^{t}$ $C^{*}(\mathcal{G}_{t})$ $%
\overset{*\text{-isomorphic}}{=}$ $\oplus _{j=1}^{t}$ $\mathcal{M}%
_{G}(H_{j}),$
\end{center}

\strut 

where $H_{j}$'s are the embedded subspaces of $H$ which are Hilbert-space
isomorphic to $p_{j}$ $H,$ for all $p_{j}$ $\in $ $C^{*}(\mathcal{G}%
_{p}^{(j)}),$ for all $j$ $=$ $1,$ ..., $t.$
\end{example}

\strut \strut \strut \strut

\begin{example}
Let $U$ be the unilateral shift on the Hilbert space $H$ $=$ $l^{2}(\Bbb{N}%
_{0}),$ and let $k_{1}$ $>$ $k_{2}$ $\in $ $\Bbb{N}.$ We can have the shift
operators $U^{k_{1}}$ and $U^{k_{2}}$ on $H.$ Then they generate the $C^{*}$%
-subalgebras $C^{*}(\{U^{k_{1}}\})$ and $C^{*}(\{U^{k_{2}}\})$ of $B(H),$
which are $*$-isomorphic to the classical Toeplitz algebra $\mathcal{T}(H).$
By hypothesis, $\pi (U^{k_{1}},$ $U^{k_{2}})$ $\neq $ $0$ on $H.$ Therefore, 
$C^{*}\left( \{U^{k_{1}},\text{ }U^{k_{2}}\}\right) $ is $*$-isomorphic to
the classical Toeplitz algebra $\mathcal{T}(H),$ by Section 4.2, since $%
H_{\{U^{k_{1}},\text{ }U^{k_{2}}\}}$ is Hilbert-space isomorphic to $H$ $=$ $%
l^{2}(\Bbb{N}_{0}).$ Indeed, if $\mathcal{G}_{W}$ $=$ $\{U^{k_{1}},$ $%
U^{k_{2}}\},$ then $\mathcal{G}_{f}^{(s:\max )}$ $=$ $\{U^{k_{2}}\}.$
Therefore, we can get the desired result.
\end{example}

\strut \strut \strut \strut \strut

\begin{example}
Let $k$ $\in $ $\Bbb{N}$ and $U,$ the unilateral shift on

\strut 

\begin{center}
$H$ $=$ $l^{2}(\Bbb{N}_{0})$ $\oplus $ $l^{2}(\Bbb{N}_{0})$ $\overset{denote%
}{=}$ $H_{1}$ $\oplus $ $H_{2}.$
\end{center}

\strut 

Define the operator $V$ on $H$ by

\strut 

\begin{center}
$V$ $\overset{def}{=}$ $\left( 
\begin{array}{ll}
0_{H} & 0_{H} \\ 
1_{H} & 0_{H}
\end{array}
\right) $ on $
\begin{array}{l}
H_{1} \\ 
\oplus  \\ 
H_{2}
\end{array}
.$
\end{center}

\strut 

Consider the family $\mathcal{G}$ $=$ $\{x_{1},$ $x_{2}\},$ where $x_{1}$ $=$
$U^{k}$ and $x_{2}$ $=$ $V.$ We can check that

\strut 

\begin{center}
\strut \strut $i_{*}(x_{1})$ $=$ $(0,$ $0,$ $k,$ $0)$ and $i_{*}(0,$ $0,$ $%
\infty ,$ $0)$ in $(\Bbb{N}_{0}^{\infty })^{4}.$
\end{center}

\strut 

Therefore, the $C^{*}$-subalgebras $\mathcal{A}_{x_{k}}$ $=$ $%
C^{*}(\{x_{k}\})$ of $B(H),$ for $k$ $=$ $1,$ $2,$ satisfy that

\strut 

\begin{center}
$\mathcal{A}_{x_{1}}$ $\overset{*\text{-isomorphic}}{=}$ $\mathcal{T}(H)$,
\end{center}

and

\begin{center}
$\mathcal{A}_{x_{2}}$ $\overset{*\text{-isomorphic}}{=}$ $\left( \Bbb{C}%
\cdot 1_{H_{1}}\right) $ $\otimes _{\Bbb{C}}$ $M_{2}(\Bbb{C}),$
\end{center}

\strut 

by Section 3.1. Also, we can have the family $\mathcal{G}_{G}$ $=$ $%
\{G_{x_{1}},$ $G_{x_{2}}\}$ of the corresponding graphs $G_{x_{k}}$ of $%
x_{k},$ for $k$ $=$ $1,$ $2.$ The graphs are

\strut 

\begin{center}
$G_{x_{1}}$ $=\quad \underset{x_{1}^{*}x_{1}}{\bullet }\overset{x_{1}}{%
\longrightarrow }\underset{x_{1}x_{1}^{*}}{\bullet }\overset{x_{1}^{(2)}}{%
\longrightarrow }\underset{x_{1}^{2}x_{1}^{\,2*}}{\bullet }\overset{%
x_{1}^{(3)}}{\longrightarrow }\underset{x_{1}^{3}x_{1}^{3*}}{\bullet }%
\overset{x_{1}^{(4)}}{\longrightarrow }\underset{x_{1}^{4}x_{1}^{\,4*}}{%
\bullet }\overset{x_{k}^{(5)}}{\longrightarrow }$
\end{center}

and

\begin{center}
$G_{x_{2}}$ $=$\quad $\underset{x_{2}^{*}x_{2}}{\bullet }\overset{x_{2}}{%
\longrightarrow }\underset{x_{2}x_{2}^{*}}{\bullet }.$
\end{center}

\strut 

We can check that

\strut 

\begin{center}
$\pi (x_{1}^{n},$ $x_{2}^{m})$ $=$ $(x_{1}^{n\,*}$ $x_{1}^{n})$ $(x_{2}^{m}$ 
$x_{2}^{m\,*})$ $=$ $\left\{ 
\begin{array}{ll}
x_{2}^{m}\text{ }x_{2}^{m\,*}\neq 0 & \text{if }m=1 \\ 
0 & \text{if }m\neq 1,
\end{array}
\right. $
\end{center}

\strut 

since $x_{2}^{k}$ $=$ $0,$ for all $k$ $\in $ $\Bbb{N}$ $\setminus $ $\{1\},$
and

\strut 

\begin{center}
$\pi (x_{1}^{n\,*},$ $x_{2}^{m\,*})$ $=$ $(x_{1}^{n}$ $x_{1}^{n\,*})$ $%
(x_{2}^{m\,*}$ $x_{2}^{m})$ $=$ $\left\{ 
\begin{array}{ll}
(1_{H})(1_{H_{1}}\oplus 0)=1_{H_{1}} & \text{if }m=1 \\ 
0 & \text{if }m\neq 1,
\end{array}
\right. $
\end{center}

\strut 

since $x_{2}^{k}$ $=$ $0,$ for all $k$ $\in $ $\Bbb{N}$ $\setminus $ $\{1\},$
for all $n,$ $m$ $\in $ $\Bbb{N}.$ Thus we can get the conditional iterated
glued graph $G$ $=$ $G_{x_{1}}$ $\#_{\pi }$ $G_{x_{2}}$ as the graph which
is graph-isomorphic to the following graph:

\strut 

\begin{center}
$G$ $=$ $\quad 
\begin{array}{lllllllllll}
\bullet  & \longrightarrow  & \bullet  & \longrightarrow  & \bullet  & 
\longrightarrow  & \bullet  & \longrightarrow  & \bullet  & \longrightarrow 
& \cdots  \\ 
\downarrow  &  & \downarrow  &  & \downarrow  &  & \downarrow  &  & 
\downarrow  &  &  \\ 
\bullet  &  & \bullet  &  & \bullet  &  & \bullet  &  & \bullet  &  & \cdots 
\end{array}
$
\end{center}

\strut 

By Section 4.2, the $C^{*}$-algebra $\mathcal{A}$ $=$ $C_{\alpha }^{*}(\Bbb{G%
})$ is $*$-isomorphic to

\strut 

\begin{center}
$\mathcal{A}$ $\overset{*\text{-isomorphic}}{=}$ $\overline{\left( \mathcal{T%
}(H)\right) *_{a\lg }\left( (\Bbb{C}\cdot 1_{H_{1}})\otimes _{\Bbb{C}}\text{ 
}M_{2}(\Bbb{C})\right) }.$
\end{center}

\strut \strut 
\end{example}

\strut

\begin{example}
Let $H$ $\overset{\text{Hilbert-Space}}{=}$ $l^{2}(\Bbb{N}_{0})$ having its
Hilbert basis $\mathcal{B}_{H}$ $=$ $\{\xi _{n}$ $:$ $n$ $\in $ $\Bbb{N}%
_{0}\},$ where

\strut 

\begin{center}
$\xi _{0}$ $=$ $(1,$ $0,$ $0,$ $0,$ ...$)$
\end{center}

\strut and

\begin{center}
$\xi _{n}$ $=$ $\left( \underset{n\text{-times}}{\underbrace{0,\text{ ..., }0%
}}\text{ , }1,\text{ }0,\text{ }0,\text{ ...}\right) ,$ for all $n$ $\in $ $%
\Bbb{N}.$
\end{center}

\strut 

Define an operator $x$ $\in $ $B(H)$ by the operator having its infinite
matricial form,

\strut 

\begin{center}
$\left( 
\begin{array}{llllllll}
0 & 0 & 0 & 0 & 0 & 0 & \cdots  &  \\ 
\frame{1} & 0 & 0 & 0 & 0 & 0 & \cdots  &  \\ 
0 & 0 & 0 & 0 & 0 & 0 & \cdots  &  \\ 
0 & \frame{1} & 0 & 0 & 0 & 0 & \cdots  &  \\ 
0 & 0 & 0 & 0 & 0 & 0 & \cdots  &  \\ 
0 & 0 & \frame{1} & 0 & 0 & 0 & \cdots  &  \\ 
&  &  & \,\,\vdots  &  &  & \ddots  &  \\ 
&  &  &  &  &  &  & 
\end{array}
\right) .$
\end{center}

\strut 

i.e., the operator $x$ satisfies that

\strut 

\begin{center}
$(t_{0},$ $t_{1},$ $t_{2},$ $t_{3},$ ....$)$ $\overset{x}{\longmapsto }$ $(0,
$ $t_{0},$ $0,$ $t_{1},$ $0,$ $t_{2},$ $0,$ $t_{3},$ $0,$ ...$),$
\end{center}

\strut 

for all $(t_{0},$ $t_{1},$ $t_{2},$ ...$)$ $\in $ $H.$ The adjoint $x^{*}$
of $x$ has its infinite matricial form,

\strut 

\begin{center}
$\left( 
\begin{array}{llllllll}
0 & 0 & \frame{1} & 0 & 0 & 0 & \cdots  &  \\ 
0 & 0 & 0 & \frame{1} & 0 & 0 & \cdots  &  \\ 
0 & 0 & 0 & 0 & \frame{1} & 0 & \cdots  &  \\ 
0 & 0 & 0 & 0 & 0 & 0 & \cdots  &  \\ 
0 & 0 & 0 & 0 & 0 & 0 & \cdots  &  \\ 
0 & 0 & 0 & 0 & 0 & 0 & \cdots  &  \\ 
&  &  & \,\,\vdots  &  &  & \ddots  &  \\ 
&  &  &  &  &  &  & 
\end{array}
\right) ,$
\end{center}

\strut 

so, it satisfies that

\strut 

\begin{center}
$(t_{0},$ $t_{1},$ $t_{2},$ $t_{3},$ $t_{4},$ $t_{5},$ ...$)$ $\overset{x^{*}%
}{\longmapsto }$ $(t_{1},$ $t_{3},$ $t_{5},$ ...$),$
\end{center}

\strut 

for all $(t_{0},$ $t_{1},$ $t_{2},$ ...$)$ $\in $ $H.$

\strut 

Instead of observing the Wold decomposition of $x,$ we will consider $x$ by
the operator determined by the (finite dimensional) shifts $%
(y_{n})_{n=0}^{\infty }$. First, we define the following subspaces $K_{n}$'s
of $H$:

\strut 

\begin{center}
$K_{n}$ $\overset{def}{=}$ $\Bbb{C}$ $\cdot $ $\xi _{n}$ $\overset{\text{%
Hilbert-Space}}{=}$ $\Bbb{C}$ $\overset{\text{Subspace}}{\subset }$ $H,$ for
all $n$ $\in $ $\Bbb{N}_{0}.$
\end{center}

\strut 

Now, define the operators $y_{n}$ by

\strut 

\begin{center}
$y_{0}$ $:$ $(t,$ $0,$ $0,$ $0,$ ...$)$ $\longmapsto $ $(0,$ $t,$ $0,$ $0,$ $%
0,$ ...$)$
\end{center}

\strut and

\begin{center}
$y_{n}$ $:$ $\left( \underset{n\text{-times}}{\underbrace{0,\text{ ...., }0}}%
,\text{ }t,\text{ }0,\text{ }0,\text{ ...}\right) $ $\mapsto $ $\left( 
\underset{(2n+1)\text{-times}}{\underbrace{0,\text{ ......., }0}},\text{ }t,%
\text{ }0,\text{ }0,\text{ ...}\right) $
\end{center}

\strut 

on $H,$ for all $n$ $\in $ $\Bbb{N}$ and $t$ $\in $ $\Bbb{C}.$ Then each
operator $y_{n}$ is a partial isometry on $H$ having its initial space and
its final space as follows:

\strut 

\begin{center}
$H_{init}^{y_{0}}$ $=$ $K_{0}$ and $H_{fin}^{y_{0}}$ $=$ $K_{1}$
\end{center}

and

\begin{center}
$H_{init}^{y_{n}}$ $=$ $K_{n}$ and $H_{fin}^{y_{n}}$ $=$ $K_{2n+1},$ for all 
$n$ $\in $ $\Bbb{N}.$
\end{center}

\strut 

Then the given operator $x$ is defined by

\strut 

\begin{center}
$x$ $\xi _{k}$ $\overset{def}{=}$ $y_{k}$ $\xi _{k},$ for all $k$ $\in $ $%
\Bbb{N}_{0},$ $\xi _{k}$ $\in $ $\mathcal{B}_{H}.$
\end{center}

\strut 

i.e., we can define

\strut 

\begin{center}
$x\left( \underset{k=0}{\overset{\infty }{\sum }}\text{ }t_{k}\text{ }\xi
_{k}\right) $ $\overset{def}{=}$ $\underset{k=0}{\overset{\infty }{\sum }}$ $%
t_{k}$ $y_{k}(\xi _{k})$ $=$ $\underset{k=0}{\overset{\infty }{\sum }}$ $%
t_{k}$ $\xi _{2k+1},$
\end{center}

\strut 

for all $\sum_{k=0}^{\infty }$ $t_{k}$ $\xi _{k}$ $\in $ $H,$ where $t_{k}$ $%
\in $ $\Bbb{C}.$ Now, construct the subspaces $(H_{n})_{n=0}^{\infty }$ of $H
$ by

\strut \strut 

\begin{center}
$H_{n}$ $=$ $K_{n}$ $\oplus $ $K_{2n+1},$ for all $n$ $\in $ $\Bbb{N}_{0}.$
\end{center}

\strut 

Note that $H_{n}$ $\overset{\text{Hilbert-Space}}{=}$ $\Bbb{C}^{\oplus \,2},$
for all $n$ $\in $ $\Bbb{N}_{0}.$ Then the operators $y_{n}$ can be
understood as shifts on $H_{n}.$ (Notice that they are regarded as shifts on
finite dimensional space $\Bbb{C}^{\oplus \,2}.$) Moreover, the operators $%
y_{n}$ on $H_{n}$ have their matricial forms,

\strut \strut 

\begin{center}
$y_{n}$ $=$ $\left( 
\begin{array}{ll}
0 & 0 \\ 
1 & 0
\end{array}
\right) $ on $H_{n}$ $=$ $
\begin{array}{l}
\,\,K_{n} \\ 
\,\oplus  \\ 
K_{2n+1}
\end{array}
,$
\end{center}

\strut 

for all $n$ $\in $ $\Bbb{N}_{0}.$ Therefore, the $C^{*}$-algebras $%
C^{*}(\{y_{n}\})$ generated by the shifts $y_{n}$ are $*$-isomorphic to $%
M_{2}(\Bbb{C})$ $=$ $B(H_{n}),$ as embedded $C^{*}$-subalgebras of $B(H),$
for all $n$ $\in $ $\Bbb{N}_{0}.$ Similar to Section 3.3, if we let $%
\mathcal{G}$ $=$ $\{y_{n}$ $\in $ $B(H)$ $:$ $n$ $\in $ $\Bbb{N}_{0}\}$ (by
regarding $y_{n}$'s as operators on $H$), then we get the $\mathcal{G}$%
-space $H_{\mathcal{G}}$ $=$ $H,$ and we can have the $\mathcal{G}$-action $%
\alpha $ of the $\mathcal{G}$-groupoid $\Bbb{G},$ where

\strut 

\begin{center}
$G_{y_{n}}$ $=\quad \underset{y_{n}^{*}y_{n}}{\bullet }\overset{y}{%
\longrightarrow }$ $\underset{y_{n}y_{n}^{*}}{\bullet }$ $\in \mathcal{G}%
_{G},$\quad for all $n$ $\in $ $\Bbb{N}_{0}.$
\end{center}

\strut 

Define a map $f$ $:$ $\Bbb{N}_{0}$ $\rightarrow $ $\Bbb{N}_{0}$ by

\strut 

\begin{center}
$f(m)$ $=$ $2m$ $+$ $1,$ for all $m$ $\in $ $\Bbb{N}_{0}.$
\end{center}

\strut 

We will denote the iterated compositions\ $\underset{k\text{-times}}{%
\underbrace{f\circ \text{......}\circ f}}$ \ by $f^{(k)},$ for all $k$ $\in $
$\Bbb{N}.$ Then we can have the subset $X_{(n)}$ of $\Bbb{N}$ by

\strut \strut 

\begin{center}
$X_{(n)}$ $\overset{def}{=}$ $\{n\}$ $\cup $ $\{f^{(k)}(n)$ $:$ $k$ $\in $ $%
\Bbb{N}\}$ $\subset $ $\Bbb{N}_{0}$,
\end{center}

\strut 

for all $n$ $\in $ $\Bbb{N}_{0}.$ It is easy to check that if $n$ $\in $ $%
\Bbb{N}_{0},$ then, for any $k$ $\in $ $\Bbb{N},$ the sets $X_{(f^{(k)}(n))}$
are contained in $X_{(n)}.$ i.e.,

\strut 

\begin{center}
$X_{(f^{(k)}(n))}$ $\subset $ $X_{(n)},$ for all $n$ $\in $ $\Bbb{N}_{0}$
and $k$ $\in $ $\Bbb{N}.$
\end{center}

\strut 

For instance, $X_{(7)}$ $\subset $ $X_{(3)}$ $\subset $ $X_{(1)}$ $\subset $ 
$X_{(0)},$ etc. This means that we can take chains under the usual
set-inclusion on the collection $\mathcal{Y}$ $=$ $\{X_{(n)}$ $:$ $n$ $\in $ 
$\Bbb{N}_{0}\}$ of $X_{(n)}$'s. Denote $\mathcal{X}$ by the collection of
all maximal elements of $\mathcal{Y},$ under the partial ordering $\subset $ 
$.$ For instance, $X_{(0)},$ $X_{(2)},$ $X_{(4)},$ $X_{(6)}$ are the first
five elements of $\mathcal{X}.$ Let $I$ be the subset of $\Bbb{N}_{0}$
defined by

\strut 

\begin{center}
$I$ $=$ $\{n$ $\in $ $\Bbb{N}_{0}$ $:$ $X_{(n)}$ $\in $ $\mathcal{X}\}$ $%
\subset $ $\Bbb{N}_{0}.$
\end{center}

\strut \strut 

Let $\mathcal{G}$ $=$ $\{y_{0},$ $y_{1},$ $y_{2},$ $y_{3},$ ...$\}$ and let $%
\mathcal{G}_{G}$ $=$ $\{G_{(n)}$ $:$ $n$ $\in $ $I\}.$ Then we can define
the conditional iterated glued graph

\strut 

\begin{center}
$G$ $=$ $\underset{n\in I}{\,\#_{\pi }}$ $G_{(n)},$
\end{center}

\strut where

\begin{center}
$
\begin{array}{ll}
G_{(n)} & =\text{ }G_{n}\text{ }\#_{\pi }\left( \underset{k=n}{\overset{%
\infty }{\,\#_{\pi }}}G_{y_{2k+1}}\right)  \\ 
&  \\ 
& =\quad \underset{y_{n}^{*}y_{n}}{\bullet }\overset{y_{n}}{\longrightarrow }%
\underset{y_{n}y_{n}^{*}=y_{2n+1}^{*}y_{2n+1}}{\bullet }\overset{y_{2n+1}}{%
\longrightarrow }\bullet \overset{y_{2(2n+1)+1}}{\longrightarrow }\bullet 
\text{ }\cdots ,
\end{array}
$
\end{center}

\strut 

for all $n$ $\in $ $I.$ Clearly, we can create the corresponding graph
groupoid $\Bbb{G}$ of $G.$ Then we can determine $H_{\mathcal{G}}$ $=$ $H$
and the groupoid action $\alpha $ of $\Bbb{G}$ on $H,$ like Section 3.3. We
can check that

\strut 

\begin{center}
$C^{*}(\mathcal{G})$ $\overset{*\text{-isomorphic}}{=}$ $C_{\alpha }^{*}(%
\Bbb{G}),$
\end{center}

\strut 

by Section 3.4. By Section 4, we can check that the graph groupoid $\Bbb{G}%
_{(n)}$ of $G_{(n)}$ generates the $C^{*}$-algebra $C_{\alpha }^{*}(\Bbb{G}%
_{(n)})$ satisfies that

\strut 

\begin{center}
$C_{\alpha }^{*}(\Bbb{G}_{(n)})$ $\overset{*\text{-isomorphic}}{=}$ $%
\underset{k=1}{\overset{\infty }{*_{top}^{r}}}$ $\mathcal{M}_{k},$ with $%
\mathcal{M}_{k}$ $=$ $B(\Bbb{C}\xi _{k}\oplus \Bbb{C}\xi _{2k+1}),$
\end{center}

\strut 

for all $k$ $\in $ $\Bbb{N},$ and $n$ $\in $ $I.$ Notice that $\mathcal{M}%
_{k}$ $\overset{*\text{-isomorphic}}{=}$ $M_{2}(\Bbb{C}),$ for all $k$ $\in $
$\Bbb{N}.$
\end{example}

\strut \strut \strut \strut

\strut

\strut

\section{\strut Extensions of Unbounded Operators}

\strut

\strut

\strut

In this Section, we will consider an application of our results to Unbounded
Operator Theory, in particular, the Cayley Transform Theory. Let $H$ be a
countably infinite dimensional Hilbert space and let $\mathcal{D}$ be a
dense subspace of $H.$ A linear unbounded operator $S$ defined on $D$ is
said to be \emph{Hermitian} if

\strut

(6.1)$\qquad \qquad \qquad \qquad <S$ $\xi _{1},$ $\xi _{2}>$ $=$ $<\xi
_{1}, $ $S$ $\xi _{2}>$

\strut

holds for all $\xi _{1},$ $\xi _{2}$ $\in $ $\mathcal{D}.$ Notice that, the
relation (6.1) means that the dense subspace $\mathcal{D}$ is contained in
the domain of the adjoint operator $S^{*}$ of $S,$ and that

\strut

\begin{center}
$S$ $\xi $ $=$ $S^{*}$ $\xi ,$ for all $\xi $ $\in $ $\mathcal{D}.$
\end{center}

\strut

In other words, the graph $\mathcal{G}(S)$ of $S$ is contained in $\mathcal{G%
}(S^{*})$ of $S^{*}.$ Here, the graph $\mathcal{G}(S)$ of $S$ means that the
subset $\{(\xi ,$ $S$ $\xi )$ $:$ $\xi $ $\in $ $\mathcal{D}\}$ of $H$ $%
\times $ $H.$ We denote this relation by

\strut

\begin{center}
$S$ $\subseteq $ $S^{*}.$
\end{center}

We are interested in finding Hermitian elements $T$ of $S$; i.e., Hermitian
operators $T$ such that

\strut

(6.2)$\qquad \qquad \qquad \qquad \qquad \quad \mathcal{G}(S)$ $\subseteq $ $%
\mathcal{G}(T)$.

\strut

\begin{definition}
The Cayley transform $V$ of $S$ is a unbounded operator

\strut 

(6.3)$\qquad \qquad \qquad \qquad \quad V$ $\overset{def}{=}$ $\left(
S+i\right) (S-i)^{-1}.$

i.e.,

(6.4)$\quad \qquad \qquad \quad V\left( S\text{ }\xi -i\xi \right) $ $=$ $S$ 
$\xi $ $+$ $i\xi ,$ for all $\xi $ $\in $ $\mathcal{D},$

\strut 

where $i$ of (6.3) means $i1_{H},$ with the imaginary number $i$ satisfying $%
i^{2}$ $=$ $-1$ in $\Bbb{C},$ where $1_{H}$ means the identity operator on $%
H.$
\end{definition}

\strut

It is immediate that the given unbounded operator $S$ is Hermitian if and
only if the Cayley transform $V$ is a partial isometry. Moreover, the
orthogonal compliment $\mathcal{E}_{+}$ of the operator $V^{*}$ $V$ satisfies

\strut (6.5)

\begin{center}
$
\begin{array}{ll}
\mathcal{E}_{+} & =\ker \left( S^{*}-i\right) \\ 
&  \\ 
& =\left\{ \xi _{+}\in dom(S^{*})\left| S^{*}\xi _{+}=i\xi _{+}\right.
\right\} ,
\end{array}
$
\end{center}

\strut

where $dom(T)$ means the domain of an operator $T.$ Also, the orthogonal
compliment $\mathcal{E}_{-}$ of $V$ $V^{*}$ satisfies

\strut

(6.6)

\begin{center}
$
\begin{array}{ll}
\mathcal{E}_{-} & =\ker \left( S^{*}+i\right) \\ 
&  \\ 
& =\left\{ \xi _{-}\in dom(S^{*})\left| S^{*}\xi _{-}=-i\text{ }\xi
_{-}\right. \right\} .
\end{array}
$
\end{center}

\strut

\begin{definition}
The subspaces $\mathcal{E}_{+}$ and $\mathcal{E}_{-}$ are called the defect
spaces or the deficiency space of $V.$ And the numbers

\strut 

(6.7)$\qquad \qquad \qquad \quad n_{+}$ $\overset{def}{=}$ $\dim \mathcal{E}%
_{+}$ and $n_{-}$ $\overset{def}{=}$ $\dim \mathcal{E}_{-}$

\strut 

are called the deficiency indices of $V.$
\end{definition}

\strut

\begin{lemma}
(See [48]) A Hermitian symmetric operator $S$ has its self-adjoint
extensions $T$ if and only if $n_{+}$ $=$ $n_{-}.$ In this case, all
self-adjoint extensions $T$ satisfy that

\strut 

(6.8) $\qquad \qquad \qquad \qquad \qquad S$ $\subseteq $ $T$ $\subseteq $ $%
S^{*}.$

\strut 

Moreover, there is an one-to-one correspondence between all the unitary
extensions $U$ of the Cayley transform $V$ and all the self-adjoint
extensions $T$ of $S$:

\strut 

\quad (i) If $U$ is a unitary extension of $V,$ i.e., $\mathcal{G}(V)$ $%
\subseteq $ $\mathcal{G}(U),$ then

\strut 

(6.9)$\qquad \qquad \qquad \qquad T$ $\overset{def}{=}$ $i\left(
U+1_{H}\right) \left( U-1_{H}\right) ^{-1}$

\strut 

\quad is a self-adjoint extension of $S.$

\strut 

\quad (ii) If conversely $T$ is a self-adjoint extension of $S,$ then

\strut 

(6.10)$\qquad \qquad \qquad \qquad U$ $=$ $\left( T+i\right) \left(
T-i\right) ^{-1}$

\strut 

\quad is a unitary extension of $V.$

\strut 

\quad (iii) In general, the correspondence between (i) and (ii) holds for
the family of all Hermitian extensions $V_{T}$ of $V,$ and all the Hermitian
extensions $T$ of $S$; and the correspondence is decided by

\strut 

(611)$\qquad \qquad \qquad \quad T$ $=$ $i\left( V_{T}+1_{H}\right) \left(
V_{T}-1_{H}\right) ^{-1}$

and

(6.12)$\qquad \qquad \qquad \qquad V_{T}=\left( T+i\right) (T-i)^{-1},$

\strut 

respectively. $\square $
\end{lemma}

\strut

We provide several examples to illustrate our purpose.

\begin{example}
(1) Let $H$ $=$ $L^{2}(0,$ $1)$ and let

\strut 

(6.13)$\qquad \qquad \qquad \qquad \qquad S$ $\overset{def}{=}$ $\frac{1}{i}$
$\frac{d}{dx}$

with

\begin{center}
$\mathcal{D}$ $=$ $dom(S)$ $=$ $C_{0}^{1}(0,$ $1)$ $=$ all $C^{1}$-functions 
$\varphi $ on $[0,$ $1]$
\end{center}

such that

\begin{center}
$\varphi (0)$ $=$ $\varphi (1)$ $=$ $0.$
\end{center}

\strut 

Then the operator $S$ has its deficiency indices

\strut 

\begin{center}
$n_{+}$ $=$ $1$ $=$ $n_{-}$
\end{center}

\strut 

and the defect spaces are

\strut 

(6.14)$\qquad \qquad \qquad \qquad \qquad \mathcal{E}_{+}$ $=$ $\Bbb{C}\cdot
e^{-x}$

and

(6.15)$\qquad \qquad \qquad \qquad \qquad \mathcal{E}_{-}$ $=$ $\Bbb{C}\cdot
e^{x}.$

\strut 

(2) Every $\theta $ $\in $ $[0,$ $2\pi )$ determines a unique self-adjoint
extension $T_{\theta }$ of $S$ as follows. Set

\strut 

\begin{center}
$\mathcal{B}_{\theta }$ $\overset{def}{=}$ $\left\{ e^{i(\theta +2\pi n)x}%
\text{ }:\text{ }n\in \Bbb{Z}\right\} .$
\end{center}

Then

\begin{center}
$\mathcal{B}_{\theta }$ $\subseteq $ $dom\left( T_{\theta }\right) $
\end{center}

and

\begin{center}
$T_{\theta }$ $e^{i(\theta +2\pi \,n)x}$ $=$ $\left( \theta +2\pi \,n\right) 
$ $e^{i\left( \theta +2\pi \,n\right) x},$ for $n$ $\in $ $\Bbb{Z}.$
\end{center}

i.e.,

(6.16)$\qquad \qquad \qquad \quad spec\left( T_{\theta }\right) =\theta
+2\pi $ $\Bbb{Z}$ in $\Bbb{C},$

\strut 

where $spec(T_{\theta })$ means the spectrum of $T_{\theta }.$
\end{example}

\strut

The following proposition is the main result of this Section.

\strut

\begin{proposition}
Let $S$ be a Hermitian symmetric operator with its dense domain in a Hilbert
space $H,$ and assume that $S$ has the deficiency indices $(1,$ $1).$ Then
all the self-adjoint extensions of $S$ have Cayley transform

\strut 

(6.17) $\qquad \qquad \qquad \qquad \qquad U$ $=$ $V$ $\oplus $ $W,$

\strut 

where $V$ is the Cayley transform of $S$ determined in (6.3), and where $W$
has $W^{*}$ $W$ and $W$ $W^{*}$ of rank $1.$ Moreover, there exists $\alpha $
$\in $ $\Bbb{C}$, $\left| \alpha \right| $ $<$ $1,$ such that

\strut 

(6.18)$\qquad \qquad \qquad \qquad \qquad \quad W^{2}$ $=$ $\alpha $ $W.$
\end{proposition}

\strut

\begin{proof}
The defect spaces $\mathcal{E}_{+}$ and $\mathcal{E}_{-}$ in (6.5) and (6.6)
are the initial and final spaces of the partial isometries $W$ from the
previous lemma which extend the Cayley transform $V$ of $S.$ So,

\strut $\strut $

\begin{center}
$W^{*}$ $W$ $=$ $\mathcal{E}_{+}$ and $W$ $W^{*}$ $=$ $\mathcal{E}_{-}.$
\end{center}

\strut

Since $\dim \mathcal{E}_{\pm 1}$ $=$ $1,$ there are vectors $e_{\pm 1}$ $\in 
$ $\mathcal{D}$ such that $\left\| e_{\pm 1}\right\| $ $=$ $1$ and

\strut

\begin{center}
$\mathcal{E}_{+}$ $=$ $\Bbb{C}$ $\cdot $ $e_{+}$ and $\mathcal{E}_{-}$ $=$ $%
\Bbb{C}\cdot e_{-}.$
\end{center}

\strut \strut

It follows that

\strut

(6.19)$\qquad \qquad \qquad \qquad \qquad W$ $=$ $\left| e_{-}><e_{+}\right|
,$

\strut

where we use the Dirac notation,

\strut

(6.20)

\begin{center}
$\left| \xi ><\eta \right| \zeta $ $\overset{def}{=}$ $<\eta ,$ $\zeta >$ $%
\xi ,$ for all $\zeta ,$ $\eta ,$ $\xi $ $\in $ $H.$
\end{center}

\strut

Hence we have that

\strut

\begin{center}
$W^{2}$ $=$ $<e_{+},$ $e_{-}>$ $W$
\end{center}

and

\begin{center}
$\alpha $ $=$ $<e_{+},$ $e_{-}>$
\end{center}

works in (6.18).

\strut

If $\left| \alpha \right| $ $<$ $1,$ then, by Schuarz, it would follow that $%
\mathcal{E}_{+}$ $=$ $\mathcal{E}_{-}$, which contradicts

\strut

\begin{center}
$S^{*}$ $e_{\pm 1}$ $=$ $\pm i\,e_{\pm 1},$
\end{center}

\strut

from (6.5) to (6.6).\strut
\end{proof}

\strut

By the previous proposition, we can get the following corollary.

\strut

\begin{corollary}
Let $W$ be as in the previous proposition. Then $W^{n}$ $\rightarrow $ $0,$
as $n$ $\rightarrow $ $\infty .$
\end{corollary}

\strut

\begin{proof}
By (6.18) and by Induction, we have that

\strut

\begin{center}
$W^{n+1}$ $=$ $\alpha ^{n}$ $W,$ for all $n$ $\in $ $\Bbb{N}.$
\end{center}

\strut

Thus $W^{n}$ $\rightarrow $ $0,$ as $n$ $\rightarrow $ $\infty $.
\end{proof}

\strut \strut

\strut

\strut

\strut \textbf{References}

\strut

\begin{quotation}
\strut

{\small [1] \ \ A. G. Myasnikov and V. Shapilrain (editors), Group Theory,
Statistics and Cryptography, Contemporary Math, 360, (2003) AMS.}

{\small [2] \ }$\ ${\small B. Solel, You can see the arrows in a Quiver
Operator Algebras, (2000), preprint.\strut }

{\small [3] \ \ \strut C. W. Marshall, Applied Graph Theory, ISBN:
0-471-57300-0 (1971) John Wiley \& Sons}

{\small [4] \ \ D. A. Lind, Entropies of Automorphisms of a Topological
Markov Shift, Proc. AMS, vo 99, no 3, (1987) 589 - 595.}

{\small [5]\ \ \ D. A. Lind and B. Marcus, An Introduction to Symbolic
Dynamics and Coding, (1995) Cambridge Univ. Press.}

{\small [6]\ \ \ D.W. Kribs, Quantum Causal Histories and the Directed Graph
Operator Framework, arXiv:math.OA/0501087v1 (2005), Preprint.}

{\small [7]\ \ \ G. C. Bell, Growth of the Asymptotic Dimension Function for
Groups, (2005) Preprint.}

{\small [8]\ \ \ I. Cho, Direct Producted }$C^{*}$-{\small Probabiliy Spaces
and Corresponding Amalgamated Free Stochastic Integration, B. of KMS, 44,
(2007), 131 - 150. \ }

{\small [9]\ \ \ I. Cho, Graph von Neumann Algebras, ACTA Applied Math., 95,
(2007) 95 - 134.}

{\small [10]\ I. Cho, Characterization of Amalgamated Free Blocks of a Graph
von Neumann Algebra, Compl. An. Op. Theor., 1, (2007) 367 - 398.}

{\small [11] I. Cho, Measures on Graphs and Groupoid Measures, Comp. An. Op.
Theor., (2007) To Appear.}

{\small [12] I. Cho, Vertex-Compressed Algebras in a Graph von Neumann
Algebra, (2007), Submitted to ACTA Appled Math.}

{\small [13] I. Cho and P. E. T. Jorgensen, }$C^{*}${\small -Algebras
Generated by Partial Isometries, JAMC, (2007), To Appear.}

{\small [14] R. Exel, A new Look at the Crossed-Product of a }$C^{*}${\small %
-algebra by a Semigroup of Endomorphisms, (2005) Preprint.}

{\small [15] R. Gliman, V. Shpilrain and A. G. Myasnikov (editors),
Computational and Statistical Group Theory, Contemporary Math, 298, (2001)
AMS. }

{\small [16] R. Speicher, Combinatorial Theory of the Free Product with
Amalgamation and Operator-Valued Free Probability Theory, AMS Mem, Vol 132 ,
Num 627 , (1998).}
\end{quotation}

\begin{quote}
{\small [17] V. Vega, }$W^{*}$-{\small Algebras, Correspondences and Finite
Directed Graphs, (2007) Ph. D thesis, Univ. of Iowa.}

{\small [18] F. Radulescu, Random Matrices, Amalgamated Free Products and
Subfactors of the $C^{*}$- Algebra of a Free Group, of Noninteger Index,
Invent. Math., 115, (1994) 347 - 389.}

{\small [19]\ P. D. Mitchener, }$C^{*}${\small -Categories, Groupoid
Actions, Equivalent KK-Theory, and the Baum-Connes Conjecture,
arXiv:math.KT/0204291v1, (2005), Preprint.}

{\small [20]\ A. Gibbons and L. Novak, Hybrid Graph Theory and Network
Analysis, ISBN: 0-521-46117-0, (1999) Cambridge Univ. Press.}

{\small [21]\ F. Balacheff, Volume Entropy, Systole and Stable Norm on
Graphs, arXiv:math.MG/0411578v1, (2004) Preprint.}

{\small [22]\ I. Raeburn, Graph Algebras, CBMS no 3, AMS (2005).}

{\small [23] R. Scapellato and J. Lauri, Topics in Graph Automorphisms and
Reconstruction, London Math. Soc., Student Text 54, (2003) Cambridge Univ.
Press.}

{\small [24] S. H. Weintraub, Representation Theory of Finite Groups:
Algebra and Arithmetic, Grad. Studies in Math, vo. 59, (2003) AMS.}

{\small [25] W. Dicks and E. Ventura, The Group Fixed by a Family of
Injective Endomorphisms of a Free Group, Contemp. Math 195, AMS.}

{\small [26] N. I. Akhiezer and I. M. Glazman, Theory of Linear Operator in
Hilbert Space, vol II, Monograph \& Studies in Math, 10, ISBN 0-273-08496-8
(1981) Boston, Mass.}

{\small [27] L. A. Coburn, The }$C^{*}${\small -algebra Generated by an
Isometry, Bull. of AMS, 73, (1967) 722-726.}

{\small [28] \strut D.Voiculescu, K. Dykemma and A. Nica, Free Random
Variables, CRM Monograph Series Vol 1 (1992).\strut }

{\small [29] O. Bratteli, Inductive Limits of Finite Dimensional }$C^{*}$%
{\small -Algebras, Trans. AMS, 171, (1972) 195-234.}

{\small [30] O. Bratteli and D. W. Robinson, Operator Algebras and
Quantum-Statistical Mechanics, Texts \& Monographs in Physics, ISBN
0-387-10381-3, (1981) Springer-Verlag, NY.}

{\small [31] J. Cuntz, Simple }$C^{*}${\small -Algebras Generated by
Isometries, Comm. Math. Phys., 57, no 2, (1977) 173-185.}

{\small [32] M. Fannes, B Nachtergaele and R. F. Wener, Finitely Correlated
States on Quantum Spin Chains, Comm. Math. Phys., 144, no 3, (1992) 443-490.}

{\small [33] T. Kandelaki, On the Universal }$C^{*}${\small -Algebra
Generated by a Partial Isometry, Georgian Math. J., 5, no 4, (1998) 333-340.}

{\small [34] K. Kraus, Operations and Effects in the Hilbert Space
Formulation of Quantum Theory, Lect. Notes in Phys. vol 29, Springer,
Berlin, (1974) 206-229.}

{\small [35] A. V. Legedev and A. Odzievich, Extensions of }$C^{*}${\small %
-Algebras by Partial Isometries, Mat. Sb., 195, no 7, (2004) 37-70.}

{\small [36] J-M. Vallin, Multiplicative Partial Isometries and Finite
Quantum Groupoids, IRMA Lect. Math. Theor. Phys., vol 2, (2003) 189-227.}

{\small [37] K. Thomsen, Duality in Equivalent }$KK${\small -Theory, Pacific
J. Math., 222, no 2, (2005) 365-397.}

{\small [38] B. Sz.-Nagy and C. Foia, Harmonic Analysus of Operators on a
Hilbert Space, (1970) North-Holland Publ. Co., Amsterdam.}

{\small [39] P. E. T. Jorgensen, An Optimal Spectral Estimator for
Multi-Dimensional Time Series with an Infinite Number of Sample Points,
Math. Z., 183, (1983) 381 - 398.}

{\small [40] P. E. T. Jorgensen, Analytic Continuation of Local
Representations of Lie Groups, Pac. J. of Math., vol. 125, no. 2, (1986) 397
- 408.}

{\small [41] P. E. T. Jorgensen, Unitary Dialations and the }$C^{*}${\small %
-Algebra }$\mathcal{O}_{2}${\small , Israel J. of Math., vol. 56, no. 2,
(1986) 129 - 142.}

{\small [43] P. E. T. Jorgensen and P. Muhly, Self Adjoint Extensions
Satisfying the Weyl Operator Commutation Relations, J. D'anlyse Math., vol
37, (1980) 46 - 99.}

{\small [44] P. E. T. Jorgensen, Iterated Function Systems, Representations,
and Hilbert Spaces, Internat. J. of Math, vol. 15, no. 8, (2004) 813 - 832.}

{\small [45] R. Bott and R. Duffin, Impedance Synthesis without Use of
Transforms, J. of Appl. Phys., vol. 20, (1949) 816.}

{\small [46] J. Kigami, Analysis on Fractals, Cambridge Tracts in
Mathematics, vol. 143, ISBN: 0-521-79321-1, (2001) Cambridge Univ. Press.}

{\small [47] J. E. Hutchison, Fractals and Self-Similarity, Indiana Univ.
Math. J., vol. 30, (1981) 713 - 747.}

{\small [48] J. von Neumann, \r{U}ber Adjungierte Funktional Operatoren,
Ann. of Math., (2), vol. 30, no. 2, (1932) 294 - 310.}

{\small [49] R. T. Powers, Heisenberg Model and a Random Walk on the
Permutation Group, Lett. Math. Phys., vol. 1, no. 2, (1975) 125 - 130.}

{\small [50] R. T. Powers, Heisenberg Model, Comm. Math. Phys., vol. 51, no.
2, (1976) 151- 156.}

{\small [51] R. T. Powers, Resistence Inequalities for the Isotropic
Heisenberg Ferromagnet, J. of Math. Phys., vol. 17, no. 10, (1976) 1910 -
1918.}

{\small [52] R. T. Powers, Resistence Inequalities for the Isotropic
Heisenberg Model, }$C^{*}${\small -Algebras and Applications to Physics,
Lecture Notes in Math., vol. 650, Springer, (1978) 160 - 167.}

{\small [53] R. S. Strichartz, Differential Equations on Fractals, Princeton
Univ. Press, ISBN: 978-0-691-12731-6; 0-691-12731-X, (2006).}

{\small [54] L. A. Coburn, The }$C^{*}${\small -Algebra Generated by an
Isometry, Bull. of AMS, vol. 73, (1967) 722 - 726.}

{\small [55] F. K. Anderson and S. Silvestrov, The Mathematics of Internet
Search Engines, (2007) Preprint: Submitted to ACTA Appl. Math.}
\end{quote}

\end{document}